\documentclass[moor]{informs3}
\usepackage[dvipsnames]{xcolor}
\usepackage{array}
\usepackage{xy}
\usepackage{changepage}
\usepackage{epsfig}
\usepackage{setspace}
\usepackage{graphics}
\usepackage{wrapfig}
\usepackage{mathrsfs}
\usepackage{ifthen}
\usepackage{verbatim}
\usepackage{color}
\usepackage{mathrsfs}
\usepackage{bbm}
\usepackage{comment}
\usepackage{fancyvrb}
\usepackage{lmodern}
\usepackage{natbib}
 \def\bibfont{\small}%
 \bibpunct[, ]{[}{]}{,}{n}{}{,}
\usepackage[colorlinks=true,breaklinks=true,bookmarks=true,urlcolor=blue,
     citecolor=blue,linkcolor=blue,bookmarksopen=false,draft=false]{hyperref}

\usepackage{apptools}
\usepackage{bm}
\usepackage{hyperref}
\usepackage{thmtools}

\usepackage{setspace}

\usepackage{endnotes}
\let\footnote=\endnote

\EquationsNumberedBySection
\TheoremsNumberedBySection    % Preferred (Theorem 1, Lemma 1, Theorem 2)
\AtAppendix{\counterwithin{theorem}{section}}

\renewcommand{\P}{\mathbf{P}}

\newcommand{\B}{\mathcal{B}}
\newcommand{\V}{\mathcal{V}}
\newcommand{\D}{\mathbb{D}}
\newcommand{\BV}{\mathbb{BV}}
\newcommand{\AC}{\mathbb{AC}}
\newcommand{\CON}{\mathbb{CNCV}}

\newcommand{\BVT}{{\mathbb{BV}[0,T]}}
\newcommand{\BVM}{{\mathbb{BV}[0,M]}}

\newcommand{\E}{\mathbf{E}}
\newcommand{\I}{\mathbbm{1}}
\newcommand{\ind}{\mathbbm{1}}
\newcommand{\R}{\mathbb{R}}
\newcommand{\N}{\mathbb{N}}

\newcommand{\T}{\mathcal{T}}

\newcommand{\bigO}{\mathcal{O}}

\newcommand{\one}{\mathbbm{1}}

\newcommand{\cmt}[1]{{\color{red}[#1]}}

\newcommand{\eop}{\hfill\Halmos}
\renewcommand{\d}{\mathrm d}

\usepackage[normalem]{ulem}

\theoremstyle{plain}
\newtheorem{result}{Result}[section]

\newif\ifnavigationlinks
\newif\ifnotationindex
\newif\ifrebuttal
\newif\ifdiscussion
\newif\ifprintaddressedcomment

\newcommand{\opacity}{40}

\ifdiscussion
    \newcommand{\CR}[1]{ {\color{RedViolet}[CR: #1]} } % Comment from Chang-Han
    \newcommand{\MB}[1]{ {\color{RoyalBlue}[MB: #1]} } % Comment from Mihail
    \newcommand{\BZ}[1]{ {\color{OliveGreen}[BZ: #1]} } % Comment from Bert
    \newcommand{\JB}[1]{ {\color{RawSienna}[JB: #1]} } % Comment from Jose
    \ifprintaddressedcomment % Print addressed comments
        \newcommand{\MBA}[1]{ {\color{RoyalBlue!\opacity}[MB: #1]} } % Addressed Comment from Mihail
        \newcommand{\CRA}[1]{ {\color{RedViolet!\opacity}[CR: #1]} } % Adderssed Comment from Chang-Han
        \newcommand{\BZA}[1]{ {\color{OliveGreen!\opacity}[BZ: #1]} } % Addressed Comment from Bert
        \newcommand{\JBA}[1]{ {\color{RawSienna!\opacity}[JB: #1]} } % Addressed Comment from Jose
    \else % Do not print addressed comments
        \newcommand{\CRA}[1]{} 
        \newcommand{\MBA}[1]{} 
        \newcommand{\BZA}[1]{} 
        \newcommand{\JBA}[1]{} 
    \fi
\else
    \newcommand{\CR}[1]{}
    \newcommand{\MB}[1]{}
    \newcommand{\BZ}[1]{}
    \newcommand{\JB}[1]{}
    \newcommand{\CRA}[1]{} 
    \newcommand{\MBA}[1]{} 
    \newcommand{\BZA}[1]{} 
    \newcommand{\JBA}[1]{} 
\fi

\ifrebuttal
    \newcommand{\rvoutopacity}{100}
    \newcommand{\rvin}[1]{{\color{red}{#1}}}
    \newcommand{\rvout}[1]{{\color{red!\rvoutopacity}{\sout{#1}}} }
    \renewcommand{\rvout}[1]{{\color{black!\rvoutopacity}{\ifmmode\text{\sout{\ensuremath{#1}}}\else\sout{#1}\fi}} }
    \newcommand{\rvoutm}[1]{{\color{red!\rvoutopacity}{\sout{#1}}} }
    \renewcommand{\rvoutm}[1]{{\color{black!\rvoutopacity}{\ifmmode\text{\sout{\ensuremath{\displaystyle#1}}}\else\sout{#1}\fi}} }
\else
    \newcommand{\rvin}[1]{#1}
    \newcommand{\rvout}[1]{}
    \newcommand{\rvoutm}[1]{}
\fi

\newcommand{\rvtxt}[3]{\linkdest{#1}\hyperlink{back to #2}{\color{black}#3}%
}
\newcommand{\rvidx}[2]{\linkdest{link to #1}\hyperlink{#1}{#2}}
\newcommand{\A}[1]{{\color{NavyBlue}{#1}}\bigskip}

\usepackage[shortlabels]{enumitem}
\newlist{thmdependence}{itemize}{10}
\setlist[thmdependence]{nosep,label=-,leftmargin=1em}

\newcommand{\thmtreenode}[5]{\item[#1] \linkdest{location, thm tree #3} {#2}~\ref{#3} \linktopf{#3} \thmsum{#4}{#5}}

\newcommand{\thmtreeref}[2]{\item[\elsewhere] {{\hyperlink{location, thm tree #2}{\color{gray}#1}}}~\ref{#2}\thmsum{0.5}{}}

\ifnavigationlinks
\newcommand{\linksinthm}[1]{\emph{\linkdest{location, #1}\linktopf{#1} \linktothmtree{location, thm tree #1} }}
\newcommand{\linksinthmwopf}[1]{\emph{\linkdest{location, #1} \linktothmtree{location, thm tree #1} }}
\newcommand{\linksinpf}[1]{\linkdest{location, proof of #1}\linktothm{#1} \linktothmtree{location, thm tree #1} }
\else
\newcommand{\linksinthm}[1]{}
\newcommand{\linksinthmwopf}[1]{}
\newcommand{\linksinpf}[1]{}
\fi

\ifnotationindex
\newcommand{\notationdef}[2]{\linkdest{location, notation definition of #1}\hyperlink{location, notation index of #1}{#2}}
\newcommand{\notationidx}[2]{\linkdest{location, notation index of #1}\hyperlink{location, notation definition of #1}{#2}}
\newcommand{\linktonotation}[2]{\hyperlink{location, notation definition of #1}{#2}}
\else
\newcommand{\notationdef}[2]{#2}
\newcommand{\notationidx}[2]{}
\newcommand{\linktonotation}[2]{}
\fi

\newcommand{\notationindexlocationphrase}{location: notation index}
\newcommand{\linktopf}[1]{\hyperlink{location, proof of #1}{\pflinksymbol}}
\newcommand{\linktothm}[1]{\hyperlink{location, #1}{\thmlinksymbol}}
\newcommand{\linktothmtree}[1]{\hyperlink{#1}{\thmtreelinksymbol}}
\newcommand{\thmlinksymbol}{{\tiny [Theorem]}}
\newcommand{\pflinksymbol}{{\tiny [Proof]}}
\newcommand{\thmtreelinksymbol}{{\tiny [ThmTree]}}

\newcommand{\elsewhere}{}

\newcommand{\thmsum}[2]{\quad{\color{gray}\begin{minipage}[t]{#1\linewidth}{#2}\vspace{0.5\baselineskip}\end{minipage}}}

\makeatletter
\newcommand{\linkdest}[1]{\Hy@raisedlink{\hypertarget{#1}{}}}
\makeatother

\begin{document}

%%%%%%%%%%%%%%%%%%%%% \ifrebuttal %%%%%%%%%%%%%%%%%%%%%%%%%%%
\ifrebuttal
\newpage
% \bigskip
\pagestyle{plain}
\linkdest{location of rebuttal}
% \begin{center}
\hfil{\LARGE Second Round Revision Report}\hfil\\

% \medskip

\hfil{\large M. Bazhba, C.-H. Rhee, J. Blanchet, B. Zwart}\hfil\\
% \medskip
% \end{center}

\noindent
% \begin{minipage}[t]{\linewidth}
We thank the reviewer and the editors for carefully checking our paper again. 
We have addressed all the comments, and we attached our point-by-point responses to the reviewer's comments below. 
In each comment, we provided hyperlinks to the corresponding texts in the paper and color-coded the revised texts with a red font.
% : the added texts is in bright red, and the deleted texts are in light red. 
\BZA{I find this whole file too fancy. All this theoremtree and stuff needs removal when submitted. The rebuttal should be uploaded as separate file. Also, simply cross the deleted text. Shaded text can be confusing.}\CRA{I changed the macro so that the deleted text is crossed out and in gray font. This looks good to me, but I can easily customize the macro further so that the deleted texts look whatever way you want. The rebuttal and navigation links (theorem trees, toc, notation index, etc) can also be easily suppressed by simply commenting out lines 159--162 of the current tex source. When we submit the revised paper, I will comment out lines 159--162 to get the revised paper and comment out lines 159--161 to get the rebuttal). I attached such samples in the email I just sent out.}\BZA{ok great. I think we should leave the crossed-out text in the original color to make it visible on printed paper}
To facilitate the navigation, we also embedded a hyperlink in each revised text so that you can click the revised text to jump back to the corresponding item in the rebuttal. 
% \end{minipage}
\medskip

\noindent 
0) \linkdest{back to 0)} 
\begin{minipage}[t]{0.965\linewidth}
\emph{The structure of the paper has been completely changed and improved. The authors have taken into account my comments of the previous report. I recommend nevertheless that the authors read their paper again, and correct typos, mistakes or unproved claims in the new version of the paper.}\\

\A{We apologize for the mistakes that were still remaining in the submitted manuscript. We went through our manuscript again and corrected a good number of typos and mistakes. We also added extra explanations and references. 

\BZA{I have read on this moment [1130 AM tuesday dutch time] the first 3 sections. I will implement minor grammatical mistakes (i.e. adding punctiating, fixing trivial typos, other elementary grammar issues) without marking them. There are a couple of dozens of those. More significant issues that require attention and I didn't fix are:

1. we sometimes write "with (good) rate function" and sometimes "with THE (good) rate function". Similar for speed, or speed function or the speed or the speed function. I propose we will be consistent, and perhaps we should follow a convention like the one in  DZ. 

2. Definitions are not always introduced with the definition operator, especially in section 2 there are a bunch of those.

3. stylistic stuff: point-wise or pointwise? etc etc etc
}\CRA{Since DZ use the phrase ``... satisfies the LDP in $C([0,1])$ with the good rate function $I=...$,'' I added `the' to both `speed' and `rate function' whenever such phrase appears in our manuscript. I also checked pointwise/point-wise in DZ. They use pointwise. So I changed `point-wise' to `pointwise' in our manuscript. I also replaced the equalities to triangleq's whenever I saw definitions.}

% Please follow this \rvidx{(a)}{[link]} to see the revised text.
}
\end{minipage}

\noindent 
1) \linkdest{back to 1)}
\begin{minipage}[t]{0.965\linewidth}
\emph{In Theorem 2.3 the authors formulate the result concerning the asymptotic of the probability $\P_\pi(\cdot)$ (for stationary process), while in the proof of Theorem 2.1, they apply Theorem 2.3 for the probability $\P_0(\cdot)$ (for the process with the starting point $X_0 = 0$). And the proof of Theorem 2.3 is given also for the probability $\P_0(\cdot)$. These annoying typos come from the previous version of the paper where Theorem 2.3 was formulated and proved for the probability $\P_\pi(\cdot)$:}\\

\A{
We fixed these typos and rewrote the statement of Theorem~\ref{subexpbusyperiod-steadystate}.
Please follow this 
% [%
\rvidx{1-1}{[link]}
% ] 
to see the revised text.
}
\end{minipage}

\noindent 
2) \linkdest{back to 2)}
\begin{minipage}[t]{0.965\linewidth}
\emph{In Theorem 2.3 : For $b = 0, 1/(nb)\alpha$ is not defined. Maybe, the authors would like to say that for any $b \geq 0$,
$$\lim_{n\to\infty} \frac 1 {n^\alpha} \log \P(V_n \geq nb) = - \mathcal B_0^*b^\alpha.$$}

\A{
Thank you for pointing this out. We made changes as suggested. Please follow this \rvidx{2}{[link]} to see the revised text.
}
\end{minipage}

\noindent 
3) \linkdest{back to 3)}
\begin{minipage}[t]{0.965\linewidth}
\emph{page 10, line 8 and 9:
    If Theorem 2.3 holds for any $b \geq 0$, that is if (??) \BZA{are these question marks intentional?}\CRA{This is exactly how this comment appeared in the review report. My guess is that the reviewer must have put question marks as placeholder there, but never changed it to actual equation number before submitting the review. I suppose that the reviewer's intention was to refer to \eqref{formula-subexpbusyperiod-steadystate}.}\BZA{ok, lets keep it as is then}%
     holds for any $b \geq 0$, I do not see why the authors need to introduce $\epsilon > 0$.}    
    \\
    
\A{
We suppose that the reviewer refers to equation \eqref{formula-subexpbusyperiod-steadystate} by ``(??).'' 
If that's the case, we use $\epsilon$ to make sure that 
$$\left\{\frac1n V_n \geq 0\vee a + \epsilon\right\}\bigg\backslash\left\{\frac1n V_n \geq b\right\} \subseteq \left\{\frac1n V_n \in B\right\}$$
since $B= (a,b)\cap \R_+$ does not include $a$. 
We realize that this may have been confusing due to a typo in our previous statement that ${\P(\frac1n V_n \geq b)}/{\P(\frac1n V_n \geq 0\vee a)}\to 0$, where $+\epsilon$ was missing in the denominator. 
We fixed this typo and added one line to the math display to make the reason for introducing $\epsilon$ here clearer. 
Please follow this \rvidx{3}{[link]} to see the revised text.
}
\end{minipage}

\noindent 
4) \linkdest{back to 4)}
\begin{minipage}[t]{0.965\linewidth}
\emph{page 11, line 1: the authors should explain why $\tau$ is light-tailed (or give a reference to a result proving that $\tau$ is light-tailed).}
\\
    
\A{
We added a short explanation. We also added more details to the justification of \eqref{light-tailasy}. Please follow this \rvidx{4}{[link]} to see the revised text.
}
\end{minipage}

\noindent 
5) \linkdest{back to 5)}
\begin{minipage}[t]{0.965\linewidth}
\emph{page 37, Result B.1 : It seems that the same quantity is denoted by two
different notations: $\overline{\mathcal K}_n$ and $\overline{ K}_n^0$.}\\

\A{
Thank you for catching this. We fixed the typo. Please follow this \rvidx{5}{[link]} to see the revised text.
}
\end{minipage}

\noindent 
6) \linkdest{back to 6)}
\begin{minipage}[t]{0.965\linewidth}
\emph{I suggest also to simplify the notation by replacing $f(X_i)$ by $X_i^p$. Indeed, throughout the paper it is sometimes used as $f(X_i)$ and sometimes as $X_i^p$. I cannot see why the notation $f(X_i)$ is needed because in any case, the authors only consider the function $f(x) = x^p$.}\\

\A{Thank you for the suggestion. We changed $f(X_i^p)$ to $X_i^p$ throught the paper except for the introduction. 
Please follow these links [%
\rvidx{6-1}{1},% 
\rvidx{6-2}{2},%
\rvidx{6-3}{3},%
\rvidx{6-4}{4},%
\rvidx{6-5}{5},%
\rvidx{6-6}{6},%
\rvidx{6-7}{7},%
\rvidx{6-8}{8},%
\rvidx{6-9}{9}%
] to see the revised texts.
In the introduction, we compare our result to the standard Donsker-Varadhan type results in the literature, where $f$ cannot take the form $f(x) = x^p$ since $f$ has to be bounded, and we believe that the comparison is much more clear when we specifically refer to $f$ instead of calling it ``a (bounded) function.''
}
\end{minipage}

\newpage
\fi
%%%%%%%%%%%%%%%%%%%%% end of \ifrebuttal %%%%%%%%%%%%%%%%%%%%

\pagestyle{headings}
\setcounter{page}{1}

\RUNAUTHOR{M. Bazhba, C-H. Rhee, J. Blanchet, B. Zwart }

% Title or shortened title suitable for running heads. Sample:
% \RUNTITLE{Bundling Information Goods of Decreasing Value}
% Enter the (shortened) title:
\RUNTITLE{Sample-path LDP for unbounded additive functionals of the RRW}

% Full title. Sample:
% \TITLE{Bundling Information Goods of Decreasing Value}
% Enter the full title:
\TITLE{Sample-path large deviations for unbounded additive functionals of the reflected random walk}

% Block of authors and their affiliations starts here:
% NOTE: Authors with same affiliation, if the order of authors allows,
%   should be entered in ONE field, separated by a comma.
%   \EMAIL field can be repeated if more than one author
\ARTICLEAUTHORS{%
	\AUTHOR{Mihail Bazhba }
	\AFF{Quantitative Economics, University of Amsterdam, Netherlands, \EMAIL{M.Bazhba@uva.nl}} %, \URL{}}
	\AUTHOR{Jose Blanchet}
	\AFF{Management Science and Engineering, Stanford University, USA, \EMAIL{jose.blanchet@stanford.edu }}
	\AUTHOR{Chang-Han Rhee}
	\AFF{Industrial Engineeiring and Management Sciences, Northwestern University, USA, \EMAIL{chang-han.rhee@northwestern.edu}}
	\AUTHOR{Bert Zwart}
	\AFF{Stochastics Group, CWI and Eindhoven University of Technology, Netherlands, \EMAIL{bert.zwart@cwi.nl}}
	% Enter all authors
} % end of the block

\ABSTRACT{
	We prove a sample-path large deviation principle (LDP) with sub-linear speed for unbounded functionals of certain Markov chains induced by the Lindley recursion. The LDP holds in the Skorokhod space $\D[0,1]$ equipped with the $M_1'$ topology. Our technique hinges on a suitable decomposition of the Markov chain in terms of regeneration cycles. Each regeneration cycle denotes the area accumulated during the busy period of the reflected random walk. We prove a large deviation principle for the area under the busy period of the MRW, and we show that it exhibits a heavy-tailed behavior.  
}	

\maketitle
\thispagestyle{empty}

\section{Introduction}

In this paper we develop sample-path large deviation principles (LDP) for
additive functionals of a Markov chain which is important in Operations
Research (OR), namely, Lindley's recursion. This Markov chain describes the
waiting time sequence in a single-server queue under a FIFO discipline and
under independent and identically distributed (i.i.d.) inter-arrival times
and service times. We focus on the case in which the input is light-tailed, 
i.e. the service times and inter-arrival times have a finite moment
generating function in a neighborhood of the origin.

While the model that we consider is vital to many OR applications, and
therefore important in its own right, our main contributions are also fundamental from a methodological standpoint. We contribute, as we shall
explain, to the development of key tools in the study of sample-path large
deviations for additive functionals of light-tailed and geometrically 
ergodic Markov chains.

A rich body of theory, pioneered by Donsker and Varadhan in classical work
which goes back over forty years (see, for example, \cite
{donsker1975asymptotic}) provides powerful tools designed to study large
deviations for additive functionals of light-tailed and geometrically
ergodic Markov chains. Roughly speaking, these are chains which converge
exponentially fast to stationarity and whose stationary distribution is light-tailed.

%tail distribution decreases exponentially fast. % in the asymptotic regime. 

Unfortunately, despite remarkable
developments in the area, including the more recent contributions in \cite%
{kontoyiannis2003spectral}, the prevailing assumptions in the literature are
often not applicable to natural functionals of well-behaved geometrically
ergodic models, such as Lindley's recursion with light-tailed input.

In particular, every existing general result describing sample-path large
deviations of functionals of a process such as Lindley's recursion, must
assume the function of interest to be bounded. Hence, the current
state-of-the-art rules out very important cases, such as the
sample-path behavior of the empirical average of the waiting time sequence
in single-server queue over large time scales.
Our development allows one to study sample-path large
deviations for the cumulative waiting time sequence of a single-server queue.
In particular, we provide methodological ideas which, we believe, will
be useful in further development of the general theory of sample-path
large deviations for additive functionals of geometrically ergodic Markov
processes. More precisely, our contributions are summarized as follows,

%\smallskip

A) Let $\left\{ X_{k},k\geq 0\right\} $ follow  Lindley's recursion.
Assume that the associated increments have a finite moment generating
function in a neighborhood of the origin, \rvout{and}the traffic intensity is less
than one, 
% \rvtxt{6-13}{6)}{\rvout{
and let $f\left(x\right) =x^{p}$\ for any
% }\rvin{and}} 
$p>0$.
We establish a sample-path large deviation principle for $\bar{Y}_{n}\left( \cdot\right) =\sum_{k=1}^{\left\lfloor n\cdot \right\rfloor }
% \rvtxt{6-14}{6)}{\rvout{
f\left( X_{k}\right)
% }\rvin{X_i^p}}
/n$ as $n\rightarrow \infty $ w.r.t.\ the $M_{1}^{\prime }$ topology on %\todo{define BERT: I do not see what needs to be defined. $M_1'$ is already defined.} 
$\D[0,T]$ with a
good rate function and a sublinear speed function, all of which are fully
characterized in Theorem \ref{SPLDPYN}. Though our result only pertains to a specific Markov chain, they can be extended to more general stochastic recursions, and diffusions; this will be pursued
in future work. Related work, covering  the case of one-dimensional Langevin diffusions can be seen in \cite{LDadditiveFunctionalsDiffusion}.
%\smallskip

B) %As one of our contributions, we highlight 
We believe that our overall strategy for establishing
Theorem~\ref{SPLDPYN} can be applied generally to the
sample-path large deviations analysis of additive functionals of geometrically
ergodic Markov chains. Our strategy involves splitting the sample path in
cycles, roughly corresponding to returns to a compact set (in the case of
the Lindley recursion, the origin). Then, we show that the additive functional
in a cycle has a Weibullian tail. Finally, we use ideas similar to those
developed in \cite{bazhba2017sample}, involving sample-path large deviations for
random walks with Weibullian increments for the analysis. The result in \cite{bazhba2017sample}, however, cannot be applied directly to our setting here because of two reasons.
First, the cycle in progress at the end of the time interval is different from the rest. 
Second, the number of cycles (and thus the number of terms in the decomposition) is random.

%\bigskip

\smallskip

%C) Finally, we highlight a number of technical results of independent
%interest. For example, we derive a classical sample path large deviations
%principle for light-tailed random walks in the $J_{1}$ topology having a
%finite moment generating function in a neighborhood of the origin. The
%result, reported in Proposition~\ref{extensionmogulskildp}, extends the
%classical results in (i) \cite{mogul1976large}, which assumes that moment generating
%functions are finite everywhere, (ii) \cite{mogulskii1993large}, which focuses
%on L\'{e}vy processes, and (iii) \cite{puhalskii1995large}, which studies
%non-decreasing random walks in $\D[0,\infty)$. 

%We also develop an extension of the contraction
%principle, see Lemma~\ref{extendedexecontractionprinciple}.
%\bigskip

The sublinear speed of convergence highlighted in A) underscores the main
qualitative difference between our result and those traditionally obtained
in the Donsker-Varadhan setting. In our setting, as hinted in B), the
large deviations behavior of $\bar{Y}_{n}$ is characterized by heavy-tailed
phenomena (in the form of Weibullian tails) which arise when studying the
tails of the additive functional over a given busy period. Our choice of 
% \rvtxt{6-15}{6)}{\rvout{
$f\left( \cdot \right)$
% }\rvin{$f(X_i) = X_i^p$}} 
underscores the frailty of the boundedness assumptions required to apply the Donsker-Varadhan type theory.
Note that although $f(\cdot)$ grows slowly when $p \approx 0$, 
just a small amount of growth derails the application of the standard theory.

The choice of topology is an important aspect of our result. In {\cite{bazhba2017sample} it is argued that }$M_{1}^{\prime }$ is a natural topology to
consider for developing a full sample-path large deviation principle for random
walks with Weibullian increments. It is explained that such a result is
impossible in the context of the $J_{1}$ topology in $\D[0,T]$. Actually, to
be precise, the topology that we consider is a %slightly 
stronger variation
of the one considered by \cite{puhalskii1997functional} and \cite%
{puhalskii1998functional}, who introduced the $M_{1}^{\prime }$ topology in 
%\todo{define BERT: what? $M_1'$ is defined.} 
$\D[0,\infty )$, but in such a way that its direct projection onto $\D[0,T]$
loses important continuous functions (such as the maximum of the path in
the interval). The key aspect in our variation is the evaluation of the
metric at the right endpoint. 
%The differences between these two alternative interpretations of the $M_{1}^{\prime }$ topology on compact intervals disappears if the limiting functions are continuous, as those considered in \cite{puhalskii1997functional} and \cite{puhalskii1998functional}.
The version that we consider merges the
jumps, in the same way in which it is done at the left endpoint in the
standard $M_{1}^{\prime }$ description. This variation results in a stronger
topology when restricted to functions on compact intervals and it includes
the maximum as a continuous function. An important reason for using the
$M_{1}^{\prime }$ topology is that it allows merging jumps. This seems to
be particularly relevant given that in our setting the large deviations
behavior will eventually merge the increments within the busy periods.

In addition to the two elements mentioned in B), which make the result in \cite{bazhba2017sample}
not directly applicable, our choice of a strong topology also makes the approach in \cite{bazhba2017sample} 
difficult to use. In fact, in contrast to \cite{bazhba2017sample}, in this paper, we use a projective limit strategy to 
obtain our large deviation principle. A direct approach we explored, using the result in \cite{bazhba2017sample}, consisted
in replacing the random number of busy periods by its fluid limits (for which there is a
large deviations companion with a linear speed rate). Then, we tried to verify that this replacement
results in an exponentially good approximation. This would have been a successful strategy if we had 
used the version of the $M_{1}^{\prime }$ topology considered by \cite{puhalskii1997functional}, but 
unfortunately such exponential approximation does not hold in the presence of our stronger topology.

The development of Theorem \ref{SPLDPYN} highlights interesting and somewhat
surprising qualitative insights. For example, consider the case $f\left(
x\right) =x$, corresponding to the area drawn under the waiting time
as a curve. As we show, deviations of order $O(1)$ upwards from the typical behavior of the process $\bar{Y}_{n}\left( \cdot \right)$ occur due to extreme behavior in a single busy period of duration $O(n^{1/2})$. A somewhat surprising insight involves the busy period in process at time $n$, which is split into two parts of size $O(n^{1/2})$ involving the age and forward life time of the cycle (the former contributes to the area calculations, while the latter does not). This asymmetry, relative to the other busy periods during the time horizon $[0,n]$, which are completely accounted for inside the area calculation, raises the question of whether a correction in the LDP is needed, due to this effect, at the end of the time horizon. The answer is no; the contribution to the current busy period and the ones inside the time horizon are symmetric. This result is highlighted in Theorems~\ref{subexpbusyperiod-zero} and \ref{subexpbusyperiod-steadystate}, which characterize the variational problem 
governing extreme busy periods.

There are several related works that deal with large deviations for the area
under the waiting time sequence in a busy period. But they focus on queue
length as in \cite{blanchet2013large}, or assume that the moment generating
function of the increment is finite everywhere, as in \cite{duffy2014large}.
None of these works obtain  sample-path results. Instead, we do not
assume that the moment generating function of the service times or
inter-arrival times is finite everywhere. To handle this level of generality, we 
employ recently developed sampled-path LDP's \cite{borovkov2013large, borovkov2014large, vysotsky}.
%Because of our more general framework, we derive the extensions to the classical results mentioned in C).
This level of generality requires us to put in a substantial amount of work to rule 
out discontinuous solutions of the functional optimization problems that appear in the large deviations
analysis. %{\bf In addition, it is no longer possible to use direct methods from the calculus of variations as the local rate function may not be superlinear in this case (\cite{mogulskii1993large}). }

Another hurdle in developing tail asymptotics
for the additive functional in a busy period (reported in Theorems~\ref{subexpbusyperiod-zero} and \ref{subexpbusyperiod-steadystate}) is the fact
that the functional describing the area under the busy period is not
continuous. To deal with this, 
we exploit path properties of the most probable---in asymptotic sense--trajectories of the busy period along with the continuity of the area functional over a fixed time horizon. In particular, we rigorously show how to approximate the area over the busy period (which has a random endpoint) with the area over a large, {\em fixed} horizon. This is counter-intuitive at first, because the former approach allows one to remove the reflection operator. However, the latter approach does not have a first passage time (which is a discontinuous
function) as horizon, and this turns out to carry more weight.
%{\color{blue}In an earlier version of our proof, we attempted to exploit that the most likely path leading to a large area is concave, as the area functional is continuous at such paths. {\bf Though we obtain concavity as a by-product of our
%analysis (extending an idea from  \cite{duffy2010most}, who study a related problem), our method does not rely on this property, [our method relies on concavity; we use it to prove Lemma %\ref{B_y_star-equal-V_y_T_star}] and can therefore be extended to situations where the optimal path is not concave; we will exhibit this in more detail in a forthcoming
%study where we consider additive functionals of diffusion processes. 
%}} \cmt{BZ: their proof is not rigorous. Do we want to mention it? And we don't use concavity in lemma \ref{B_y_star-equal-V_y_T_star} but we prove it, and it is a non-essential byproduct, %which could be exploited for computational purposes}.  \\
%{\bf We exploit that the most likely path leading to a large area is concave, as the area functional is continuous at such paths.  Our method relies on this property and we use it to prove %that the optimal values of some key variational problems  are equal (Lemma~\ref{B_y_star-equal-V_y_T_star}).} \cmt{BZ: I didnt write this. I reiterate my comment above.}
%This paper is organized as follows. We give a detailed model description as
%well as the main results in Section~\ref{SMDMRPI}. %Section~\ref{Contnvar}
%presents  several
%structural properties of our key variational problem and its associated optimal solution. 
The proof of our sample-path LDP is provided in Section~\ref{samplepathproofs}. 
Section~\ref{sec:firstproofsection} focuses on the technical details behind deriving the tail asymptotics for the area under a busy period. 
%Section 5 explains how
%one may extend Theorems~\ref{subexpbusyperiod-zero} and \ref%
%{subexpbusyperiod-steadystate} (and therefore also Theorem~\ref{SPLDPYN}) to
%a larger class of Markov chains. 
The paper is closed with two appendices
covering auxilary duality results for Markov chains as well as large deviations
results.

%Sample path Large Deviation Principles (LDP) yield estimations for the probability tha`t
%the path of a random process $(Y)$ lies in a given set of paths. Moreover, when studying more general functionals of $Y$ it becomes natural to %consider logarithmic asymptotics, as is common in large
%deviations theory, cf.\ \cite{dembo2010large}.  From what we are aware of, sample path large deviation principles have been mainly developed for %processes {\color{red}...........}  \cite{mogulskii1993large}.

	\section{Model description and main results}\label{SMDMRPI}
	\subsection{Preliminaries}
	We consider the time-homogeneous Markov chain $\{X_n, n \geq 0 \}$ that is induced by the Lindley recursion, i.e.\, $ \notationdef{nota-Xk}{X_{n+1}} \triangleq [X_n+U_{n+1}]^{+}, \, n\geq 0,$
	and $X_0=0$. Note that the r.v.'s \notationdef{nota-Ui}{$\{U_i, i \geq 1\}$}  are i.i.d.\ copies of a random variable $U$ such that \notationdef{nota-mu}{$\mu \triangleq \E(U) <0$}. The state space of the Markov chain 
	$\{X_n, n \geq 0 \}$
	is the half-line of non-negative real numbers. We make the following technical but necessary assumptions:
\begin{assumption}\label{finitemgf}	
Let \notationdef{nota-theta-pm}{$\theta_+$ and $\theta_-$} be the supremum and infimum of the set $\{\theta: \E(e^{\theta U})< \infty \}$ respectively.
We assume that $-\infty \leq \theta_{-}<0<\theta_+\leq \infty$.
% \text{such that} \ \E(e^{\lambda U})< \infty \ \text{for every} \ \lambda \in (\theta_-,\theta_+),$$
\end{assumption}
\begin{assumption}\label{regularitycondition}
For $\theta_+$ and $\theta_-$ in Assumption~\ref{finitemgf},
$$\lim_{n \to \infty}\frac{\log \P(U \geq n)}{n} = -\theta_+, \lim_{n \to \infty}\frac{\log \P(-U \geq n)}{n} = \theta_-. $$
\end{assumption}
\begin{assumption}\label{obvious-condition}
$\P(U>0) > 0$.
\end{assumption}
    The purpose of this paper is to prove a sample-path LDP for $\bar Y_n = \{ \bar Y_n(t), t\in [0,T]\}$, where
	$$
	\notationdef{nota-Y-bar}{\bar Y_n(t)} \triangleq \frac{1}{n}\sum_{i=1}^{\lfloor n t \rfloor}\rvtxt{6-1}{6)}{\rvout{f(X_i)}\rvin{X_i^p}},\quad t \in [0,1], 
 	\rvtxt{6-1}{6)}{\rvout{\quad f(x)=x^p,}}
	$$
and 
$p > 0$ is a fixed constant.
We introduce basic notions that are used in the statement of one of our main results (Theorem~\ref{SPLDPYN}). 	
First, we set \notationdef{nota-alpha}{$\alpha \triangleq 1/(1+p)$}.
Let \notationdef{nota-D}{$\D[0,T]$} denote the Skorokhod space---the space of càdlàg paths from $[0,T]$ to $\R$. 
We sometimes also consider the space $\D[0,\infty)$ of càdlàg paths from $[0,\infty)$ to $\R$. 
%When the domain is omitted, we always refer to $[0,T]$. 

Let \notationdef{nota-TM1'}{$\mathcal{T}_{M_1'}$} denote the $M_1'$ Skorokhod topology, whose precise definition is provided below.  
Unless specified otherwise, we  assume that $\D[0,T]$ is equipped with $\mathcal T_{M_1'}$ throughout the rest of this paper.
	\begin{definition}
		For $\xi\in \D[0,T]$, define the extended completed graph \notationdef{nota-Gamma'}{$\Gamma'(\xi)$} of $\xi$ as
		$$
		\Gamma'(\xi) \triangleq \{(u,t)\in \R\times [0,T]: u\in [\xi(t-)\wedge \xi(t),\ \xi(t-)\vee \xi(t)]\},
		$$
		where $\xi(0-) \triangleq 0$.
		Define an order on the graph $\Gamma'(\xi)$ by setting $(u_1,t_1) < (u_2,t_2)$ if either
			 $t_1 < t_2$, or
			 $t_1 = t_2$ and $|\xi(t_1-)-u_1| < |\xi(t_2-) - u_2|$.
		We call a continuous nondecreasing function $(u,t) =\big((u(s), t(s)),\,s \in [0,T]\big)$ from $[0,T]$ to $\Gamma'(\xi)$ a parametrization of $\Gamma'(\xi)$ if $\Gamma'(\xi) = \{(u(s), t(s)): s\in [0,T]\}$. 
		We also call such $(u,t)$ a parametrization of $\xi$, and we denote the set of all parametrizations of $\xi$ with \notationdef{nota-Pi'}{$\Pi'(\xi)$}.
		%We call a pair of continuous functions $(u,t) \triangleq ((u(s), t(s)), s \in [0,1])$ such that $t(s)$ is nondecreasing with $t(0)= 0$ and $t(1)=1$ a parametrization of $\Gamma'(\xi)$ if $\Gamma'(\xi) = \{u(s), t(s)): s\in [0,1]\}$.
	\end{definition}
	
	\begin{definition}
		Define the $M_1'$ metric \notationdef{nota-dM1'}{$d_{M_1'}$} on $\mathbb D$ as follows
		$$
		d_{M_1'}(\xi,\zeta)
		\triangleq
		\inf_{\substack{(u,t) \in \Pi'(\xi)\\(v,r) \in \Pi'(\zeta)}}
		\{
		\|u-v\|_\infty + \|t-r\|_\infty
		\}.
		$$
	\end{definition}
		We say that $\xi\in \D[0,T]$ is a pure jump path if $\xi = \sum_{i=1}^\infty x_i\I_{[u_i,T]}$  such that $x_i\in \R$ and $u_i\in[0,T]$ for each $i \geq 1$ and the $u_i$'s are all distinct. Let \notationdef{nota-D-p-uparrow}{$\D_{\leqslant \infty}^\uparrow[0,T]$} be the subspace of $\D[0,T]$ consisting of non-decreasing pure jump paths that assume non-negative values at the origin.
%	Let $\mu_1 =\E(W_1).$
%
%Furthermore, define
%$\D_{\leqslant 1}(1)=\{\xi \in \D: \xi \in \D_{\leqslant 1}, \xi = x_1 \cdot \mathbbm{1}_{\{1\}}, x_1 \in [0,\infty)\}$ i.e; the subspace of step functions with one jump at the end of the time horizon.
%Let $$\mathbbm{C}_1^{\mu_1} \triangleq \{\zeta \in \D[0,1]: \zeta =\mu_1 \cdot t +\xi_1+ \xi_2, \xi_1 \in \D_p^{\uparrow}, \xi_2 \in \D_{\leqslant 1}(1)\}.$$
Let \notationdef{nota-BV-T}{$\BV[0,T]$} be the subspace of $\D[0,T]$ consisting of paths with finite variation. Every $\xi \in \BV[0,T]$ has a Lebesgue decomposition with respect to the Lebesgue measure. 
That is, $\xi=\xi^{(a)}+\xi^{(s)}$, where \notationdef{nota-xia}{$\xi^{(a)}$} denotes the absolutely continuous part of $\xi$, and \notationdef{nota-xis}{$\xi^{(s)}$} denotes the singular part of $\xi$. 
Subsequently, using Hahn's decomposition theorem we can decompose $\xi^{(s)}$ into its non-decreasing singular part \notationdef{nota-xiu}{$\xi^{(u)}$} and non-increasing singular part \notationdef{nota-xid}{$\xi^{(d)}$} so that $\xi^{(s)} = \xi^{(u)}+\xi^{(d)}$. 
Without loss of generality (w.l.o.g.), we assume that $\xi^{(s)}(0) = \xi^{(u)}(0) = \xi^{(d)}(0) = 0$.
We will sometimes also consider \notationdef{nota-BV-infty}{$\BV[0,\infty)$} - the subspace of $\D[0,\infty)$ consisting of paths that are of bounded variation on any compact interval.

\subsection{Sample-path large deviations}
In this subsection, we present the sample-path large deviation principle for $\bar{Y}_n$ and the main ideas of its proof. 
We start with a few definitions. 
Let $\Psi$ be the reflection map defined by $\notationdef{nota-Psi}{\Psi(\xi)(t)}\triangleq\xi(t)-\inf_{s\in[0,t]}\{\xi(s) \wedge 0\}, \ \forall  t \geq 0$.
Define $\notationdef{nota-cal-T}{\T(\xi)} \triangleq \inf\{t>0: \Psi(\xi)(t) = 0\}$, 
$\notationdef{nota-By}{B_y} \triangleq \left\{ \xi \in \BV[0,\infty): \xi(0)=y, \ \int_{0}^{\T(\xi)}\Psi(\xi)(s)^pds \geq 1 \right\}$, 
and $\notationdef{nota-Lambda-star}{\Lambda^*(y)}\triangleq \sup_{\theta \in \R}\{\theta y -\log \E(e^{\theta U})\}.$ 
Set
$$
\notationdef{nota-Iy}{I_{y} (\xi)}
\triangleq
\begin{cases}
\int_{0}^{ \T(\xi)}  \Lambda^*(\dot{\xi}^{(a)}(s))ds + \theta_+\xi^{(u)}(\T(\xi)) + \theta_-\xi^{(d)}(\T(\xi))
& \text{if } \xi(0) = y \text{ and }\xi \in \BV[0,\infty),
\\
\infty
& \text{otherwise},
\end{cases}
$$
and denote with \notationdef{nota-cal-B-y-star}{${\mathcal{B}_{y}^{\ast}}$} the optimal value of the  variational problem
\[
\B_y^* \triangleq \inf_{\xi \in B_{y}} I_y(\xi).
\notationdef{nota-cal-B-y}{\tag{$\mathcal{B}_{y}$} \label{By}}
\]
Similarly, denote with \notationdef{nota-cal-B-y-pi}{${\mathcal{B}_{\pi}^{\ast}}$}  the optimal value of the variational problem
\begin{equation}
\B_\pi^* \triangleq \inf_{y \in [0,\infty),\, \xi \in B_{y}} \left\{ \beta y+ I_y(\xi) \right\}, 
\tag{$\mathcal{B}_{\pi}$} \label{Bpi}
\end{equation}
where $\notationdef{nota-beta}{\beta} \triangleq \sup\left\{\theta \geq 0: \E(e^{\theta  U}) \leq 1 \right\}$ is the decay rate of the steady state distribution  $\pi$ of the reflected random walk (see Result~\ref{steadystateasymptotics}). Note that $\beta \leq \theta_+$  and $\beta$ is strictly positive in view of Assumption 1 and the assumption that $\mu<0$.
Note also that $\mathcal B_\pi^* = \inf_{y\in [0,\infty)} \left\{\beta y + \mathcal B_y^*\right\}$.

Let $T_0 \triangleq 0$ and $ \notationdef{nota-T-i}{T_i} \triangleq \inf\{ k> T_{i-1}: X_k = 0\}$ for $i\geq 1$, and subsequently, define $\notationdef{nota-lambda}{\lambda} \triangleq \E(\sum_{i=1}^{T_1}X_i^p)/\E(T_1)$. Define $\notationdef{nota-D-lambda}{\D^{(\lambda)}[0,T]}\triangleq \{\xi \in \D[0,T]: \xi(t) =\lambda t + \zeta(t),\,\forall t\in[0,T],\, \zeta \in \D_{\leqslant \infty}^{\uparrow}[0,T]\}$, i.e., the subspace of increasing paths with slope $\lambda$ and countable upward jumps.
Recall that $\alpha = 1/(1+p)$.%
\begin{theorem}\label{SPLDPYN}
\linksinthm{SPLDPYN}
	The stochastic process $\bar{Y}_{n}$ satisfies a large deviation principle in $(\D[0,T], \mathcal{T}_{M_1'})$ with the speed $n^{\alpha}$ and the rate function $I_{\textbf{Y}}: \D[0,T] \to \R_+$, defined as
	\begin{equation}
	\notationdef{nota-I-Y}{I_{\textbf{Y}}(\zeta)}
	\triangleq
	\begin{cases}
    	{\mathcal{B}^*_{0}} \sum_{t: \zeta(t) \neq \zeta(t-)}\left(\zeta(t)-\zeta(t-)\right)^{\alpha}
    	&
    	\text{if } \zeta \in  \D^{(\lambda)}[0,T], 
    	\\
    	\infty
    	&
    	\text{otherwise}.
	\end{cases}
	\end{equation}
    That is, for any measurable set $A$,
    \begin{equation}
    -\inf_{A^{\circ}}I_{\textbf{Y}}(\xi)
    \leq
    \liminf_{n \to \infty}
    \frac{\log\P\left( \bar{Y}_n \in A\right)}{n^{\alpha}}
    \leq
    \limsup_{n \to \infty}
    \frac{\log\P\left( \bar{Y}_n \in A\right)}{n^{\alpha}}
    \leq
    -\inf_{\bar{A}}I_{\textbf{Y}}(\xi).
    \end{equation}	
%where the interior and closure in the infimum are w.r.t.\ $\mathcal T_{M_1'}$.
\end{theorem}
%{\color{red}From the form of the rate function, taking into account a possible jump at the end-point of $[0,1]$ it is clear that the usage of the Skorokhod $J_1$ topology is infeasible. In fact, it may be possible that $ {\mathcal{B}_{\pi}^{\ast}}^{(p)} < {\mathcal{B}_{0}^{\ast}}^{(p)} $ so that this value will dominate.
%Thus, the heuristics from the theory of heavy tails for i.i.d. Weibullian random variables, in which a large sum occurs by a single
%large value taking place uniformly across the interval, may break down in this case.}
%	
The full proof of Theorem~\ref{SPLDPYN} is deferred to Section~\ref{samplepathproofs}. The strategy relies on a suitable  representation for $\bar Y_n$  using renewal theory, which is presented next. 
\iffalse
% BERT: WHAT FOLLOWS IS SO STANDARD THAT IT CAN BE OMITTED, YOU CAN PUT IT BACK IN FOR THE THESIS
Consequently, the Markov chain $X_{k}$ is regenerative with respect to the sequence $\{T_j: j \geq 0\}$. In particular,
		\begin{itemize}
		\item $\{ X_{T_n}, \ldots, X_{T_{n+1}-1}\}_{n \geq 1}$ are i.i.d.\ and
		\item $\{ X_{T_n}, \ldots, X_{T_{n+1}-1}\}_{n \geq 1}$ are independent of $X_0,\ldots, X_{T_{n}-1}$.
	\end{itemize}
\fi
	The sequence $\{T_j, j\geq 1\}$ induces a renewal process $\{N(t), t\geq 0\}$ defined by
	$
	\notationdef{nota-N-t}{N(t)} \triangleq \max\{k \geq 0: T_k \leq t\}, \ t\geq 0.
	$
	We decompose the process $\bar Y_n$ as follows. For fixed $t\geq 0$:
    \begin{equation}\label{eq:representationYbarandom}
	\bar{Y}_{n}(t) = \frac 1 n\sum_{j=1}^{N(nt)}\sum_{i=T_{j-1}+1}^{T_j} \rvtxt{6-2}{6)}{\rvout{f(X_i)}\rvin{X_i^p}}+ \frac 1 n\sum_{i=T_{N(nt)}+1}^{\lfloor n t \rfloor}\rvtxt{6-2}{6)}{\rvout{f(X_i)}\rvin{X_i^p}},
	\end{equation}
	with the convention that $\sum_{i=T_{N(nt)}}^{\lfloor n t \rfloor}\rvtxt{6-3}{6)}{\rvout{f(X_i)}\rvin{X_i^p}}$   is zero in case the superscript $\lfloor n t \rfloor$ is
	strictly smaller than the subscript $T_{N(nt)}$.
	We introduce some notation for the analysis of $\bar Y_n$. Define
		\begin{itemize}
		\item  $\notationdef{nota-tau-j}{\tau_j} \triangleq T_{j}-T_{j-1}, j \geq 1$, the inter-arrival times of the renewal process $N$,
		\item  $\notationdef{nota-W-j}{W_{j}} \triangleq \sum_{i=T_{j-1}+1}^{T_j}\rvtxt{6-4}{6)}{\rvout{f(X_i)}\rvin{X_i^p}}, j\geq 1$, the ``area'' under $\rvtxt{6-4}{6)}{\rvout{f(X_i)}\rvin{X_i^p}}$ during the $j$\textsuperscript{th} busy period of $X_n$,
		\item  $\notationdef{nota-bar-Z-n}{\bar{Z}_n(t)} \triangleq \frac 1 n\sum_{j=1}^{N(nt)}W_j, t\in [0,1]$, the aggregate process (excluding the last regeneration cycle),
		\item  $\notationdef{nota-bar-cal-R-n}{\bar{R}_n(t)} \triangleq \frac{1}{n}\sum_{i=T_{N(n)}+1}^{\lfloor nt \rfloor}\rvtxt{6-5}{6)}{\rvout{f(X_i)}\rvin{X_i^p}}, t\in [0,1]$, the process that equals $\bar Y_n$ during the last regeneration cycle,
		\item  $\notationdef{nota-V-n}{V_n} \triangleq \sum_{i=T_{N(n)}+1}^{n}\rvtxt{6-6}{6)}{\rvout{f(X_i)}\rvin{X_i^p}}$, the ``area'' under $\rvtxt{6-6}{6)}{\rvout{f(X_i)}\rvin{X_i^p}}$ during the last regeneration cycle,
		\item  $\notationdef{nota-bar-S-n}{\bar{S}_n(t)} \triangleq \frac1n V_n\mathbbm{1}_{\{1\}}(t)$, the process with one jump, which aggregates the ``area'' under $\rvtxt{6-7}{6)}{\rvout{f(X_i)}\rvin{X_i^p}}$ during the last regeneration cycle.
	\end{itemize}

	We often refer to the functions associated with these quantities by dropping the argument; for example by writing $\bar{R}_n$ we refer to $\{\bar{R}_n(t), t\in [0,T]\}$.
	The strategy to prove our main result, Theorem~\ref{SPLDPYN},  %Section \ref{samplepathproofs}, 
	 builds on tail estimates for $W_1$ and $V_n$,  which are presented in Theorems~\ref{subexpbusyperiod-zero} and \ref{subexpbusyperiod-steadystate} below.
%		\  $\bar{Z}_n$ and $\bar{S}_n$ satisfy an LDP in $(\D[0,T], \mathcal{T}_{M_1'})$. 
%		\item $\bar{Z}_n+\bar{S}_n $ satisfies an LDP in $(\D[0,T],\mathcal{T}_{M_1'})$ with the rate function $I_{\mathbf{Y}}$.
%		\item $\bar{Z}_n+\bar{S}_n $ and $\bar{Y}_n$ are exponentially equivalent in the $M_1'$ topology.
%	\end{itemize}
%	\cmt{BZ: this needs to be restructured, rewritten. I will do it tomorrow.}
%	Regarding the first step, we first discuss the logarithmic asymptotics of $\bar V_n$ and $W_1$, which are presented in 
 %	For the sample path LDP of $\bar{Y}_n$, we prove the exponential equivalence of $\bar{\mathcal{Y}}_n$ and $\bar{\mathcal{S}}_n$ in Lemma~\ref{exp-equiv-ldp-sn} by pushing the last cycle $\bar{R}_n$ to the end of the time horizon. Consequently, the LDP of $\bar{R}_n$ is deduced, due to the LDP of $\bar{S}_n$ in $(\D[0,1],\mathcal{T}_{M_1'})$.
 			Using Theorem \ref{subexpbusyperiod-steadystate}, we derive an LDP for $\bar{S}_n$ in Lemma~\ref{exp-equiv-ldp-sn}. For the LDP of $\bar{Z}_n$ in $\D[0,T]$, we start by obtaining an LDP for the finite projections of $\bar{Z}_n$ in Lemmas~\ref{tails-asy-for-k-dim},~\ref{ldp-for-non-overlapping-partial-sums},~and \ref{fin-dim-LDP-process}. Then, the finite-dimensional LDP is lifted in the standard projective limit to the pointwise convergence topology in Lemma~\ref{LDP-Z-n-pointwise-convergence}, and finally extended to the $M_1'$ topology using the continuity of the identity map in the subspace of increasing c\'adl\'ag paths in Lemma~\ref{lemma-ldp-Z}. In the last step, we infer an LDP for  $\bar{Z}_n+\bar{S}_n $ through the use of a continuous mapping approach, and hence, we obtain the LDP for $\bar{Y}_n$.
 	% 	We derive an LDP for $\bar{Z}_n$ in $\D$ with respect to the $M_1'$ topology by obtaining an LDP with the point-wise convergence topology which is strengthened to the $M_1'$ topology using the continuity of the identity map in the subspace of increasing c\'adl\'ag paths. In the last step, we infer an LDP for  $\bar{Z}_n+\bar{S}_n $ through the use of a continuous mapping approach, and hence, we obtain the LDP for $\bar{Y}_n$.
		For the sample-path LDP of $\bar{Y}_n$, we prove the exponential equivalence of $\bar{{Y}}_n$ and $\bar{Z}_n+\bar{S}_n$ in Lemma~\ref{exponentialequivalence-Y-R-Z} by considering the $M_1'$ distance of $\bar{Y}_n$ with the aggregate process $\bar{Z}_n$ and the last regeneration cycle $\bar{R}_n$ pushed  to the end of the time horizon. Consequently, the LDP of $\bar{{Y}}_n$ is deduced, due to the LDPs of $\bar{S}_n$ and $\bar{Z}_n$ in $(\D[0,T],\mathcal{T}_{M_1'})$.

	Before embarking on the execution of this technical program, it is worth commenting on the role of $\bar{R}_n$, since this element will allow us to expose the importance of a careful analysis involving the area during a busy period. 
	As mentioned in the introduction, one may wonder if the contribution of $\bar{R}_n$  may end up counting different in the form of the LDP. 
	The typical path for $\bar{Y}_{n}$ is a straight line with drift equal to the steady-state workload. 
	Our developments indicate that most likely large deviations behavior away from the most likely path occur due to isolated busy periods which exhibit extreme behavior. 
	For example, in the case \rvtxt{6-8}{6)}{\rvout{$f(x)=x$}\rvin{$p=1$}}, substantially extreme busy periods (leading to large deviations of order $O(n)$) have a duration of  $O\left( n^{1/2}\right)$ and exhibit excursions of order $O\left(n^{1/2}\right)$, therefore accumulating an area of order $O\left( n\right) $. 
	
	%Our results in the next section characterize the variational problem which governs such extreme busy periods. But
	Each busy period, including the one in progress at the end of the time horizon contributes the same way in the rate function. This follows from Theorems~\ref{subexpbusyperiod-zero} and \ref{subexpbusyperiod-steadystate} given below, but may be somewhat remarkable. The reason is that when the cycle in progress at the end of the time horizon is extreme, its duration is of order $O\left( n^{1/2}\right)$. This suggests that the remainder of the cycle is also of order $O\left( n^{1/2}\right) $. It turns out that this long time duration has no significant contribution to the total area: while the remainder of the cycle in progress may be large, the position of the chain is actually of order $o\left( n^{1/2}\right)$ from the end of the time horizon, so the total contribution to the area of the remaining portion of the cycle is negligible. This calculation is exposed in Proposition \ref{miscellaneous-proposition-for-tail-asymptotics} below, and a time-reversal argument given in Appendix~\ref{genresults}.
	
%	[From Jose: the following two paragraphs need to be used to emphasize intuition, either here or in the next section. I need to streamline to avoid repetition. Currently, emphasizing this type of intuition would provide good connective tissue with the rest of the description.] 
	
%	That is, since $\{W_i\}_{i\geq 1}$ have heavy-tailed Weibull asymptotics, we  deduce that $\bar{Z}_n$ satisfies the same large deviation principle as with the random walk with heavy tailed Weibull increments.

%BERT: Jose, many thanks. I am ok with it as it is now. 

 \subsection{Busy period asymptotics}\label{BPnThnPrnId}

 It is clear that a large deviations analysis of the area under a busy period is indispensable for deriving the sample-path LDP of $\bar Y_n$ in Theorem~\ref{SPLDPYN}. Our next two theorems provide the asymptotic estimation for the tails of $W_1$ and $V_n$, showing that they exhibit Weibull behavior. We will discuss their statements and defer their proofs to Section \ref{sec:firstproofsection}.
\begin{theorem}
\label{subexpbusyperiod-zero}
\linksinthm{subexpbusyperiod-zero}
	Recall that $W_1 = \sum_{k=1}^{T_1} X_k^p$ and $\alpha = 1/(1+p)$. It holds that
	\begin{equation}\label{W-1-asymptotics}
	\lim_{t \to \infty}\frac{1}{ t^{\alpha}}\log\P\left(W_1 \geq t\right) = -{\mathcal{B}_0^{\ast}}.
	\end{equation}
\end{theorem}

For $V_n$, our analysis points to  Weibull-like asymptotic behavior similar to $W_1$ except that the pre-factor associated with $V_n$ is ${\mathcal{B}^{\ast}_{\pi}}$ (instead of $\mathcal B_0^*$). It turns out that (see Proposition~\ref{miscellaneous-proposition-for-tail-asymptotics}) the pre-factor $\mathcal{B}_{\pi}^*$ is equal to $\mathcal{B}_0^*$.  This leads to the conclusion that every busy period, including the one in progress at the end of the time horizon, has the same tail asymptotics.  	
% \begin{remark}\label{nojumpatzero}
% 	Since $\theta_+ \geq \beta$, Proposition~\ref{bp-ison-bzero} implies that the optimal solution of ${\mathcal{B}_y}, y>0,$ has no discontinuity at time zero.
% \end{remark}
\begin{theorem}
\label{subexpbusyperiod-steadystate}
\linksinthm{subexpbusyperiod-steadystate} 
For the area of the \rvtxt{1-1}{1)}{\rvin{last}} busy period\rvtxt{1-1}{1)}{\rvout{ starting from steady state (with distribution $\pi$)}}, \BZA{ADD: active at time $n$, }\CRA{I changed `For the area of the busy period' to `For the area of the last busy period'}%
we have \rvtxt{1-1}{1)}{\rvout{that}\rvin{the following tail asymptotics:}} for any $b\geq 0$,
\begin{align}\label{formula-subexpbusyperiod-steadystate}
    &
    \lim_{n \to\infty} \frac{1}{\rvtxt{2}{2)}{\rvout{(nb)^\alpha}\rvin{n^{\alpha}}}} \log \rvtxt{1-1}{1)}{\rvout{\P_{\pi}}\rvin{\P}} (V_n \geq nb) = -{\mathcal{B}^{\ast}_{0}}\rvtxt{2}{2)}{\rvin{\cdot b^{\alpha}}}.
    %  \\
    %  &
    %  \lim_{t\to\infty} \frac{1}{t^{\alpha}} \log \P_{0} (Y_0(p)>t) = -\mathcal{V}^*_{0}(p).
    \end{align}
\end{theorem}
The tail asymptotics for $W_1$ and $V_n$ are derived using a recently developed  LDP for random walks with light-tailed increments due to
\cite{borovkov2013large, borovkov2014large, vysotsky}, cf.\ Result~\ref{extensionmogulskildp} below. %and several versions (extensions) of the contraction principle.
Specifically, the tail probability of $W_1$ can be written as the tail probability of the image $\Phi(\bar K_n^0)$ of the unrestricted scaled random walk $\{\bar K_n^0(t), t\in [0,T]\}$, where $\notationdef{nota-bar-K-n-y}{\bar K_n^y (t)}\triangleq y + \frac{1}{n}K_{\lfloor nt \rfloor}$,  $\notationdef{nota-K-n}{K_n} \triangleq\sum_{i=1}^n U_i$, and the functional $\Phi$ is defined as $\notationdef{nota-Phi}{\Phi(\xi)}\triangleq \int_{0}^{\mathcal{T}(\xi)}(\Psi(\xi)(s))^p ds$.
Note that $\Phi:\D[0,\infty) \to \R_+$ is not continuous, and hence, the  proof for the tail asymptotics of $W_1$ gets more involved than simply applying the contraction principle. 
We derive large deviations upper and lower bounds and show that they coincide. 

For the upper bound, we replace the hitting time $T_1$ with a sufficiently large value $T$. 
\BZA{$T$ is used in the introduction, and in the first lemma of section 3... do we want to modify or we leave it?}\MBA{Thanks for the comment. I suggest to leave it as it is. I do not believe the reader will be confused by it and it is good keep notation at a minimum.}\CRA{I'd ideally change $T_i$'s to something else, but I vote for leaving as is, since there are too many of them throughout the manuscript, and we don't have many symbols left. I changed the function $T$ to \hyperlink{TtoUpsilon}{$\Upsilon$}, though.}%
This enables us to  study  the area  of $X_n$  over the finite time horizon $[0,T]$. 
For $T$ large enough, we show that the area of the reflected random walk over the whole time horizon $[0,T]$ serves as an asymptotic upper bound  for $W_1$ (see the proof of Proposition \ref{technical-proposition-for-tail-asymptotics}), and it is expressed as a  functional of $\bar K_n$ (see Lemma \ref{intermediate-lemma}). This functional is shown to be uniformly continuous in the $M_1'$ topology on level sets of the rate function associated with the LDP for $\bar K_n$ (cf.\ Lemma \ref{uniformcontinuity}). Invoking Result \ref{extensionmogulskildp}, recently established in \cite{vysotsky}, we get a large deviation upper bound.

For the lower bound, we confine the functional of the area under the busy period, over a fixed time horizon by imposing an extra condition, see the proof of Proposition \ref{technical-proposition-for-tail-asymptotics}. Subsequently, we exploit some regularity properties of a variational problem associated with the lower bound to show that $\mathcal{B}_0$ has the same optimal value as the variational problem associated with the large deviation upper and lower bounds. We have organized the presentation in such a way that analytic details associated to variational problems are gathered in Section \ref{additionalproofsection}.
%We show the equivalence of the aforementioned variational problems by using a structural property of their optimal paths (concavity) which we prove through the use of functional optimization techniques.  

For $V_n$, we follow the same approach with some slight modifications. In order to carry out our analysis  for $V_n$, we associate the tail of $V_n$  with the tail of  $W_1$ through Lemmas~\ref{timereversibility} and \ref{lastcycleequivalence}.
	In particular, we prove that %$\bar V_n$ has similar tail asymptotics to that of $W_1$, initialized from the steady state of $X_n$ i.e;
	$$
	\lim_{n\to \infty}\frac{\log\P_0(V_n > nx)}{n^{\alpha}}
	 =
	  \lim_{n\to \infty}\frac{\log\P_{\pi}(W_1 > nx)}{n^{\alpha}}.
	$$
To show this, we also rely on Result~\ref{steadystateasymptotics} from \cite{nuyenszwart}, describing the asymptotic behavior of the invariant measure $\pi$. Lastly, we repeat similar steps as in the analysis of $W_1$.

The constant $\mathcal{B}_0^*$ appears in all our theorems, and is the solution of a variational problem. We show in Proposition \ref{eq:ineqformonotonicity} that $\mathcal{B}_0^* \in (0,\infty)$. 
This property is all that is needed for our main sample-path large deviations results. 
Nevertheless, it is of interest to compute $\mathcal{B}_0^*$. This can be done by solving a suitable variational problem which in turn is typically done using the associated Euler-Lagrange equations. 
However, for the Euler-Lagrange equations to characterize the solution, it must be shown that an optimizing path $\xi^*$ exist which is sufficiently smooth, i.e.\ not just absolutely continuous, but differentiable with continuous derivative.  
In general, showing a priori sufficient smoothness of an optimizer is a non-trivial task, but we explain how to execute this for the case $p=1$, using a framework presented in \cite{Cesari}. 
The details can be found in Appendix \ref{appendix-C}.

\section{Proof of the sample-path LDP}
%Proofs for the sample path LDPs of $\bar{Y}_n-\bar{R}_n$ and %$\bar{R}_n$}
\label{samplepathproofs}

In this section, we prove  Theorem~\ref{SPLDPYN}. For notational convenience, we take $T=1$ throughout this section. 
We begin our analysis with a lemma that establishes the large-deviations behavior of the area under the busy period active at time $n$. To this end, define $\D[0,1]^{\leqslant 1} \triangleq \{\xi\in \D[0,1]: \xi = x\mathbbm1_{\{1\}} \text{ for some $x\geq 0$}\}$. Recall that $\bar S_n = \frac1n V_n \mathbbm1_{\{1\}}$ and $V_n =  \sum_{i=T_{N(n)}+1}^n \rvtxt{6-9}{6)}{\rvout{f(X_i)}\rvin{X_i^p}}$. 
\begin{lemma}
\label{exp-equiv-ldp-sn}
\linksinthm{exp-equiv-ldp-sn}
	$\bar{S}_n$ satisfies the LDP in $(\D[0,1], \mathcal{T}_{M_1'})$ with the speed $n^{\alpha}$ and the rate function $I_S: \D[0,1] \to \R_+$, where
		\begin{equation}
		\notationdef{nota-I-S}{I_S(\zeta)} \triangleq
		\begin{cases}
		{\mathcal{B}^*_{0}}(\zeta(1)-\zeta(1-))^{\alpha}  &  \mathrm{if} \ \zeta \in \D[0,1]^{\leqslant 1}, \\
		\infty & \ \text{otherwise}.
		\end{cases}
		\end{equation}
%	for every compact set $F \subseteq \R_+$
%	
%	\begin{align*}
%	\limsup_{n \to \infty}\frac{1}{n^{\alpha}}\log\P_0\left(\bar V_n \in F\right) \leq -\inf_{y \in F}\mathcal{I}_v(y),
%	\end{align*}
%	and for any open set $G \subseteq \R_+$
%	
%	\begin{equation*}
%	\liminf_{n \to \infty}\frac{1}{n^{\alpha}}\log\P_0\left(\bar V_n \in G\right) \geq -\inf_{y \in G}\mathcal{I}_v(y).
%	\end{equation*}
\end{lemma}

\proof{Proof.}
\linksinpf{exp-equiv-ldp-sn}
Define a function \linkdest{TtoUpsilon}$\Upsilon:\R_+  \to {\D[0,1]}^{\leqslant 1}$ as
$
\Upsilon(x)\triangleq x \cdot \mathbbm{1}_{\{1\}}
$.
Then, $\bar S_n = \Upsilon(\frac 1n V_n)$ and it is straightforward to see that $\Upsilon$ is a continuous function w.r.t. the $M_1'$ topology. \BZA{$T$ has been used notation-wise in other ways before, and we do not explain why this function is continuous contrary to other cases where we make continuity claims}\CRA{I added that it is straightforward to show the continuity. I think it will do. I also chaged $T$ to $\Upsilon$}
Therefore, the desired LDP follows from the contraction principle if we prove that $\frac1n V_n$ satisfies an LDP in $\R_+$ with the sub-linear speed $n^{\alpha}$ and the good rate function $I_V:\R_+ \to \R_+$, where
$
	\notationdef{nota-I-V}{I_V(x)}\triangleq{\mathcal{B}^*_{0}} \cdot x^{\alpha}.
$
	To prove the LDP for $\frac1n V_n$, note first that $\P(\frac1n V_n \in \cdot)$ is exponentially tight (w.r.t.\ the speed $n^{\alpha}$) from Theorem~\ref{subexpbusyperiod-steadystate}.
	Therefore, it is enough to establish the weak LDP.
	For the weak LDP, we start with showing that for any $a,b\in \R$, $B\triangleq(a,b)\cap \R_+$ satisfies
	$\limsup_{n \to \infty}\frac{\log\P\left(\frac1n V_n \in B\right)}{n^{\alpha}}
	=
	\liminf_{n \to \infty}\frac{\log\P\left(\frac1n V_n \in B\right)}{{n^{\alpha}}}$.
	Since this holds trivially with value $-\infty$ if $0\geq b$ or $a\geq b$, we assume that $0\vee a < b$.
	Note that, from Theorem~\ref{subexpbusyperiod-steadystate}, 
	\begin{align*}
%	&
	\limsup_{n \to \infty}\frac{\log\P\left(\frac1n V_n \in B\right)}{n^{\alpha}}
%	=
%	\limsup_{n \to \infty}\frac{\log \Big(\P_0\left(\bar V_n > a\right)-\P_0\left(\bar V_n \geq b\right)\Big)}{n^{\alpha}}
%		\\
%	&
%	\leq
%	\limsup_{n \to \infty}\frac{\log\left\{\P\left(\bar V_n \geq 0\vee a\right)\left(1-\frac{\P\left(\bar V_n \geq b\right)}{\P\left(\bar V_n \geq 0\vee a\right)}\right)\right\}}{n^{\alpha}}
%	\\
	&
	\leq
	\limsup_{n \to \infty}\frac{\log \P\left(\frac1n V_n \geq 0\vee a\right)}{n^{\alpha}}
	\leq -{\mathcal{B}_0^*}\cdot (0\vee a)^{\alpha}.
	\end{align*}
	\rvtxt{3}{3)}{\rvout{For small enough $\epsilon>0$,}\rvin{For the lower bound, note that  Theorem~\ref{subexpbusyperiod-steadystate} implies that} \rvout{and by using the property} ${\P(\frac1n V_n \geq b)}/{\P(\frac1n V_n \geq 0\vee a\rvin{+\epsilon})}\to 0$ \rvin{for any $\epsilon>0$. Therefore}}, we obtain
	\begin{align*}
	\liminf_{n \to \infty}\frac{\log\P\left(\frac1n V_n \in B\right)}{n^{\alpha}}
	&
	\rvtxt{3-1}{3)}{\;\rvin{
	    \geq
	    \liminf_{n \to \infty}\frac{\log\Big(\P\left(\frac1n V_n \geq 0\vee a+\epsilon\right) -\P\left(\frac1n V_n \geq b\right)\Big)}{n^{\alpha}}
	}}
	\\
	&
	\geq
	\liminf_{n \to \infty}\frac{\log\left\{\P\left(\frac1n V_n \geq 0\vee a+\epsilon\right)\left(1-\frac{\P\left(\frac1n V_n \geq b\right)}{\P\left(\frac1n V_n \geq 0\vee a+\epsilon\right)}\right)\right\}}{n^{\alpha}}
	\\
	&
	=
	\liminf_{n \to \infty}\frac{\log \P\left(\frac1n V_n \geq 0\vee a+\epsilon\right)}{n^{\alpha}}=-{\mathcal{B}_0^*}\cdot (0\vee a+\epsilon)^{\alpha}.
	\end{align*}
	Taking $\epsilon\to 0$, we see that the limit supremum and the limit infimum coincide.
	Since  $\{(a,b)\cap \R_+: \ a, b \in \R, \, a \leq b\}$ forms a base of the Euclidean topology on $\R_+$,  Theorem 4.1.11 of \cite{dembo2010large} applies, and hence, proves the desired weak LDP with the rate function $I_V$. This concludes the proof. 
	\Halmos
\endproof

We next work towards a  sample-path LDP for $\bar{Z}_n$. We employ a well-known technique, based on the projective limit theorem by Dawson and G\"artner; see Theorem 4.6.1 in \cite{dembo2010large}. 
The following three lemmas lead to the first key step in this approach, which consists of obtaining the finite-dimensional LDP for $\bar{Z}_n$.   

\begin{lemma}
\label{tails-asy-for-k-dim}
\linksinthm{tails-asy-for-k-dim}
	For any given $0=t_0<t_1<t_2<\ldots< t_k$, let $\Delta t_i \triangleq t_i - t_{i-1}$ for $i=1,\ldots,k$. Then,	 
	\begin{equation}\label{finite-dif-dim-up-bound}
	\limsup_{n \to \infty} \frac{1}{{n^{\alpha}} }\log \P\left(\sum_{j=1}^{N(nt_1)}W_{j} \geq n a_1, \ldots
		,\sum_{j=N(nt_{k-1})+1}^{N(nt_k)}W_j \geq n a_k
		\right)
	\leq -{\mathcal{B}_0^*}\left(\sum_{i=1}^{k}(a_i-\lambda\Delta t_i)_+^{\alpha}\right),
	\end{equation}
	\begin{equation}\label{finite-dif-dim-low-bound}
	\liminf_{n \to \infty} \frac{1}{{n^{\alpha}} }\log \P\left(\sum_{j=1}^{N(nt_1)}W_{j} \geq n a_1, \ldots
		,\sum_{j=N(nt_{k-1})+1}^{N(nt_k)}W_j \geq n a_k
		\right)
	\geq -{\mathcal{B}_0^*}\left(\sum_{i=1}^{k}(a_i-\lambda\Delta t_i)_+^{\alpha}\right),
	\end{equation}
	where $(x)_+ \triangleq x\vee 0$.
\end{lemma}
\proof{Proof.}
\linksinpf{tails-asy-for-k-dim}
Recall that $\tau_j = T_{j}-T_{j-1}$, where $T_j$ is the $j$\textsuperscript{th} hitting time of 0. 
Fix an arbitrary $\epsilon>0$ and let $E_i^{(n)}(\epsilon) \triangleq n[l_i, u_i]$, where $l_i \triangleq t_i/\E \tau_1 - \epsilon$ and $u_i \triangleq t_i/\E \tau_1 + \epsilon$.
We will use this notation throughout the proof of this lemma. 	
For the upper bound in Equation (\ref{finite-dif-dim-up-bound}), note that 
	\begin{align*}
	&  
	\P\left(\sum_{j=1}^{N(nt_1)}W_j \geq na_1
	,\ldots,
	\sum_{j=N(nt_{k-1})+1}^{N(nt_k)}W_j \geq na_k
	\right) \leq \underbrace{\sum_{i=1}^{k}\P\left(N(nt_i) \notin E_i^{(n)}(\epsilon)\right)}_{=\mathrm{(I)}} 
	\\
	& 
	+
	\underbrace{\P\left(\sum_{j=1}^{N(nt_1)}W_j \geq na_1, 
		\ldots,
		\sum_{j=N(nt_{k-1})+1}^{N(nt_k)}W_j \geq na_k, 
		N(nt_i) \in E_i^{(n)}(\epsilon) \ \mathrm{for} \ i=1,\ldots, k \right)}_{=\mathrm{(II)}}.
	\end{align*}
% {\color{red} A proof of this claim can be found in an ongoing paper by Nikolaos Karapiperakhs (to be appeared in arxiv). More details will be given shortly.}
\rvtxt{4}{4)}{
\rvout{Since $\tau$ is light-tailed, we use the definition of a renewal process and  Cram\'ers theorem to conclude that}\rvin{Note also that $\tau_1$ is light-tailed because
% \begin{align*}
$
\P(\tau_1 \geq k) 
= \P\left(\sum_{i=1}^j U_i \geq 0,\ j =1,\ldots,k\right)
\leq \P\left(\frac1k\sum_{i=1}^k U_i \geq 0\right),
$
and $\P\left(\frac1k\sum_{i=1}^k U_i \geq 0\right)$ decays at a geometric rate as $k\to\infty$ due to our assumptions that $\E U_1  <0$ and $\theta_+>0$, along with Cram\'er's Theorem. 
% \end{align*} 
Moreover,
\begin{align*}
\Big\{N(nt_i) \notin E_i^{(n)}\Big\} 
% &= \Big\{N(nt_i) < \lceil ns_i \rceil\Big\}\cup \Big\{N(nt_i) > \lfloor nu_i \rfloor\Big\}
% \\
&= \Big\{N\big(n(l_i+\epsilon)\E\tau_1\big) < \lceil nl_i \rceil\Big\}\cup \Big\{N\big(n(u_i-\epsilon)\E\tau_1\big) > \lfloor nu_i \rfloor\Big\}
\\
&\subseteq \left\{\frac{\sum_{i=1}^{\lceil nl_i\rceil} \tau_i }{\lceil nl_i \rceil} > \frac{n(l_i+\epsilon)}{\lceil nl_i \rceil}\E \tau_1\right\}
\cup
\left\{\frac{\sum_{i=1}^{\lceil nu_i\rceil} \tau_i }{\lceil nu_i \rceil} < \frac{n(u_i-\epsilon)}{\lceil nu_i \rceil}\E \tau_1\right\}.
\end{align*}
Therefore, again, by Cram\'er's theorem,}}
% 	we use the definition of a renewal process and Cram\'ers theorem. Since $\tau$ is light-tailed,
% \\
	\begin{equation}\label{light-tailasy}
	\limsup_{n \to \infty}\frac{\log\mathrm{(I)}}{n^{\alpha}}=-\infty.
	\end{equation}
	Shifting our attention to (II),
	\begin{align*}
	& 
	\P\left(\sum_{j=1}^{N(nt_1)}W_j \geq na_1, 
	\ldots,
	\sum_{j=N(nt_{k-1})+1}^{N(nt_k)}W_j \geq na_k, 
	N(nt_i) \in E_i^{(n)}(\epsilon) \ \mathrm{for} \ i=1,\ldots, k \right)
	\\
	&
	\leq
	\sum_{i_1=\lceil n(t_1/\E\tau_1-\epsilon)\rceil}^{\lfloor n(t_1/\E\tau_1+\epsilon)\rfloor}\cdots \sum_{i_k=\lceil n(t_k/\E\tau_1-\epsilon) \rceil}^{\lfloor n(t_k/\E\tau_1+\epsilon) \rfloor}
		\P\left(\sum_{j=1}^{i_1}W_j \geq na_1
		,\, 
		\ldots
		,\,
		\sum_{j=i_{k-1}+1}^{i_k}W_j \geq na_k,\,
		N(nt_l)=i_l  \mathrm{ \ for \ }  l=1,\ldots, k \right)
	\\
	&
	\leq 
	\sum_{i_1=\lceil n(t_1/\E\tau_1-\epsilon)\rceil}^{\lfloor n(t_1/\E\tau_1+\epsilon)\rfloor}\cdots \sum_{i_k=\lceil n(t_k/\E\tau_1-\epsilon) \rceil}^{\lfloor n(t_k/\E\tau_1+\epsilon) \rfloor}
		\P\left(\sum_{j=1}^{i_1}W_j \geq na_1
		,\, 
		\ldots
		,\,
		\sum_{j=i_{k-1}+1}^{i_k}W_j \geq na_k \right)I(i_1\leq \ldots\leq i_k)
	\\
	&
	=
	\sum_{i_1=\lceil n(t_1/\E\tau_1-\epsilon)\rceil}^{\lfloor n(t_1/\E\tau_1+\epsilon)\rfloor}\cdots \sum_{i_k=\lceil n(t_k/\E\tau_1-\epsilon) \rceil}^{\lfloor n(t_k/\E\tau_1+\epsilon) \rfloor}
		\P\left(\sum_{j=1}^{i_1}W_j \geq na_1\right)
		\cdots
		\P\left(\sum_{j=i_{k-1}+1}^{i_k}W_j \geq na_k \right) 
	\\
	&
	\leq 
	(2\epsilon n)^k \, \P\left(
	\sum_{j=1}^{\lfloor n(t_1/\E\tau_1+\epsilon) \rfloor}W_j \geq na_1 
	\right)
	\cdots
	\P\left(
	\sum_{j=\lceil n(t_{k-1}/\E\tau_1-\epsilon) \rceil}^{\lfloor n(t_k/\E\tau_1+\epsilon) \rfloor}W_j \geq na_k 
	\right).
	\end{align*}
Now, we have that from Result~\ref{sp-ldp-m1'-top} and Theorem~\ref{subexpbusyperiod-zero},
	\begin{align*}%\label{1stlogasy-proj-lmt-fin-dim}
	\limsup_{n \to \infty}\frac{1}{n^{\alpha}}\log \mathrm{(II)} \nonumber
	& 
	\leq
	\sum_{i=1}^{k}
	\limsup_{n \to \infty}\frac{1}{n^{\alpha}}\log 
	\P\left(
	\sum_{j=\lceil n(t_{i-1}/\E\tau_1-\epsilon) \rceil}^{\lfloor n(t_i/\E\tau_1+\epsilon) \rfloor}W_j
	\geq
	na_i
	\right) 
	+
	\limsup_{n \to \infty}\frac{\log (2\epsilon n)^k}{n^{\alpha}}
	\nonumber
	\\
	&
	\leq 
	-{\mathcal{B}_0^*}\sum_{i=1}^{k}(a_i-\lambda(\Delta t_i + 2\epsilon \E \tau_1)  )_+^{\alpha}.
	\end{align*}
Taking $\epsilon \to 0$, we arrive at
	\begin{equation}
	\label{1stlogasy-proj-lmt-fin-dim}
	\limsup_{n \to \infty}\frac{1}{n^{\alpha}}\log \mathrm{(II)}
	\leq 
	-{\mathcal{B}_0^*}\sum_{i=1}^{k}(a_i-\lambda \Delta t_i  )_+^{\alpha}.	
	\end{equation}
In view of \eqref{light-tailasy} and \eqref{1stlogasy-proj-lmt-fin-dim},
	\begin{align*}
	&
	\limsup_{n \to \infty}\frac{1}{n^{\alpha}}\log\P\left(\sum_{j=1}^{N(nt_1)-1}W_j \geq na_1
	, 
	\ldots,
	\sum_{j=N(nt_{i-1})}^{N(nt_i)-1}W_j \geq na_i,
	\ldots,
	\sum_{j=N(nt_{k-1})}^{N(nt_k)-1}W_j \geq na_k
	\right) 
	\\
	&
	\leq 
	\max\left\{
	\limsup_{n \to \infty}\frac{\log\mathrm{(I)}}{n^{\alpha}}, \limsup_{n \to \infty}\frac{\log\mathrm{(II)}}{n^{\alpha}}\right\}
	\leq 
	-{\mathcal{B}_0^*}\sum_{i=1}^{k}(a_i-\lambda\Delta t_i)_+^{\alpha}.
	\end{align*}
For the lower bound in Equation (\ref{finite-dif-dim-low-bound}), notice that 
	\begin{align}\label{projlimitLB}
& 	
	\P\left(\sum_{j=1}^{N(nt_1)-1}W_j > na_1
	,\ldots,
%	\sum_{j=N(nt_{i-1})}^{N(nt_i)-1}W_j > na_i
%	,\ldots,
	\sum_{j=N(nt_{k-1})}^{N(nt_k)-1}W_j > na_k
	\right) \nonumber
	\\
	& 
	\geq
	\P\left(
	\sum_{j=1}^{N(nt_1)-1}W_j > na_1
	, 
	\ldots,
	\sum_{j=N(nt_{k-1})}^{N(nt_k)-1}W_j > na_k
	,
	N(nt_i) \in E_i^{(n)}(\epsilon)  \ \text{for} \ i=1,\ldots,k   \right)  \nonumber
	\\
	&
	\geq 
	\P\left(\sum_{j=1}^{\lfloor n(t_1/E\tau_1-\epsilon) \rfloor-1} W_j > na_1,\ldots,
	\sum_{j=\lceil n(t_{k-1}/E\tau_1+\epsilon) \rceil}^{\lfloor n(t_k/E\tau_1-\epsilon) \rfloor-1}W_j > na_k
	, N(nt_i) \in E_i^{(n)}(\epsilon)  \ \text{for} \ i=1,\ldots,k   \right)  \nonumber
	\\
	&
	\geq
	\P\left(\sum_{j=1}^{\lfloor n(t_1/E\tau_1-\epsilon) \rfloor-1} W_j > na_1, \ldots,
		\sum_{j=\lceil n(t_{k-1}/E\tau_1+\epsilon) \rceil}^{\lfloor n(t_k/E\tau_1-\epsilon) \rfloor-1}W_j > na_i 
	\right)
	-\mathrm{(I)}
	\nonumber
	\\
%	\end{align}
%	The first term can be decomposed as follows:
%	\begin{align}\label{indi-Weib-asy}
%	\text{(III)}
	&
	=
	\P\left(\sum_{j=1}^{\lfloor n(t_1/E\tau_1-\epsilon) \rfloor-1} W_j > na_1\right)
	 \prod_{i=2}^{k}\P\left(
	\sum_{j=\lceil n(t_{i-1}/E\tau_1+\epsilon) \rceil}^{\lfloor n(t_i/E\tau_1-\epsilon) \rfloor-1}W_j > na_i 
	\right) - \mathrm{(I)}
	 \nonumber 	\\
	&	
	=
	\underbrace{
	\P\left(\sum_{j=1}^{\lfloor n(t_1/E\tau_1-\epsilon) \rfloor-1} W_j > na_1\right)
	 \prod_{i=2}^{k}\P\left(
	\sum_{j=1}^{\lfloor n(t_i/E\tau_1-\epsilon) \rfloor-\lceil n(t_{i-1}/E\tau_1+\epsilon) \rceil}W_j > na_i 
	\right)
	}_{=\mathrm{(III)}}
	-\mathrm{(I)}.
	\end{align}
	From Theorem~\ref{subexpbusyperiod-zero}, Result~\ref{sp-ldp-m1'-top} and \eqref{light-tailasy}, we get $\frac{\text{(I)}}{\text{(III)}} \to 0$ as $n\to\infty$.
	Therefore, \eqref{projlimitLB} leads to
	\begin{align*}
	&
	\liminf_{n \to \infty}\frac{1}{n^{\alpha}}\log 	\P\left(\sum_{j=1}^{N(nt_1)-1}W_j > na_1
	,\ldots,
	\sum_{j=N(nt_{k-1})}^{N(nt_k)-1}W_j > na_k
	\right)
	\\
&	
\geq 
	\liminf_{n \to \infty}\frac{1}{n^{\alpha}}\log\left\{
	\text{(III)}\left(1-\frac{\text{(I)}}{\text{(III)}}\right)
	\right\}
	=
	\liminf_{n \to \infty}\frac{1}{n^{\alpha}}\log
	\text{(III)}
	\\
	&
	=
	-{\mathcal{B}_0^*}\sum_{i=1}^{k}(a_i-\lambda(\Delta t_i - 2\epsilon \E \tau_1))_+^{\alpha}.
	\end{align*}
Taking $\epsilon \to 0$, we arrive at \eqref{finite-dif-dim-low-bound} concluding the proof.
		\Halmos
\endproof

Our next lemma establishes the LDP for $\left(\frac{1}{n}\sum_{j=1}^{N(nt_1)}W_j, \ldots
,\frac{1}{n}\sum_{j=N(nt_{k-1})+1}^{N(nt_k)}W_j
\right)$.

	\begin{lemma}
	\label{ldp-for-non-overlapping-partial-sums}
	\linksinthm{ldp-for-non-overlapping-partial-sums}
	For any given $\mathbf t=(t_1,\ldots,t_k)$ such that $0=t_0\leq t_1<\ldots<t_k\leq 1$, the probability measures $\mu_n$ of $\left(\frac{1}{n}\sum_{j=1}^{N(nt_1)}W_j, \ldots
	,\frac{1}{n}\sum_{j=N(nt_{k-1})+1}^{N(nt_k)}W_j
	\right)$ satisfy the LDP in $\R^k_+$ w.r.t.\ Euclidean topology with the speed $n^{\alpha}$ and the good rate function $I_{\mathbf t}:\R^{k}_+ \to \R_+$:
	\begin{equation}
	\notationdef{nota-I-t}{I_{\mathbf t}(x_1,\ldots,x_k)}\triangleq
	\begin{cases}
	{\mathcal{B}_0^*}\sum_{i=1}^k (x_i-\lambda\Delta t_i)^{\alpha} & \text{if} \ x_i\geq \lambda \Delta t_i,\ \forall i=1,\ldots, k,\\
	\infty & \ \text{otherwise}.  
	\end{cases}
	\end{equation}
\end{lemma}

\proof{Proof.}
\linksinpf{ldp-for-non-overlapping-partial-sums}
Note that it is straightforward from \eqref{finite-dif-dim-up-bound} of Lemma~\ref{tails-asy-for-k-dim} to see that $\left(\frac{1}{n}\sum_{j=0}^{N(nt_1)}W_j, \ldots,\frac{1}{n}\sum_{j=N(nt_{k-1})+1}^{N(nt_k)}W_j\right)$ is exponentially tight by considering compact sets $\prod_{i=1}^k[0,a_i]$ for sufficiently large $a_i$'s.
Also, we claim that 
$\left(\frac{1}{n}\sum_{j=0}^{N(nt_1)}W_j, \ldots
	,\frac{1}{n}\sum_{j=N(nt_{k-1})+1}^{N(nt_k)}W_j
	\right)$
	satisfies a weak LDP.
    Once this claim is established, Lemma 1.2.18 of \cite{dembo2010large} applies, showing that the full LDP is satisfied. 
\CRA{Are we claiming/using a result as follows: ``Weak LDP in Polish space + goodness of the rate function $\implies$ exponential tightness''? Is there a reference?}%
\CRA{I tried to fix it with minimal change using the upper bound \eqref{finite-dif-dim-up-bound} of Lemma~\ref{tails-asy-for-k-dim}. Please let me know if this looks fine.}% 
	Now, to prove the claimed weak LDP, we start by showing that
\begin{equation}\label{LALAequality}
	\underbrace{\limsup_{n \to \infty}\frac{\log \mu_n(A)}{n^{\alpha}}}_{\triangleq\overline{\mathcal L}_A}
	=
	\underbrace{\liminf_{n \to \infty}\frac{\log \mu_n(A)}{n^{\alpha}}}_{\triangleq\underline{\mathcal L}_A}
\end{equation}
for every $A \in\mathcal{A}\triangleq\{\prod_{i=1}^{k}\big((a_i,b_i)\cap \R_+\big) : \,  a_i < b_i \}$. 
Let $${\mathcal L}_A \triangleq
\begin{cases}
-{\mathcal B_0^*}\sum_{i=1}^k (a_i - \lambda \Delta t_i)_+^\alpha & \mathrm{if}\ b_i\geq \lambda \Delta t_i \ \mathrm{for}\ i=1,\ldots,k,\\
-\infty & \ \text{otherwise}.
\end{cases}
$$
%Since the above equality is trivial if $b_i \leq 0 $ or $a_i \geq b_i$ for any $i$, we assume that $0\vee a_i < b_i$ for $i=1,\ldots, k$.
%
%We will show both the upper bound $\overline{\mathcal L}_A  \leq {\mathcal L}_A$ and the lower bound ${\mathcal L}_A  \leq \underline{\mathcal L}_A $ by considering the following two cases seperately: 
%
We will prove \eqref{LALAequality} by showing that $\overline{\mathcal L}_A  \leq {\mathcal L}_A  \leq \underline{\mathcal L}_A $.
We consider two cases: 
\begin{itemize}[itemindent=40pt]
\item[Case 1.] $b_i \geq \lambda \Delta t_i$ for $i=1,\ldots, k$,
\item[Case 2.] $b_i < \lambda \Delta t_i$ for some $i \in \{1,\ldots, k\}$.
\end{itemize}
Let $A \triangleq \prod_{i=1}^k \big( (a_i, b_i)\cap \R_+\big)$ and $a_i < b_i$ for $i=1,\ldots, k$.
We start with Case 1. 
	Since $A \subseteq \prod_{i=1}^k[a_i, b_i)$,
	\begin{align}\label{overlineLAlessthanLA}
	\overline{\mathcal L}_A
	&
	\leq
	\limsup_{n \to \infty}\frac{1}{n^{\alpha}}\log	\P
	\left(
	\sum_{j=1}^{N(nt_1)}W_j \geq n a_1, \ldots
	,  \sum_{j=N(nt_{k-1})+1}^{N(nt_k)}W_j \geq n a_k
	\right)
%	\\
%	&
	\leq -{\mathcal{B}_0^*}\sum_{i=1}^{k}(a_i-\lambda\Delta t_i)^{\alpha} = \mathcal L_A,
	\end{align}
	where the second inequality is from \eqref{finite-dif-dim-up-bound}.
	Since $ \prod_{i=1}^k[a_i+\epsilon, b_i)\subseteq A$ for small enough $\epsilon>0$,
	\begin{align}\label{lowerbound-weakldp}
	\underline{\mathcal{L}_A}
	&
	\geq
\liminf_{n \to \infty}\frac{1}{n^{\alpha}}\log	\P
	\left(
	\left(\frac{1}{n} \sum_{j=0}^{N(nt_1)}W_j, \ldots
	,\frac{1}{n} \sum_{j=N(nt_{k-1})+1}^{N(nt_k)}W_j
	\right)
	\in 
	\prod_{i=1}^{k}[a_i+\epsilon,b_i)
	\right) \nonumber
	\\
	&
	\geq 
\liminf_{n \to \infty}\frac{1}{n^{\alpha}}\log 	
	\left\{
	\P
	\left(
	\frac{1}{n} \sum_{j=0}^{N(nt_1)}W_j > a_1+\epsilon, \ldots
	,\frac{1}{n} \sum_{j=N(nt_{k-1})+1}^{N(nt_k)}W_j > a_k+\epsilon
	\right) \right.
	\nonumber
	\\
	&
	\hspace{130pt}-		
	\left.
		\sum_{l=1}^k
		\P
		\left(
		\sum_{j=N(nt_{i-1})+1}^{N(nt_i)}W_j \geq na_i \, \forall \, i \neq l 
		, \sum_{j=N(nt_{l-1})+1}^{N(nt_l)}W_j \geq nb_l
		\right)
	\right\}
	\nonumber
	\\
	&
	\geq 
\liminf_{n \to \infty}\frac{1}{n^{\alpha}}\log 	\P
	\left(
	\frac{1}{n} \sum_{j=0}^{N(nt_1)}W_j > a_1+\epsilon, \ldots
	,\frac{1}{n} \sum_{j=N(nt_{k-1})+1}^{N(nt_k)}W_j > a_k+\epsilon
	\right) \nonumber
	\\
	&
	\hspace{10pt}+
	\liminf_{n \to \infty}\frac{1}{n^{\alpha}}\log 
	\left(
	1-
	\frac {
		\sum_{l=1}^{k}\P
		\left(
		\frac1n\sum_{j=N(nt_{i-1})+1}^{N(nt_i)}W_j \geq a_i+\epsilon \  \forall i \neq l 
		,\ \frac1n\sum_{j=N(nt_{l-1})+1}^{N(nt_l)}W_j \geq b_l
		\right)
		}{
		\P
		\left(
		\frac1n\sum_{j=1}^{N(nt_1)}W_j > a_1+\epsilon, \ldots
		, \frac1n\sum_{j=N(nt_{k-1})+1}^{N(nt_k)}W_j > a_k+\epsilon
		\right)
		} 
	\right).
	\end{align}	
	Note that, due to the logarithmic asymptotics of Lemma~\ref{tails-asy-for-k-dim}, for every $l \in \{1,\ldots,k\}$,
	\begin{equation*}%\label{lowerboundtrick}
	\frac{\P
		\left(
		\frac1n\sum_{j=N(nt_{i-1})+1}^{N(nt_i)}W_j \geq a_i+\epsilon\ \text{for} \ i \in \{1,\ldots,k\} \setminus l 
		,\ \frac1n\sum_{j=N(nt_{j-1})+1}^{N(nt_l)}W_j \geq b_l
		\right)}
	{\P
		\left(
		\frac1n\sum_{j=1}^{N(nt_1)}W_j >  a_1+\epsilon, \ldots
		,\ \frac1n\sum_{j=N(nt_{k-1})+1}^{N(nt_k)}W_j > a_k+\epsilon
		\right)} 
		\to 0,
	\end{equation*}
	and hence, the second term of \eqref{lowerbound-weakldp} disappears. Therefore,
	\begin{align*}
	\underline{\mathcal L}_A
	&\geq 
	\liminf_{n \to \infty} \frac{\log \P
		\left(
		\frac1n\sum_{j=1}^{N(nt_1)}W_j > a_1+\epsilon, \ldots
		,\ \frac1n\sum_{j=N(nt_{k-1})+1}^{N(nt_k)}W_j > a_k+\epsilon
		\right) }{n^{\alpha}}
	\\
	&\geq 
	-{\mathcal{B}_0^*}\sum_{i=1}^{k}(a_i+\epsilon-\lambda\Delta t_i)^{\alpha}.
	\end{align*}
	Taking $\epsilon\to0$, we arrive at $\underline{\mathcal L}_A \geq \mathcal L_A$,
	which, together with \eqref{overlineLAlessthanLA},
	proves \eqref{LALAequality}  for Case 1. 

	For Case 2, note that by Result~\ref{sp-ldp-m1'-top},
	$$
	\overline{\mathcal L}_A \leq 
	\limsup_{n \to \infty}\frac{1}{n^{\alpha}}\log	\P
	\left(
	\sum_{j=N(nt_{i-1})+1}^{N(nt_i)}W_j < n b_i
	\right)	
	= -\infty,
	$$
	and hence, $\overline{\mathcal L}_A = \underline{\mathcal L}_A = \mathcal L_A = -\infty$.

Now note also that 
\begin{equation}\label{I-tau-k}
I_{\mathbf t}(x_1,\ldots,x_k) = \sup \left\{-{\mathcal L}_A: A\in \mathcal A,  (x_1,\ldots,x_k)\in A\right\}.
\end{equation}
	Since $\mathcal A$ is a base of the Euclidean topology, the desired weak LDP follows from \eqref{LALAequality}, \eqref{I-tau-k}, and Theorem~4.1.11 of \cite{dembo2010large}.
		\Halmos
\endproof

The following is an immediate Corollary of Lemma~\ref{ldp-for-non-overlapping-partial-sums}.

\begin{lemma}
\label{fin-dim-LDP-process}
\linksinthm{fin-dim-LDP-process}
	For any given $\mathbf t=(t_1,\ldots,t_k)$ such that $0=t_0\leq t_1<\ldots<t_k\leq 1$, the probability measures $(\mu_n)$ of $\left(\frac{1}{n}\sum_{j=0}^{N(nt_1)}W_j, \ldots
	,\frac{1}{n}\sum_{j=0}^{N(nt_k)}W_j
	\right)$
	satisfy an LDP in $\R_+^k$ with the speed $n^{\alpha}$ and with the good rate function, $\tilde{I}_{\mathbf t}:\R^{k}_+ \to \R_+$,
	\begin{equation}
	\notationdef{nota-I-tilde-t}{\tilde{I}_{\mathbf t}(x_1,\ldots,x_k)}\triangleq
	\begin{cases}
	{\mathcal{B}_0^*}\sum_{i=1}^k (x_i-x_{i-1}-\lambda \Delta t_i)^{\alpha} & \mathrm{if} \ x_i-x_{i-1} \geq \lambda \Delta t_i, \mathrm{\ for\ } i=1,\ldots, k, \\
	\infty & \ \text{otherwise}.  
	\end{cases}
	\end{equation}
\end{lemma}

\proof{Proof.}
\linksinpf{fin-dim-LDP-process}
	The proof is an application of the contraction principle (\cite{dembo2010large}). To this end, consider the function $f:\R_+^k \to \R_+^k$, 
	$
	f(x_1,x_2,\ldots,x_k)\triangleq(x_1,x_1+x_2, \ldots, x_1+\ldots+x_k).
	$
	Notice that 
	$$\left(\frac{1}{n}\sum_{j=0}^{N(nt_1)}W_j, \ldots
	,\frac{1}{n}\sum_{j=0}^{N(nt_k)}W_j
	\right)
	=
	f\left(\frac{1}{n}\sum_{j=0}^{N(nt_1)}W_j, \ldots
	,\frac{1}{n}\sum_{j=N(nt_{k-1})+1}^{N(nt_k)}W_j
	\right),
	$$
	where $f$ is a continuous function.
	That is, $\left(\frac{1}{n}\sum_{j=0}^{N(nt_1)}W_j, \ldots
	,\frac{1}{n}\sum_{j=0}^{N(nt_k)}W_j
	\right)$ satisfies a large deviation principle with the rate function 
	$\tilde{I}_{\mathbf t}(y_1,\ldots,y_k)=\inf\{I_{\mathbf t}(x_1,\ldots,x_k): y =f(x_1,\ldots,x_k)\}$. 
	However, since $(y_1,\ldots,y_k)=f(x_1,\ldots,x_k)$,
	it is immediate that $y_1\leq y_2 \leq \ldots \leq y_k$. Therefore,

	\begin{equation*}
	\tilde{I}_{\mathbf t}(y_1,\ldots,y_k)=
	\begin{cases}
	{\mathcal{B}_0^*}\sum_{i=1}^k (y_i-y_{i-1}-\lambda \Delta t_i)^{\alpha} & \ \mathrm{if} \ y_{i+1} - y_i \geq \lambda \Delta t_i \mathrm{\ for\ } i=1,\ldots,k, \\
	\infty & \ \text{otherwise}.  
	\end{cases}
	\end{equation*}
		\Halmos
	
\endproof			

Now, for a path $\xi \in \D[0,1]$, let 
	\begin{equation}
	\notationdef{nota-I-Z}{I_Z(\xi)}\triangleq
	\begin{cases}
	{\mathcal{B}_0^*}\sum_{t: \xi(t) \neq \xi(t-)}(\xi(t) - \xi(t-))^{\alpha} & \mathrm{if} \ \xi \in \D^{(\lambda)}[0,1], \\
	\infty & \ \text{otherwise}.
	\end{cases}
	\end{equation}
\MBA{Why do we need the subscript $\alpha$? It does not serve any functionality from what I have seen so far.}\BZA{there are so many rate functions I, and they all have a subscript, so I prefer to keep it here also}\CRA{There seems to be a mistake here. Doesn't $\bar Z_n$ satisfy the LDP with $I_Z$ rather than $I_\alpha$? Is this a mistake? or am I missing something?}
Since $\bar{Z}_n$ satisfies a finite-dimensional LDP, we can show that the Dawson-G\"artner projective limit theorem  implies that $\bar{Z}_n$ satisfies an LDP in $\mathbb{D}[0,1]$ endowed with the pointwise convergence topology.
The next lemma verifies that the rate function associated with the LDP of $\bar{Z}_n$, is indeed $I_Z$.\JBA{There are a couple of issuesi involving the previous display, Lemma 3.4 and Lemma 3.5. The first issue is that Lemma 3.4 concerns a function in a k-dimensional space, where as the right hand side of Lemma 3.5 is a function in an infinite dimensional space. The corresponding projection notion applied to the left hand side of Lemma 3.5, while "obvious" it hasn't been really defined. The second issue is that the definition of $I_\alpha$ in the previous display lacks $B_0^*$. But then it is just $I_Z$.}\CRA{Thanks for catching this. I corrected the notation.}%
\begin{lemma}
\label{ratefunctionDG}
\linksinthm{ratefunctionDG}
	 Let \linkdest{temp_pin_1}$\mathbf T\triangleq\cup_{d=1}^\infty\{(t_1,\ldots,t_k): 0\leq t_1 < t_2 <\cdots < t_d \leq 1\}$ 
	 be the collection of all ordered (in the increasing order) finite subsets of $[0,1]$. 
	 \CRA{$t_1\leq t_2\leq \cdots t_k$?}% 
	 \BZA{I think so}\CRA{I changed it to DZ notation.}%
	 Then $$\sup_{\mathbf t \in \mathbf T}\tilde{I}_{\mathbf t}\big(\xi(t_1),\ldots,\xi(t_k)\big)=I_Z(\xi).$$
\end{lemma}	
\proof{Proof.}
\linksinpf{ratefunctionDG}
	The proof is essentially identical to the proof of Lemma~4 of \cite{gantert1998functional} and hence omitted. 	\Halmos
\endproof
We derive the sample-path LDP for the stochastic process $\bar{Z}_n$ w.r.t. the pointwise convergence topology, which we denote with $\notationdef{nota-cal-W}{\mathcal{W}}$. Recall that $\D^{(\lambda)}[0,1] $ denotes the subspace of increasing paths with slope $\lambda$. 
\begin{lemma}
\label{LDP-Z-n-pointwise-convergence}
\linksinthm{LDP-Z-n-pointwise-convergence}
	The stochastic process $\bar{Z}_n$ satisfies a large deviation principle in 
	$(\D[0,1], \mathcal{W})$,  
	with the speed $n^{\alpha}$ and the good rate function $I_Z$.
\end{lemma}
\proof{Proof.}
\linksinpf{LDP-Z-n-pointwise-convergence}
	This is an immediate consequence of the Dawnson and G\"artner projective limit theorem, (Theorem~4.6.1 of  \cite{dembo2010large}), Lemma~\ref{fin-dim-LDP-process} and Lemma~\ref{ratefunctionDG}. 
	\hfill\Halmos
\endproof

	Next, we  establish the sample-path LDP for the stochastic process $\bar{Z}_n$ in $(\D[0,1], \mathcal T_{M_1'})$.

\begin{lemma}
\label{lemma-ldp-Z}
\linksinthm{lemma-ldp-Z}
	The stochastic process $\bar{Z}_n$ satisfies a large deviation principle in 
	$\D[0,1]$ w.r.t.\ the $M_1'$ topology  
	with the speed $n^{\alpha}$ and the good rate function $I_Z$.
\end{lemma}

\proof{Proof.}
\linksinpf{lemma-ldp-Z}
	For the upper bound, consider a set $K_M \triangleq \{\xi \in \D[0,1]:\ \xi \ \mathrm{is\ nondecreasing,\ } \xi(0)\geq 0,\  \|\xi\|_{\infty} \leq M\}$. Let $F$ be a closed set in $(\D[0,1],\mathcal T_{M_1'})$. Then,
	\begin{align*}
	\limsup_{n \to \infty}\frac{1}{n^{\alpha}}\log\P\left(\bar{Z}_n \in F \right)
	&\leq 
	\limsup_{n \to \infty}\frac{1}{n^{\alpha}}\log	
	\left\{
		\P\left(\bar{Z}_n \in F \cap K_M \right) 
		+
		\P\left(\bar{Z}_n \in K_M^c \right)
	\right\}
	\\
	&
	\leq 
	\limsup_{n \to \infty}\frac{1}{n^{\alpha}}\log\bigg\{\P\left(\bar{Z}_n \in  F \cap K_M  \right) 
	+ \P\bigg(\sum_{j=1}^{N(nt)}W_j \geq M \bigg)
	\bigg\}.
	\end{align*}
	From Proposition $A.2$ of \cite{bazhba2017sample},  one can check that pointwise convergence in $K_M$ implies the convergence w.r.t.\ the $M_1'$ topology, and $K_M$ (and hence $F\cap K_M$ as well) is closed w.r.t.\ $\mathcal T_{M_1'}$. 
	Suppose that $\xi$ is in the closure of $F\cap K_M$ w.r.t.\ $\mathcal W$. 
	Then, because of the above mentioned properties of $K_M$, there exists a sequence of paths $\{\xi_n\}$ in $F\cap K_M$ such that  $\xi_n\to \xi$ w.r.t.\ $\mathcal T_{M_1'}$, which, in turn, implies that $\xi \in F\cap K_M$. 
	That is, $F\cap K_M$ is closed in $\mathcal W$ as well. 
	Now, applying the sample-path LDP w.r.t.\ $\mathcal W$ we have proved in Lemma~\ref{LDP-Z-n-pointwise-convergence}, and then picking $M$ large enough,
	\begin{align*}
	&
	\limsup_{n \to \infty}\frac{1}{n^{\alpha}}\log\P\left(\bar{Z}_n \in F \right) 
	\leq \max\left\{-\inf_{\xi \in  F \cap K_M}I_Z(\xi),\ -{\mathcal{B}_0^{\ast}}M\right\}
	 =
	-\inf_{\xi \in F \cap K_M}I_Z(\xi) \leq 	-\inf_{\xi \in F }I_Z(\xi).
	\end{align*}
	
	Moving on to the lower bound, let $G$ be an open set in $(\D[0,1],\mathcal{T}_{M_1'})$.  
	We assume that $I(G)< \infty$, since we have nothing to show otherwise. 
	Fix an arbitrary $\xi \in G \cap \D^{(\lambda)}[0,1]$, and let $k$ be such that an open ball of radius $\frac{1+\lambda}{k}$ around $\xi$ is inside of $G$. That is, $B_{M_1'}(\xi; \frac{1+\lambda}{k}) \triangleq \{\zeta \in \D[0,1]: d_{M_1'}(\xi, \zeta) < \frac{1+\lambda}{k}\} \subseteq G$.
	Note that, since $\xi \in \D^{(\lambda)}[0,1]$ and $\bar Z_n$ is non-decreasing, 
	$\{|\bar Z_n(i/k) - \xi(i/k)| < 1/k, \mathrm{\ for\ } i=0,\ldots,k\} \subseteq \{\bar Z_n \in B_{M_1'}(\xi; \frac{1+\lambda}{k})\}$.
Therefore, in view of Lemma~\ref{fin-dim-LDP-process}, 
	\begin{align*}
	\liminf_{n \to \infty}\frac{1}{n^{\alpha}}\log\P\left( \bar{Z}_n \in G \right)
	&
	\geq
	\liminf_{n \to \infty}\frac{1}{n^{\alpha}}\log\P\left( \bar{Z}_n \in B_{M_1'}\Big(\xi, \frac{1+\lambda}{k}\Big) \right)
	\\
	&
	\geq 
	\liminf_{n \to \infty}\frac{1}{n^{\alpha}}\log\P\left( |\bar{Z}_n(i/k) -\xi(i/k)| < 1/k, \ \mathrm{for} \ i=0,\ldots, k \right)
	\\
	&
	=
	-\inf_{(y_1,\ldots, y_{k}) \in \prod_{i=1}^{k}(\xi(i/k)-  1/k,\xi(i/k)+  1/k)}\tilde{I}_{\mathbf t}(y_1,\ldots,y_{k})
	\\
	&
	\geq 
	-{\mathcal B_0^*}^{(p)}\sum_{i=1}^{k}\Big(\xi\big(i/k\big)-\xi\big((i-1)/k\big)-\lambda/k\Big)^{\alpha}
	\\
	&
	\geq
	-{\mathcal B_0^*}^{(p)}\sum_{t:\xi(t) \neq \xi(t-)}\big(\xi(t) - \xi(t-)\big)^{\alpha}=-I_{Z}(\xi).
	\end{align*}
	Since $\xi$ was an arbitrary element of $G\cap \D^{(\lambda)}[0,1]$, we arrive at the desired lower bound:
	$$
	-\inf_{\xi \in G} I_Z (\xi) = -\inf_{\xi \in G\cap \D^{(\lambda)}[0,1]} I_Z (\xi) \leq \liminf_{n \to \infty}\frac{1}{n^{\alpha}}\log\P\left( \bar{Z}_n \in G \right).
	$$
\hfill\Halmos
\endproof

Our next lemma shows that $\bar Z_n + \bar S_n$ is exponentially equivalent to $\bar Y_n$.

\begin{lemma}
\label{exponentialequivalence-Y-R-Z}
\linksinthm{exponentialequivalence-Y-R-Z}
	$\bar{Y}_n$ and $\bar{Z}_n + \bar S_n$ are exponentially equivalent in $(\D[0,1], \mathcal{T}_{M_1'})$.
\end{lemma}

\begin{proof}{Proof.}
\linksinpf{exponentialequivalence-Y-R-Z}
Due to the construction of $\bar Y_n$, $\bar Z_n$, and $\bar S_n$, we have that for any $\delta>0$,
\begin{equation} \label{set-decomposition-M1p-dist}
\{ d_{M_1'}(\bar Y_n, \bar Z_n+\bar S_n) \geq \delta \} \subseteq \big\{(n-T_{N(n)})/n\geq\delta\big\} \cup \big\{ \exists j \leq N(n): \tau_j \geq n \delta \big\}.
\end{equation}
To bound the probability of the first set, define $D_n(\epsilon)\triangleq\{N(n)/n \geq 1/\E\tau_1-\epsilon\}$ for any $\epsilon>0$ and notice that
$
\P \left((n-T_{N(n)})/n>\delta\right) = \P\left(T_{N(n)}\leq n(1-\delta)\right) = \P\left(T_{N(n)}\leq n(1-\delta), D_n(\epsilon)\right)+\P\left(T_{N(n)}\leq n(1-\delta), D_n(\epsilon)^c\right),
$
and hence,
\begin{align}\label{eq:1partforexpequiv}
&\limsup_{n \to \infty}\frac{1}{n^{\alpha}}\log\P\left((n-T_{N(n)})/n\geq \delta\right) 
\nonumber\\
&\leq 
\limsup_{n \to \infty}\frac{1}{n^{\alpha}}\log \left\{\P\left(T_{N(n)}\leq n(1-\delta), D_n(\epsilon)\right)+\P\left( D_n(\epsilon)^c\right)\right\} \nonumber
\\
&
=
\max\left\{\limsup_{n \to \infty}\frac{1}{n^{\alpha}}\log \P\left(T_{N(n)}\leq n(1-\delta), D_n(\epsilon)\right), \, \limsup_{n \to \infty}\frac{1}{n^{\alpha}}\log \P\left( D_n(\epsilon)^c\right)\right\}.
\end{align}
%Consequently, we can infer the exponential equivalence of $\bar{S}_n$ and $\bar{R}_n$ by showing 
%\begin{equation}
%\limsup_{n \to \infty}\frac{1}{n^{\alpha}} \log \P\left(T_{N(n)}<n(1-\delta), D_n(\epsilon)\right) =-\infty.
%\end{equation} 
Letting $\epsilon < \delta/(2\E\tau_1)$, we see that from the same argument as the one preceding \eqref{light-tailasy},
\begin{align*}\label{eq:closetoone}
\limsup_{n \to \infty}\frac{1}{n^{\alpha}}\log\P\left(T_{N(n)}\leq n(1-\delta), D_n(\epsilon)\right) 
&
\leq 
\limsup_{n \to \infty}\frac{1}{n^{\alpha}}\log\P\left(T_{\left\lfloor {n\left(\frac{1}{\E\tau_1}-\epsilon\right)} \right\rfloor} \leq n(1-\delta), D_n(\epsilon)\right) \nonumber  \\
&
=
\limsup_{n \to \infty}\frac{1}{n^{\alpha}}\log\P\left(N(n(1-\delta)) \geq \left\lfloor {n\left(\frac{1}{\E\tau_1}-\epsilon\right)} \right\rfloor , D_n(\epsilon)\right)
\\
&=
-\infty.
\end{align*}
\MBA{To obtain the above result in the last equality, we have used Puhalskii's 97 paper of functional large deviations for first passage times. If the implication resulting in the last inequality could be sustained by another reference, it would be preferable.} \BZA{we do not even explain why it is true. It follows quickly from Cramers theorem as you know by now, but it needs to be explained} 
\CRA{Essentially the same argument is already provided right above \eqref{light-tailasy}. We can perhaps move the argument to the appendix and cite it in both places. If this sounds good, let me know. I will implement this.}\JBA{This argument doesn't need Puhalskii's paper, and it is sort of standard, we can either do what CH suggests or simply indicate, as done earlier in the argument preceding \eqref{light-tailasy}.}\CRA{I added that the argument is the same as the one above \eqref{light-tailasy}.}\\
Using the definition of a renewal process and  Cram\'ers theorem
we obtain
\begin{equation}\label{eq:negligiblen}
\limsup_{n \to \infty}\frac{1}{n^{\alpha}}\log \P\left( D_n(\epsilon)^c\right) =-\infty.
\end{equation}
Therefore, 
\begin{equation}\label{first-term-negative-infinite}
\limsup_{n \to \infty}\frac{1}{n^{\alpha}}\log\P\big((n-T_{N(n)})/n>\delta\big)=-\infty.
\end{equation}
Moving on to the bound for the probability of the second term in \eqref{set-decomposition-M1p-dist}, for  any $\epsilon>0$,
	\begin{align*}
	\P\left(\left\{ \exists j \leq N(n): \tau_j \geq n \delta \right\}\right)  
	& =
	\P\left( \exists j \leq N(n): \tau_j \geq n \delta, N(n)/n \leq 1/\E \tau_1 + \epsilon \right)
		+\P\left( N(n)/n > 1/\E\tau_1 + \epsilon\right)
	\\
	&
	\leq
	\P\left( \exists j \leq \left\lceil {n}/{\E(\tau_1)} +n\epsilon \right\rceil: \tau_j \geq n \delta \right)
		+\P\left( N(n)/n > 1/\E\tau_1 + \epsilon\right)
	\\
	&
	\leq
	\left\lceil {n}/{\E(\tau_1)} +n\epsilon\right\rceil \P\left( \tau_1 \geq n \delta \right)
		+\P\left( N(n)/n > 1/\E\tau_1 + \epsilon\right).
	\end{align*}
	\MBA{The equivalence seems to be correct using the $M_1'$ distance. We just add step functions that are positive. The only issue is the verification of the statement: In the $M_1'$ topology, addition is a continuous function as long as functions do not share jumps of the opposite sign at the same time. I believe this is not a statement that needs to be proven since, in my opinion, the proof will not be any different than the proof in the case of the $M_1$ topology.}\CRA{I don't understand this comment, but Mihail didn't seem to intend to call for any action regarding this. So I consider this as addressed.}
	Since $\P\left( \tau_1 \geq n \delta \right)$ and $\P\left( N(n)/n > 1/\E\tau_1 + \epsilon\right)$ both decay at an exponential rate, $$\limsup_{n\to\infty} \frac{1}{n^{\alpha}}\log\P\left(\left\{ \exists j \leq N(n): \tau_j \geq n \delta \right\}\right)    = -\infty.$$
	This, along with \eqref{first-term-negative-infinite} and \eqref{set-decomposition-M1p-dist}, proves the desired exponential equivalence.
	\hfill\Halmos
\end{proof}

Now, we have all the necessary components to prove Theorem~\ref{SPLDPYN}.

% \proof{Proof of Theorem~\ref{SPLDPYN}.}
% \subsection{Proof of Theorem~\ref{SPLDPYN}}
\begin{proof}{Proof of Theorem~\ref{SPLDPYN}.}
\linksinpf{SPLDPYN}
% \noindent Proof of Theorem~\ref{SPLDPYN}.
The preceding sequence of lemmas has resulted in LDPs of $\bar Z_n$ (Lemma \ref{lemma-ldp-Z}) and $\bar S_n$ (Lemma \ref{exp-equiv-ldp-sn}). Since $\bar Z_n$ and $\bar S_n$ are independent,
$(\bar{Z}_n, \bar{S}_n)$ satisfies an LDP in $\prod_{i=1}^{2}\D[0,1]$ with the rate function $\notationdef{nota-I-ZS}{I_{Z,S}(\zeta,\xi)}\triangleq I_{Z}(\zeta)+I_{S}(\xi)$; see, for example, Theorem 4.14 of \cite{ganesh2004big}.

Let $\phi: \prod_{i=1}^{2}\D[0,1] \to \D[0,1]$ denote the addition function $\phi(\xi,\zeta)=\xi+\zeta.$ 
Since $\phi$ is continuous on $(\xi,\zeta)$ as far as $\xi$ and $\zeta$ do not share a jump time with opposite directions (which follows from a straightforward modification of Lemma B.1 in \cite{bazhbaqueue}), $\phi$ is continuous on the effective domain of $I_{Z, S}$.
Let 
$
\notationdef{nota-I-W}{I_{\mathbf{W}}(\zeta)} \triangleq \inf\left\{ I_{Z,S}(\xi_{1},\xi_{2}): \zeta =\xi_1+ \xi_2,  \ \xi_1 \in \D^{(\lambda)}[0,1], \, \xi_2 \in \D^{\leqslant 1}[0,1] \right\}
$,
and note that it is straightforward to check that $I_W = I_Y$.
By the extended contraction principle (see \cite{puhalskii1997functional}), 
we conclude that $\bar Z_n+\bar S_n$ satisfies the sample-path LDP with the rate function $I_Y$. 

We now prove the large deviation upper bound. Let $F$ be a closed set w.r.t. the $M_1'$ topology, and let $F_\epsilon\triangleq \{\xi \in \D[0,1]: d_{M_1'}(\xi, F) \leq \epsilon\}$. Then, 
	\begin{align*}	
	&\limsup_{n \to \infty} \frac{1}{n^{\alpha}}\log \P\left( \bar{Y}_n \in F\right)
	\\
	&
	=
	\limsup_{n \to \infty} \frac{1}{n^{\alpha}}\log \bigg\{\P\left( \bar{Y}_n \in F,\ d_{M_1'}(\bar Y_n, \bar Z_n+\bar S_n)\leq \epsilon\right) + \P\Big(d_{M_1'}(\bar Y_n, \bar Z_n+\bar S_n)> \epsilon)\Big)\bigg\}
	\\
	&
	\leq
	\limsup_{n \to \infty} \frac{1}{n^{\alpha}}\log \P\left(\bar Z_n+ \bar S_n \in F_\epsilon\right)
	\leq
	-\inf_{\xi \in F_\epsilon}I_{Y}(\zeta),
	\end{align*}
where the first inequality is due to Lemma~\ref{exponentialequivalence-Y-R-Z}.
	Note that  
	$\lim_{\epsilon\to0}\inf_{\xi \in F_\epsilon}I_{Y}(\zeta) 
	=
	\inf_{\xi \in F}I_{Y}(\zeta), 
	$
	since $I_W$ is good w.r.t.\ $\mathcal T_{M_1'}$ (see Proposition A.3 of \cite{bazhba2017sample}). The desired large deviation upper bound follows by taking $\epsilon \to 0$.

	For the lower bound, let $G$ be an open set in $\mathcal T_{M_1'}$.
	We assume that $\inf_{\xi \in G} I_Y(\xi) <\infty$ since the lower bound is trivial otherwise. 
	For any given $\epsilon>0$, pick $\zeta\in G$ such that  $I(\zeta) \leq \inf_{\xi \in G} I_Y(\xi) + \epsilon$. Let $\delta>0$ be such that $B_{M_1'}(\zeta, 2\delta) \in G$.
	Then, we know from Lemma~\ref{exponentialequivalence-Y-R-Z},
	$\P\left( d(\bar Y_n, \bar Z_n + \bar S_n) < \delta\right)/\P\left(\bar Z_n+\bar S_n \in B_{M_1'}(\zeta, \delta)\right)\to 0$, and hence,
	\begin{align*}
	&
	\liminf_{n \to \infty} \frac{1}{n^{\alpha}}\log \P\left( \bar{Y}_n \in G\right)
	\\
	&
	\geq
	\liminf_{n \to \infty} \frac{1}{n^{\alpha}}\log \P\left(\bar Z_n+\bar{S}_n \in B_{M_1'}(\zeta, \delta),\  d(\bar Y_n, \bar Z_n + \bar S_n) < \delta\right)
	\\
	&
	\geq 	
	\liminf_{n \to \infty} \frac{1}{n^{\alpha}}\log \P\left(\bar Z_n+\bar{S}_n \in B_{M_1'}(\zeta, \delta)\right)
	\left\{ 1 - \frac{\P\left( d(\bar Y_n, \bar Z_n + \bar S_n) < \delta\right)}{\P\left(\bar Z_n+\bar S_n \in B_{M_1'}(\zeta, \delta)\right)}
	\right\}
	\\
	&
	=
	\liminf_{n \to \infty} \frac{1}{n^{\alpha}}\log \P\left(\bar Z_n+\bar S_n \in B_{M_1'}(\zeta, \delta)\right)
	\\
	&
	\geq
	-\inf_{\xi \in B_{M_1'}(\zeta, \delta)}I_{\mathbf{Y}}(\xi)
	\geq 
	-I_Y(\zeta) \geq -\inf_{\xi \in G} I_Y(\xi) - \epsilon
	.
	\end{align*}
	Taking $\epsilon \to 0$, we arrive at the desired lower bound.
	\Halmos

% \endproof
\end{proof}

\section{Tail asymptotics for the area of a busy period.}
%Proofs of Theorem~\ref{subexpbusyperiod-zero} and %Theorem~\ref{subexpbusyperiod-steadystate}}
%Proposition \ref{bp-ison-bzero}, Theorem \ref{subexpbusyperiod-zero}, and Theorem %\ref{subexpbusyperiod-steadystate}}
\label{sec:firstproofsection}

Our focus in this section is on proving Theorem~\ref{subexpbusyperiod-zero} and Theorem~\ref{subexpbusyperiod-steadystate}. In Section \ref{subsection-concavity}, we collect several analytic properties of key variational problems related to these theorems. The proofs of these analytic properties as well as some other analytic results are deferred to Section~\ref{additionalproofsection} to keep the focus of this section on probabilistic arguments as much as possible. In Section~\ref{subsec:proof-of-tail-asymptotics}, we state two key propositions (Proposition \ref{technical-proposition-for-tail-asymptotics} and Proposition \ref{miscellaneous-proposition-for-tail-asymptotics}) which are applied to providing the proofs of Theorem~\ref{subexpbusyperiod-zero} and Theorem~\ref{subexpbusyperiod-steadystate}. 
The rest of this section is devoted to the proofs of Proposition 
\ref{subsec:proof-of-tail-asymptotics} and \ref{miscellaneous-proposition-for-tail-asymptotics}.
The focus is again on probabilistic ideas, the substantial number of additional analytical arguments that are directly needed are stated as lemmas of which the proofs can be found in Section \ref{additionalproofsection}.

%Section~\ref{preliminary-results-By} collects preliminary results useful for the rest of the proofs. 
%Section~\ref{subsection-concavity} proves Proposition~\ref{bp-ison-bzero}.
% Section~\ref{subsec:proof-of-tail-asymptotics}  collects key technical results (Proposition~\ref{technical-proposition-for-tail-asymptotics}, Proposition~\ref{miscellaneous-proposition-for-tail-asymptotics})  and then provides the main proof of the tail asymptotics of $W_n$ and $\bar V_n$ (Theorem~\ref{subexpbusyperiod-zero} and Theorem~\ref{subexpbusyperiod-steadystate}). Section~\ref{subsec:proof-of-technical-proposition-for-tail-asymptotics} and Section~\ref{subsec:proof-of-miscellaneous-proposition-for-tail-asymptotics} provide the proofs of the key technical results, i.e., Proposition~\ref{technical-proposition-for-tail-asymptotics} and Proposition~\ref{miscellaneous-proposition-for-tail-asymptotics}, respectively. 

\subsection{Key auxiliary variational problems and related properties}
\label{subsection-concavity}

%Proofs of Proposition~\ref{bp-ison-bzero} and Lemma~\ref{F-z-T-ccompact}}

Recall that $B^{\AC}_y = B_y \cap \AC[0,\infty)$. %{\color{red}  In the following lemma we reduce the search for the optimal path of the variational problem $(\mathcal{B}_y)$ to only absolutely continuous paths, and we show that the optimal value $\mathcal{B}_y^*$ is $0$ for $y$ sufficiently large.} 
Our first lemma establishes that the infimum defining $\mathcal{B}^*_y$, taken over paths of bounded variation, can be confined to absolutely continuous paths. Its proof is given in Section \ref{pf-useful-properties-of-variational-By}.

 \begin{lemma}
 \label{useful-properties-of-variational-By}
 \linksinthm{useful-properties-of-variational-By}
Recall that $\mathcal{B}^*_y$ is the optimal value of the variational problem $(\mathcal{B}_y)$.
\begin{itemize}
\item[(i)] Let \notationdef{nota-bar-y}{$\bar{y} \triangleq (|\mu|(p+1))^{\alpha}$}. For any $y\geq \bar y$, there exists a path $\xi^* \in B_y^{\AC}$ so that $I_y(\xi^*)=0$ and
$\B_{y}^* = 0.$ 
\item[(ii)] For any $y \geq 0$, $\mathcal{B}_y^*= \inf_{\xi \in B_y^{\mathbbm{AC}}} I_y(\xi)$.
%\item[(iii)] the mapping $y \to \mathcal{B}_y^*$ is non-increasing.
\end{itemize}

\end{lemma}

{We provide the proof of the following proposition in Section \ref{pf-eq:ineqformonotonicity}. It facilitates the proof of  Proposition~\ref{miscellaneous-proposition-for-tail-asymptotics}, which is a key result for the tail asymptotics of $W_1$ and $V_n$.}

\begin{proposition}
\label{eq:ineqformonotonicity}
\linksinthm{eq:ineqformonotonicity}
	The optimal value ${\mathcal{B}^*_y}$ of   $(\mathcal{B}_y)$ satisfies the following properties:
	\begin{itemize}
		\item[(i)] $y\mapsto \B_y^*$ is non-increasing in $y$, $y \in [0,\bar y]$. 
		\item[(ii)] 
		$y\mapsto\B_y^*$ is Lipschitz continuous, $y \in [0,\bar y]$.
		\item[(iii)] for every $y \in [0,\bar y)$, $\mathcal{B}_y^* \in (0, \infty)$. 
		%and $\mathcal B_0^* < \infty$.
        %
		%\lim_{y_k \downarrow y}{\mathcal{B}^*_{y_k}} \rightarrow {\mathcal{B}^*_y}$.
		%\item[iv)] For any sequence  $y_k$, so that $y_k \downarrow y^*$ it holds that
		%\begin{equation*}
		%\lim_{k \to \infty}{\mathcal{B}^*_{y_k}}^{(p)} + \mathcal{J}(y_k)
		% \rightarrow {\mathcal{B}^*_{y}}+\beta y.
		%\end{equation*}
		
	\end{itemize}
\end{proposition}

%\todo{Jose wrote: THERE ARE TOO MANY RESULTS IN JUST A SINGLE LEMMA AND THE PROOF IS A BIT ON THE TERSE SIDE. ALSO, MAYBE WE SHOULD TRY TO ORGANIZE THE PROOFS BY ONLY STATING RESULTS THAT ARE KEY FOR THE PROOF OF THE MAIN THEOREMS WHILE RELEGATING ANY AUXILIARY RESULT AND ESPECIALLY ITS PROOF BACK TO SOME APPENDIX OR LATER IN THE CORRESPONDING SECTION. Bert: given the following lemma is now relegated to this section, I suggest to leave it as is. }
Fix  $T>0$ and consider a functional
$\Phi_T: \D[0,T] \to \R_+$, where $\notationdef{nota-Phi-T}{\Phi_T(\xi)}=\int_{0}^{T}(\Psi(\xi)(s))^pds$. 
Now, let ${\mathcal V_y^T}^*$ denote the optimal value of the  optimization problem
\begin{equation}
    \label{eq:variationalproblemV}
\notationdef{nota-cal-V-y-T-*}{{\mathcal V_y^T}^*} \triangleq  \inf_{\xi \in V_{y}^T}I_{y}^\BVT(\xi), \notationdef{nota-cal-V-y-T}{\tag{$\mathcal{V}_y^{T}$}},
 \end{equation}
where 
$$\notationdef{nota-V-y-T}{V_{y}^{T}} \triangleq \left\{ \xi \in \D[0,T]: \xi(0)=y,   \ \Phi_T(\xi) \geq 1 \right\},$$ 
and
$$
\notationdef{nota-I-y-BV}{I_y^\BVT(\xi)} \triangleq \begin{cases}
\int_{0}^{T} \Lambda^{*}(\dot{\xi}^{(a)}(s)) d s+\theta_{+} \xi^{(u)}(T)+\theta_{-} \xi^{(d)}(T)
&
\text{if } \xi(0) = y \text{ and } \xi \in \BV[0,T] ,
\\
\infty
&
\text{otherwise.}
\end{cases}
$$
The variational problem $(\mathcal{V}_y^T)$ will naturally appear in large deviations estimates. The next lemma, proved in Section \ref{pf-merging-actions-to-the-left-lemma}, summarizes several of its properties.
%, (see below for the definition of $(\mathcal{V}_y^T)$) from functions of bounded variations to absolutely continuous functions.
%Later on, we apply this lemma to the proof of Proposition~\ref{bp-ison-bzero}.
\begin{lemma}
\label{merging-actions-to-the-left-lemma}
\linksinthm{merging-actions-to-the-left-lemma}
Consider an arbitrary path $\xi \in \BV[0,T]$ and set $y \triangleq \xi(0)$. 
\begin{itemize}
\item[(i)]
There exists a path $\zeta_1 \in \BV[0,T]$ such that
\begin{itemize}
\item[i-1)]
$\zeta_1(0) = y$,
\item[i-2)]
$\Phi_T(\zeta_1) \geq \Phi_T(\xi)$,
\item[i-3)]
$I_{y}^\BVT(\zeta_1) \leq I_y^\BVT(\xi)$,
\item[i-4)]
For some $t\in [0, T]$, 
$\zeta_1$ is nonnegative over $[0,t]$ and $\zeta_1$ is linear with slope $\mu$ over $[t, T]$.
%\item[i-5)]
%$\T(\zeta_1) \leq \T(\xi)$
\end{itemize}

\item[(ii)]
There exists a path $\zeta_2 \in \AC[0,T]$ and a  $z \in[0, \xi^{(u)}(T)]$ such that %\cmt{BZ: I have put the z up here}
\begin{itemize}
\item[ii-1)]
$\zeta_2(0) = y + z$,
\item[ii-2)]
$\Phi_T(\zeta_2) \geq \Phi_T(\xi)$,
\item[ii-3)]
$\theta_+ \cdot z + I_{y+z}^\BVT(\zeta_2) \leq I_y^\BVT(\xi)$, where we interpret $\theta_+\cdot z$ as $0$ if $\theta_+ = \infty$ and $z = 0$,
\item[ii-4)]
For some $t\in [0, T]$, 
$\zeta_2$ is nonnegative over $[0,t]$ and $\zeta_2$ is linear with slope $\mu$ over $[t, T]$.
%\item[i-5)]
%$\T(\zeta_2) \leq \T(\xi)$
\end{itemize}
\end{itemize}
\begin{itemize}
\item[(iii)]
If, in addition, $\xi\in\AC[0,T]$,
there exists a path $\zeta_3\in \AC[0,T]$ such that
\begin{itemize}
\item[iii-1)]
$\zeta_3(0) = y$,
\item[iii-2)]
$\Phi_T(\zeta_3) \geq \Phi_T(\xi)$,
\item[iii-3)]
$I_{y}^\BVT(\zeta_3) \leq I_y^\BVT(\xi)$,
\item[iii-4)]
$\zeta_3$ is concave over $[0,T]$ and its derivative is bounded by $\mu$ from below.
\end{itemize}
\end{itemize}
\end{lemma}

\subsection{Proof of Theorem~\ref{subexpbusyperiod-zero} and Theorem~\ref{subexpbusyperiod-steadystate}}
\label{subsec:proof-of-tail-asymptotics}
%The proof of Theorem~\ref{subexpbusyperiod-zero} and Theorem~\ref{subexpbusyperiod-steadystate} is based on the following key propositions. 
The following propositions are instrumental.
%Lemma~\ref{technical-proposition-for-tail-asymptotics}, in particular, contains the key inequalities. 
%The proofs of the lemmas are deferred to \todo{links to Section 3.2 and Section 3.3}

%\cmt{CR: we are switching between $X(\lfloor ux\rfloor)$ and $X_{\lfloor ux\rfloor}$ without warning. We should choose one and be consistent.}
%\cmt{BZ: agree, but where has it been defined? In section 4 I have changed everything into $X(\lfloor ux\rfloor)$ in line with other continuous-time processes. }

\begin{proposition}
\label{technical-proposition-for-tail-asymptotics}
\linksinthm{technical-proposition-for-tail-asymptotics}
\begin{itemize}

\item[(i)] Recall that $T_1 = \inf\{ k> 0: X_k = 0\}$. Then,
$$
\limsup_{x\to\infty} \frac 1 x \log \P_{xy}\bigg(\int_{0}^{T_1/x} (X_{\lfloor ux \rfloor}/x)^p du \geq 1 \bigg)
\leq 
-{\B_y^*}.
$$

\item[(ii)]
Recall that $W_1=\sum_{i=1}^{T_1}X_i^p$. Then,
\begin{equation*}
\liminf_{ u \to \infty}\frac{1}{u^{\alpha}}\log\P_0\left(W_1 >u\right) \geq  -\mathcal{B}_0^*.
\end{equation*}

\end{itemize}
\end{proposition}
We provide the proof of Proposition~\ref{technical-proposition-for-tail-asymptotics} in Section~\ref{subsec:proof-of-technical-proposition-for-tail-asymptotics}.

\begin{proposition}
\label{miscellaneous-proposition-for-tail-asymptotics}
\linksinthm{miscellaneous-proposition-for-tail-asymptotics}
\begin{itemize}
\item[(i)] $\sum_{k=0}^{m-1} X_k^p > x^{1+p}$ $\iff$ $\int_0^{m/x} \left(\frac{X_{\lfloor ux \rfloor}}{x}\right)^p du>1$,

\item[(ii)]
	it holds that 
	$$\mathcal{B}_0^*= \mathcal{B}_{\pi}^*,$$

\item[(iii)] finally, 
$$\lim_{k\to \infty} \min_{i\geq 1} \Bigg\{ 
\frac{i-1}k \beta\bar y + \B_{\frac i k \bar y}^*
\Bigg\}
= \inf_{y\in [0,\infty)} \left\{\beta y + \mathcal B_y^*\right\}
=\mathcal B_\pi^*.
$$

\end{itemize}
\end{proposition}

We provide the proof of Proposition~\ref{miscellaneous-proposition-for-tail-asymptotics}  in Section~\ref{subsec:proof-of-miscellaneous-proposition-for-tail-asymptotics}.

With Proposition~\ref{technical-proposition-for-tail-asymptotics} and Proposition~\ref{miscellaneous-proposition-for-tail-asymptotics} in our hands, we are ready to prove Theorem~\ref{subexpbusyperiod-zero}.

\begin{proof}{Proof of Theorem~\ref{subexpbusyperiod-zero}.}
\linksinpf{subexpbusyperiod-zero}
For the upper bound, setting $t= x^{p+1} = x^{1/\alpha}$, 
\begin{align*}
\limsup_{t\to\infty} \frac1 {t^{\alpha}} \log \P(W_1 \geq t)
&=
\limsup_{x\to\infty}\frac 1 x  \log \P_0\left(\sum_{k=0}^{T_1-1} X_k^p \geq x^{1+p} \right)
\\
&
=
\limsup_{x\to\infty}\frac 1 x  \log
\P_{0} \left( \int_{0}^{T_1/x} \left(\frac{X_{\lfloor ux \rfloor}}{x}\right)^p du \geq 1 \right)
\\
&
\leq
-\B_0^*.
\end{align*}
We used $\sum_{i=1}^{T_1} X_i^p = \sum_{i=0}^{T_1-1} X_i^p$ to obtain the first equality,  applied part \emph{(i)} of Proposition~\ref{miscellaneous-proposition-for-tail-asymptotics} to derive the second equality, and invoked part \emph{(i)} of Proposition~\ref{technical-proposition-for-tail-asymptotics} to establish the inequality.
Together with the matching lower bound in part \emph{(ii)} of Proposition~\ref{technical-proposition-for-tail-asymptotics}, this yields the desired asymptotics \eqref{W-1-asymptotics}.
\eop
\end{proof}

The proof of Theorem~\ref{subexpbusyperiod-steadystate} is slightly more involved.

\begin{proof}{Proof of Theorem~\ref{subexpbusyperiod-steadystate}.}
\linksinpf{subexpbusyperiod-steadystate}
We start by proving the large deviation upper bound for $\frac1n V_n$. 
Denote the time-reversed Markov process of $\{X_k, k=1,\ldots,n\}$ with $\{X_k^*, k = 0,\ldots,n\}$, and let $T_1^*\triangleq\inf\{i > 0: X_i^* = 0\}$. 
Let $\bar y \triangleq (|\mu|(p+1))^{\alpha}$ and fix $b>0$. 
Setting $x^{p+1} \triangleq nb$, we obtain that 
% \rvtxt{1-1}{1)}{\rvin{From the stochastic monotonicity of Markov chain $X_n$,}}
\begin{align}
% \rvtxt{1}{1)}{
% \rvin{\P_0\left( \frac1n V_n \geq b\right)}
% }
% &
% \leq
\P_{
% \rvtxt{1}{1)}{\rvout{0}\rvin{\pi}}
0}\left( \frac1n V_n \geq b\right)
&=
\P_{
% \rvtxt{1}{1)}{\rvout{0}\rvin{\pi}}
0} \left( \sum_{i=T_{N(n)}+1}^{n} X_i^p \geq nb \right)
=
\frac{1}{\pi(0)}\P_{\pi} \left( \sum_{i=0}^{T_1^*-1} (X_i^*)^p \geq nb, X_n^*=0 \right)
\nonumber
\\
&
\leq
\frac{n+1}{\pi(0)}
\P_{\pi} \left(\, \sum_{i=0}^{T_1-1} X_i^p \geq nb \right)
=
\frac{n}{\pi(0)}\P_{\pi} \left( \int_{0}^{T_1/x} \left(\frac{X_{\lfloor ux \rfloor}}{x}\right)^p du \geq  1 \right),
\label{first-inequality-in-this-proof}
\end{align}
where the second equality follows from Lemma~\ref{timereversibility} 
with $g(y_0,...,y_n) = \I( \sum_{\max \{i\leq n: y_i=0\}}^n y_i^p > nb)$, the second inequality follows from the upper bound in Lemma \ref{lastcycleequivalence}, and the last equality follows from part \emph{(i)} of Proposition~\ref{miscellaneous-proposition-for-tail-asymptotics}.

From the tower property, we have that 
\begin{align}\label{eq:rescalingtheprobability}
&
\P_{\pi} \left( \int_{0}^{T_1/x} (X_{\lfloor ux \rfloor}/x)^p du \geq 1 \right)
\nonumber
\\
&
=
\E_\pi\Bigg[\bigg(\ind\{X_0 \geq x\bar y\} + \sum_{i=1}^k\ind\bigg\{X_0\in \bigg[\frac{i-1}k x\bar y,\, \frac ik x\bar y\bigg] \bigg\}\bigg)\,\P\bigg(\int_{0}^{T_1/x} (X_{\lfloor ux \rfloor}/x)^p du \geq 1 \bigg| X_0) \bigg)\Bigg]
\nonumber
\\
&
\leq
\E_\pi \ind\{X_0 \geq x\bar y\} + 
\sum_{i=1}^k\E_\pi\bigg[\ind\bigg\{X_0\in \bigg[\frac{i-1}k x\bar y,\, \infty\bigg)\bigg\}\,\P_{\frac ik x\bar y}\bigg(\int_{0}^{T_1/x} (X_{\lfloor ux \rfloor}/x)^p du \geq 1 \bigg)\Bigg]
\nonumber
\\
&
\leq
\pi [x\bar y,\infty) + 
\sum_{i=1}^k\pi\bigg[\frac{i-1}k x\bar y,\, \infty\bigg)\,\P_{\frac ik x\bar y}\bigg(\int_{0}^{T_1/x} (X_{\lfloor ux \rfloor}/x)^p du \geq 1 \bigg),
\end{align}
where in the first inequality we used  that the Markov chain $\{X_n, n \geq 1\}$ is monotone with respect to the initial state. %\rvin{In our case, we observe that the stochastic process $X_n$ is increasing with respect to the initial state $X_0$}.
\BZA{is this concept known enough to be invoked without clarification?}\MBA{The notion of monotonicity which we use here only uses the initial state. I made a comment about this.}\CRA{I think it suffices to make it clear that what we mean is monotonicity w.r.t.\ the initial state.}
Therefore, by the principle of the maximum term and part \emph{(i)} of Proposition~\ref{technical-proposition-for-tail-asymptotics},
\begin{align}
&
\limsup_{x\to\infty} \frac 1 x \log \P_{\pi} \left( \int_{0}^{T_1/x} (X_{\lfloor ux \rfloor}/x)^p du \geq 1 \right)
\nonumber 
\\
&
\leq
\limsup_{x\to\infty} \frac 1 x \log \pi [x\bar y,\infty)  
\nonumber
\\
&
\qquad \vee \max_{i=1,\ldots,k} \Bigg\{ 
\limsup_{x\to\infty} \frac 1 x \log \Bigg(\pi\bigg[\frac{i-1}k x\bar y,\, \infty\bigg)\,\P_{\frac ik x\bar y}\bigg(\int_{0}^{T_1/x} (X_{\lfloor ux \rfloor}/x)^p du \geq 1 \bigg)\Bigg)
\Bigg\}
\nonumber
\\
&
= 
(-\beta \bar y) \vee \max_{i=1,\ldots,k} \Bigg\{ 
-\frac{i-1}k \beta\bar y + \limsup_{x\to\infty} \frac 1 x \log \P_{\frac ik x\bar y}\bigg(\int_{0}^{T_1/x} (X_{\lfloor ux \rfloor}/x)^p du \geq 1 \bigg)
\Bigg\}.
\nonumber
\\
&
\leq
(-\beta \bar y) \vee \max_{i=1,\ldots,k} \Bigg\{ 
-\frac{i-1}k \beta\bar y - \B_{\frac ik \bar y}^*
\Bigg\}.
\nonumber
\end{align}
Note that since $\B_{y}^* = 0$ for $y \geq \bar y$ due to part \emph{(i)} of Lemma~\ref{useful-properties-of-variational-By},
$$
(-\beta \bar y) \vee \max_{i=1,\ldots,k} \Bigg\{ 
-\frac{i-1}k \beta\bar y - \B_{\frac ik \bar y}^*
\Bigg\}
=
\max_{i\geq 1} \Bigg\{ 
-\frac{i-1}k \beta\bar y -\B_{\frac i k \bar y}^*
\Bigg\}
=
-\min_{i\geq 1} \Bigg\{ 
\frac{i-1}k \beta\bar y + \B_{\frac i k \bar y}^*
\Bigg\}.
$$
Taking $k\to \infty$ and applying part \emph{(ii)} and \emph{(iii)} of Proposition~\ref{miscellaneous-proposition-for-tail-asymptotics},
\begin{align*}
\limsup_{x\to\infty} \frac 1 x \log \P_{\pi} \left( \int_{0}^{T_1/x} (X_{\lfloor ux \rfloor}/x)^p du \geq 1 \right)
\leq 
-\B_\pi^*
= -\B_0^*.
\end{align*}
From this, along with \eqref{first-inequality-in-this-proof}, we arrive at the desired upper bound:
\[
\limsup_{x\to\infty} \frac 1 {n^{\alpha}} \log \P_{0}\left( \frac1n V_n \geq b\right) 
\leq
\limsup_{x\to\infty} \frac 1 x \log \P_{\pi} \left( \int_{0}^{T_1/x} (X_{\lfloor ux \rfloor}/x)^p du \geq 1 \right)\cdot b^{\alpha}
\leq -\B_0^*\cdot b^{\alpha}.
\]
Next, for $n$ sufficiently large, using the lower bound of Lemma \ref{lastcycleequivalence} for $n\geq n_0$:
\begin{align}\label{eq:1st-lower-bound}
 \P_{0}\left( \frac1n V_n \geq b\right)
& =
\P_{0} \left( \sum_{i=T_{N(n)}+1}^{n} X_i^p \geq nb \right)
=
\frac{1}{\pi(0)}\P_{\pi} \left( \sum_{i=0}^{T_1^*-1} (X_i^*)^p \geq nb, X_n^*=0 \right)
\nonumber
\\
&\geq
\frac{\pi(0)}{2}\P_{0}\left(\,\sum_{i=1}^{T_1}{X_i}^p > nb \right) -\bigO(e^{-cn})  
=
\frac{\pi(0)}{2}\P_{0}\left(W_1 > nb \right)-\bigO(e^{-cn}) .
\end{align}
We can now apply part \emph{(ii)} of Proposition~\ref{technical-proposition-for-tail-asymptotics} to $\eqref{eq:1st-lower-bound}$ and obtain the matching lower bound:
\begin{equation*}
    \liminf_{n \to \infty}\frac{1}{n^{\alpha}}\log\P_0\left(\frac1n V_n \geq b \right) \geq -\mathcal{B}_0^* \cdot b^{\alpha}.
\end{equation*}
\Halmos
\end{proof}

\subsection{Proof of Proposition~\ref{technical-proposition-for-tail-asymptotics}}
\label{subsec:proof-of-technical-proposition-for-tail-asymptotics}

We first state a number of preliminary results. These results are analytic in nature, and their proofs can be found in Sections \ref{pf-we-can-consider-only-upto-M} and \ref{pf-B_y_star-equal-V_y_T_star}. 
%\todo{Jose wrote: THE NOTATION ${\mathcal V_y^T}^*$ IS CONFUSING TO ME. WHY DO WE NEED A $*$? IT LOOKS LIKE IT IS ATTACHED TO $T$. PLUS THE UNUSUAL LABELING ADDS TO THE CONFUSION. Bert: I am currently ok with the notation. * is used for optimality, and the $T$ and $y$ also have a clear role here. }
For a fixed $M>0$, let $\notationdef{nota-B-y-AC-M}{B^{\AC;M}_{y}} \triangleq B^{\AC}_y \cap \{\xi \in \D[0,\infty): \T(\xi) \leq M\}$ and let $\notationdef{nota-B-y-M}{B_y^M} \triangleq B_y \cap \{\xi \in \D[0,\infty): \T(\xi) \leq M\}$.

\begin{lemma}
\label{we-can-consider-only-upto-M}
\linksinthm{we-can-consider-only-upto-M}
    For any given $y\geq 0$, there exists a constant  $M=M(y)>0$ such that
    \begin{itemize}
    \item[(i)]
    for each $\xi \in B^{\AC}_{y}$, there exists a path $\zeta\in B^{\AC;M}_{y}$ satisfying $I_y(\zeta) \leq I_y(\xi)$. 
    \item[(ii)]
    It holds that 
    \begin{equation}\label{eq:we-can-consider-only-upto-M}
   \inf_{\xi \in B_y^{\AC}}I_y(\xi)=\inf_{\xi \in B_y^{\AC;M}}I_y(\xi).
    \end{equation}
    \item[(iii)]
    Moreover, $M(y) \leq c y + d$ for some $c>0$ and $d>0$.
    \end{itemize}
\end{lemma}

\begin{lemma}
\label{B_y_star-equal-V_y_T_star} 
\linksinthm{B_y_star-equal-V_y_T_star} 
Let $M>0$ be the constant in Lemma~\ref{we-can-consider-only-upto-M}. Then,
$\B_y^* = {\V_y^T}^*$ 
for any $T\geq M$.
\end{lemma}

Set
$$
\notationdef{nota-K-t}{K_t}\triangleq \left\{\xi \in \D[0, t]: \xi(0)=0, \ \int_{0}^{t}(\Psi(\xi)(s))^p ds \geq 1, \ \xi(s) \geq 0 \ \text{for} \ s \in [0, t] \right\}.
$$
The following corollary is immediate from the previous lemma and Lemma~\ref{we-can-consider-only-upto-M}. 

\begin{corollary}
\label{corollary-to-merging-actions-to-the-left-lemma-and-B_y_star-equal-V_y_T_star}
\linksinthm{corollary-to-merging-actions-to-the-left-lemma-and-B_y_star-equal-V_y_T_star}
Let $M>0$ be the constant in Lemma~\ref{we-can-consider-only-upto-M}. For any $y\geq0$,
$$\inf_{t\in [0,M]}\inf_{\xi \in K_t}I_0^{\BV[0,t]}(\xi) 
= {\V_0^M}^* = \B_0^*.$$
\end{corollary}

Next, we formulate a key preparatory lemma. This lemma is motivated by a result of \cite{vysotsky}, stated as Result \ref{extensionmogulskildp} (ii) below. To apply this result, we need to verify a uniform continuity result. The next lemma provides the desired uniform continuity, whose proof is deferred to Section~\ref{pf-uniformcontinuity2}.

Recall the function $\Phi_T: \D[0,T] \rightarrow [0,\infty)$ defined as $\Phi_T(\xi) = \int_0^T \big(\Psi(\xi) (s)\big)^pds$ and $I_0^{\BV[0,T]}$ defined as 
	\begin{equation*}
	%\label{eq:ratefunctioninmrw}
	I_0^{\BV[0,T]}(\xi) \triangleq
	\begin{cases}
	\int_{0}^{T}\Lambda^*(\dot{\xi}^{(a)}(s))ds+ \theta_+(\xi^{(u)}(T))+\theta_-|\xi^{(d)}(T)| & \ \text{if} \ \xi \in \BV[0,T] \text{ and } \ \xi(0)=0,\\
	\infty & \ \text{otherwise}.
	\end{cases}
	\end{equation*}

\begin{lemma}
    \label{uniformcontinuity2}
    \linksinthm{uniformcontinuity2}
     For each $\gamma\geq 0$, $\Phi_T$ is uniformly continuous on the set $\{\xi: I_0^{\BV[0,T]}(\xi) \leq \gamma\}$ w.r.t. the $M_1'$ metric.
\end{lemma}

We apply this lemma in our next and final preparatory lemma: 

%\cmt{Somehow, convergence of the integrals is determined by continuity points, and there can be only countable many discontinuities. This might be exploited}
\begin{lemma}
\label{intermediate-lemma}
\linksinthm{intermediate-lemma}
\begin{itemize}
\item[(i)] For any $t,y\geq0$ and $T>0$,
\begin{equation}\label{upper-bound-for-large-hitting-time}
\limsup_{x \to \infty}\frac{1}{x}\log\P_{xy }\left(T_1/x >  T\right) 
\leq t y+T\log \E e^{tU}.
\end{equation}

\item[(ii)] For any $y \geq 0$ and $T>0$,
\begin{equation}\label{upper-bound-for-large-area-upto-T}
\limsup_{x \rightarrow \infty}\frac{1}{x} \log\P_{xy} \left( \int_{0}^{T} (X_{\lfloor sx \rfloor}/x)^p ds \geq 1  \right) 
\leq
-{\V_y^T}^*.
% \quad\mbox{as\ }x\to\infty, 
\end{equation}

\end{itemize}
\end{lemma}

\begin{proof}
{Proof of Lemma \ref{intermediate-lemma}}
\linksinpf{intermediate-lemma}
For part \emph{(i)}, note that
\begin{align*}
\P_{xy}\left(T_1  > xT\right)
&
\leq
\P_{x y}\left(X_{\lfloor xT\rfloor} >0 \right)
=\P\left(\sum_{i=1}^{\lfloor xT \rfloor}U_i >  -xy\right)
\leq
e^{t x y}\E\left(e^{tU}\right)^{\lfloor xT\rfloor},
\end{align*}
where the last inequality is from the Markov inequality.
Taking logarithms, dividing both sides by $x$, and taking $\limsup$, we get \eqref{upper-bound-for-large-hitting-time}.

For part \emph{(ii)}, 
%note that the reflection map $R$ is a Lipschitz continuous map from $\D[0,T]$  to $\D[0,T]$ w.r.t.\ the ${\cal M}_1$ topology (cf.\ \cite{whitt2002stochastic}, Theorem 13.5.1). Consequently, by 
note that conditional on $X_0 = xy$, $\int_{0}^{T} (X_{\lfloor sx \rfloor}/x)^p ds = \Phi_T(\bar K_x + y)$. From Lemma~\ref{uniformcontinuity}, we know that $\Phi_T$ is uniformly continuous over the sub-level sets of the rate fuction $I_y^{\BV[0,T]}$ of $\bar K_x + y$. Hence, we can apply 
%the functional $\Phi_T:\D[0,T] \mapsto \R_{+}$, defined as
%$
%\Phi_T(\xi) = \int_0^{T} (R(\xi)(s))^p ds,
%$
%is uniformly continuous on the set  $\{I_0^{\BV[0,T]}(\xi) \leq %\alpha\}$ under the $M_1$ topology. 
%Using 
Result \ref{extensionmogulskildp} (ii) to obtain
\[
\lim_{x \rightarrow \infty}\frac{1}{x} \log\P_{xy } \left( \int_{0}^{T} (X_{\lfloor sx \rfloor}/x)^p ds \geq 1  \right)
\leq
-\inf_{a \in [1, \infty)}J_{y}(a),
% \quad\mbox{as\ }x\to\infty,
\]
where $\notationdef{nota-J-y}{J_{y}(a)}\triangleq \inf\{I_0^{\BV[0,T]}(\xi): \xi \in \D[0,T], \xi(0)=y, \Phi_T(\xi)=a\}$. %, and $I_0^{\BV[0,T]}$ is the rate function associated with the LDP of $\bar K_n = \frac{1}{n}\textstyle{\sum_{i=1}^{n}U_i}$. 
Obviously, $\inf_{a \in [1,\infty)}J_{y}(a)={\V_y^T}^*$ and \eqref{upper-bound-for-large-area-upto-T} follows.
%\todo{$I_y$ was already defined as something else at the beginning of page 6.}
\eop
\end{proof}

%Our final preliminary result is an asymptotic lower bound for the area over a cycle: 

Now we are ready to prove Proposition~\ref{technical-proposition-for-tail-asymptotics}.
\begin{proof}{Proof of Proposition~\ref{technical-proposition-for-tail-asymptotics}.}
\linksinpf{technical-proposition-for-tail-asymptotics}
For part \emph{(i)},
consider a small enough $t_0>0$ so that $\E e^{t_0U} < 1$ which is possible since $\E U < 0$ and $U$ is light-tailed.
Then, thanks to Lemma~\ref{B_y_star-equal-V_y_T_star}, we can pick a sufficiently large $T>0$ so that 
$\B_{y}^* = {\V_{y}^T}^*$ and
$
t_0 y+T\log \E e^{t_0U} < -\B_{y}^*$.
Considering the case $T_1/x \leq  T$  and $T_1/x >  T$  separately and then applying the principle of the maximum term, 
\begin{align}
&
\limsup_{x\to\infty} \frac 1 x \log \P_{xy}\bigg(\int_{0}^{T_1/x} (X_{\lfloor ux \rfloor}/x)^p du \geq 1 \bigg)
\nonumber
\\
&
\leq
\limsup_{x\to\infty} \frac 1 x \log \Bigg\{\P_{xy}\bigg(\int_{0}^{T_1/x} (X_{\lfloor ux \rfloor}/x)^p du \geq 1,\, T_1/x \leq T \bigg)
+\P_{xy}\bigg(T_1/x > T \bigg)\Bigg\}
\nonumber
\\
&
\leq
\limsup_{x\to\infty} \frac1 x \log \P_{xy} \left( \int_{0}^{T} (X_{\lfloor ux \rfloor}/x)^p du \geq 1\right) 
\vee
\limsup_{x\to\infty} \frac1 x \log \P_{xy}\left(T_1/x >  T  \right) 
\nonumber
\\
&
\leq
\Big(-{\V_{ y}^T}^*\Big) \vee \Big(t_0  y+T\log \E e^{t_0U}\Big) 
=
\Big(-\B_{y}^*\Big) \vee \Big(t_0 y+T\log \E e^{t_0U}\Big) 
=
-\B_{y}^*,
\label{some-intermediate-inequality-2}
\end{align}
where we used Lemma~\ref{intermediate-lemma} for the third inequality.

Next, we move on to part \emph{(ii)}.
For any given $t>0$, let
\begin{itemize}
	\item $\notationdef{nota-A-t-epsilon}{A_{t,\epsilon}}\triangleq\{\xi \in \D[0,t]: \xi(0)=\epsilon, \  \int_{0}^{t}\Psi(\xi)(s)^p ds > 1, \, \xi(s) >0,\ \forall \, s \, \in [0,t]\}$ and
	\item $\notationdef{nota-tilde-A-t-epsilon}{\tilde A_{t,\epsilon}}\triangleq\{\xi \in \D[0,t]: \xi(0)=\epsilon, \ \int_{0}^{t}\Psi(\xi)(s)^p ds > 1, \, \xi(s) >\epsilon/2, \, \forall s \, \in [0,t]\}$.
\end{itemize}
Set $u=x^{1+p}$.
Let $\epsilon$ be small enough such that $\P(U_1>\sqrt{\epsilon})>0$. 
Define the event 
$
B_{x,\epsilon} \triangleq \{U_i > \sqrt{\epsilon}, 
i=1,\ldots,\lceil x\sqrt{\epsilon} \rceil \}.
$
Setting $k^* = \lceil x\sqrt{\epsilon} \rceil+1$, we obtain
\begin{align*}
&\liminf_{u \to \infty}\frac{1}{u^{\alpha}}\log\P_0\left(W_1>u\right)
\\
&
=
\liminf_{x \to \infty}\frac{1}{x}\log\P_0\left(\sum_{k=0}^{T_1}X_k^p> u \right)
\\
&
\geq 
\liminf_{x \to \infty}\frac{1}{x}\log
\P_0\left(\sum_{k=k^*}^{T_1}X_k^p > x^{1+p}, \, 
B_{x,\epsilon} \right)
\\
&
=
\liminf_{x \to \infty}\frac{1}{x}\log
\left[
\P_0\left(\sum_{k=k^*}^{T_1}X_k^p > x^{1+p} \big| 
B_{x,\epsilon}\right)\P_0\left(B_{x,\epsilon}\right)	   
\right]
\\
&
\geq 
\liminf_{x \to \infty}\frac{1}{x}\log 
\left[
\P_{\epsilon x}\left(\sum_{k=0}^{T_1}X_k^p > x^{1+p} \right)\P_0\left(B_{x,\epsilon}\right) 
\right]
\\
&
=
\liminf_{x \to \infty}\frac{1}{x}\log\left[
\P_{\epsilon x}\left(\int_{0}^{T_1/x}\left(X_{\lfloor sx \rfloor}/x\right)^p ds >1 \right)\P_0\left(
B_{x,\epsilon}\right)
\right]
\\
&
\geq
\liminf_{x \to \infty}\frac{1}{x}\log\left[\P_{\epsilon x}\left(\int_{0}^{t}\left(X_{\lfloor sx \rfloor}/x\right)^p ds>1, \, T_1 > xt \right)\P_0\left(B_{x,\epsilon}\right)\right]
\\
&
=
\liminf_{x \to \infty}\frac{1}{x}\log\left[\P_{\epsilon x}\left(\int_{0}^{t}\left(X_{\lfloor sx \rfloor}/x\right)^p ds>1, \, X_{\lfloor sx \rfloor}/x>0, \,\forall s \in [0,t] \right)\P_0\left(B_{x,\epsilon}\right)\right]
\\
&
\geq 
\liminf_{x \to \infty}\frac{1}{x}\log \left[\P_{\epsilon}\left(\bar{K}_x \in A_{t,\epsilon}\right) \P_0\left(B_{x,\epsilon}\right) \right]
\\
&
\geq 
-\inf_{\xi \in (A_{t,\epsilon})^{\circ}}I_{\epsilon}^{\BV[0,t]}(\xi)
+\sqrt{\epsilon} \log \P(U_1>\sqrt{\epsilon})
\\
&
\geq 
-\inf_{\xi \in \tilde A_{t,\epsilon}}I_{\epsilon}^{\BV[0,t]}(\xi)+\sqrt{\epsilon} \log \P(U_1>\sqrt{\epsilon}),
\end{align*}
where the third equality is from part \emph{(i)} of Proposition~\ref{miscellaneous-proposition-for-tail-asymptotics}.
The second to last inequality follows from part \emph{(i)} of Result \ref{extensionmogulskildp} since the integral and the infimum are both continuous in the ${M}_1$ topology (see, respectively Theorem 11.5.1 and Theorem 13.4.1 of \cite{whitt2002stochastic}). 
Recall that
$$
K_t=\left\{\xi \in \D[0, t]: \xi(0)=0, \ \int_{0}^{t}(\Psi(\xi)(s))^p ds \geq 1, \ \xi(s) \geq 0 \ \text{for} \ s \in [0, t] \right\}.
$$
Note that for all $\epsilon>0$,
\begin{equation}\label{equal-open-toG0}
\inf_{\xi \in \tilde A_{t,\epsilon}}I_\epsilon^{\BV[0,t]}(\xi) \leq \inf_{\xi \in K_t}I_0^{\BV[0,t]}(\xi).
\end{equation}
%To show that $\inf_{\xi \in K}I_y^{\BV[0,t]}(\xi) \geq \inf_{\xi \in B_{t_1,\epsilon}}I_0^{\BV[0,T]}^{t_1}(\xi)$, 
To see this, suppose that $\xi \in K_t$. Then, $\tilde{\xi}=\epsilon+\xi$ belongs to $\tilde A_{t,\epsilon}$ and $I_\epsilon^{\BV[0,t]}(\tilde \xi)=I_0^{\BV[0,t]}(\xi)$. Since the construction holds for every $\xi \in K_t$, we have that $\inf_{\xi \in K_t}I_0^{\BV[0,t]}(\xi) \geq  \inf_{\xi \in \tilde A_{t,\epsilon}}I_\epsilon^{\BV[0,t]}(\xi)$. 
Therefore,
$$
\liminf_{u \to \infty}\frac{1}{u^{\alpha}}\log\P_0\left(W_1>u\right)
\geq 
-\inf_{\xi \in K_t}I_0^{\BV[0,t]}(\xi) + \sqrt \epsilon \log \P(U_1 > \sqrt \epsilon).
$$
Since $\epsilon$ and $t$ are arbitrary, taking $\epsilon \to 0$ and taking the infimum over $t\in[0,M]$, Corollary~\ref{corollary-to-merging-actions-to-the-left-lemma-and-B_y_star-equal-V_y_T_star} gives
$$
\liminf_{u \to \infty}\frac{1}{u^{\alpha}}\log\P_0\left(W_1>u\right)
\geq 
-\inf_{t\in [0,M]}\inf_{\xi \in K_t}I_0^{\BV[0,t]}(\xi) 
= -\B_0^*.
$$
\eop

\end{proof}

\subsection{Proof of Proposition~\ref{miscellaneous-proposition-for-tail-asymptotics}}
\label{subsec:proof-of-miscellaneous-proposition-for-tail-asymptotics}

%We start with the proof of Proposition~\ref{miscellaneous-proposition-for-tail-asymptotics}.
%\begin{proof}{Proof of Proposition~\ref{miscellaneous-proposition-for-tail-asymptotics}.}
\linksinpf{miscellaneous-proposition-for-tail-asymptotics}
For part \emph{(i)}, note that 
$$\frac{1}{x^{1+p}}\sum_{k=0}^{m-1} X_k^p 
=\frac{1}{x^{1+p}} \int_0^{m} X_{\lfloor u \rfloor}^p du 
= \frac{1}{x^{1+p}}\int_0^{m/x}xX_{\lfloor xs \rfloor}^p ds
= \int_0^{m/x}\left(\frac{X_{\lfloor xs \rfloor}}{x}\right)^p ds,$$ 
where the second equality is from the change of variable $u= xs$. 
The claimed equivalence is immediate from this.

For part \emph{(iii)}, note that 
$$\lim_{k\to \infty} \min_{i\geq 1} \Bigg\{ 
\frac{i-1}k \beta\bar y + \B_{\frac i k \bar y}^*
\Bigg\}
=
\lim_{k\to \infty} \left(\min_{i\geq 1} \Bigg\{ 
\beta\frac{i}k \bar y + \B_{\frac i k \bar y}^*
\Bigg\} - \frac1k \beta \bar y\right)
=
\lim_{k\to \infty} \min_{i\geq 1} \Bigg\{ 
\beta\frac{i}k \bar y + \B_{\frac i k \bar y}^*
\Bigg\}.
$$
Moreover, from part \emph{(ii)} of Proposition~\ref{eq:ineqformonotonicity}, 
$$
\lim_{k\to \infty} \min_{i\geq 1} \Bigg\{ 
\frac{i}k \beta\bar y + \B_{\frac i k \bar y}^*
\Bigg\}
= \inf_{y\in [0,\infty)} \left\{\beta y + \mathcal B_y^*\right\}.
$$

For part \emph{(ii)}, note that by definition, ${\mathcal{B}_0^{*}} \geq {\mathcal{B}_{\pi}^{*}}$. Therefore, we only have to prove that ${\mathcal{B}_0^{\ast}} \leq {\mathcal{B}_{\pi}^{\ast}}$. Recall that $\beta=\sup\{\theta > 0:\E(e^{\theta U}) \leq 1\}$ and $\theta_+=\sup\{\theta \in \R: \E(e^{\theta U})< \infty\}$. For the rest of this proof, let $\Lambda$ be the log-moment generating function and let \notationdef{nota-D-Lambda}{$D_{\Lambda}$} denote the effective domain of $\Lambda$, i.e.\ $D_{\Lambda}=\{x: \Lambda(x)< \infty\}$. 
%That is, $\Lambda^*$ is the convex conjugate of $\Lambda$. 
We start with a claim: 
for any $\epsilon>0$ there exists 
a $u>0$ such that 
\begin{equation}\label{claim1}
\Lambda^*(u)/u \leq \beta+\epsilon.
\end{equation}
To prove (\ref{claim1}) we distinguish between the cases $\beta < \theta_+$ and $\beta = \theta_+$.
For the first case note that $\beta \in D_{\Lambda}^{\circ}$. In view of the convexity and continuity of $\E(e^{\theta U})$,  $\E(e^{\beta U})=1$. Due to Lemma~2.2.5 $(c)$ in \cite{dembo2010large}, $\Lambda$ is a differentiable function in $D_{\Lambda}^{\circ}$ with $\Lambda'(\eta) =\frac{\E(Ue^{\eta U})}{\E(e^{\eta U})}$. Since $\beta \in D_{\Lambda}^{\circ}$ we have that $\Lambda'(\beta)=\E(Ue^{\beta U})< \infty$. In addition, $\Lambda^{'}(0)=\E(U)<0$ implies that $\Lambda(\eta)$ is decreasing for small values of $\eta$. Now, the strict convexity and differentiability of $\Lambda$ over its effective domain implies that $\Lambda'$ is increasing at $\beta$ and thus $\E(Ue^{\beta U})>0$. It can be checked that, for $u=\E(Ue^{\beta U})$,
\[
\frac{\Lambda^*(u)}{u}=\frac{\beta \E(Ue^{\beta U})-\log\E(e^{\beta U})}{\E(Ue^{\beta U})}=\beta,
\]
and hence our claim is proved. 
Consider now the case $\beta = \theta_+$. In view of Equation (5.5) in \cite{mogulskii1993large}, $\lim_{x \to \infty}\frac{\Lambda^*(x)}{x} = \theta_+$. That is, for any $\epsilon >0$ we can choose a $u$ so that $\Lambda^*(u)/u \leq \theta_+ + \epsilon= \beta +\epsilon$. We proved the claim \eqref{claim1}.

Back to the inequality ${\mathcal{B}_0^{\ast}} \leq {\mathcal{B}_{\pi}^{\ast}} $, we will show that for any given $\epsilon>0$ and any given path $\xi \in B_{y}$, we can construct a path $\zeta \in B_0$ so that $I_0(\zeta) \leq I_{y}(\xi)+\beta y + \epsilon$. 
To this end, let $u>0$ be such that $\Lambda^*(u)/u \leq \beta + \epsilon/y$ and set
\[
\zeta(s)\triangleq u s \mathbbm{1}_{\{s \leq y/u\}}+\xi(s-y/u)\mathbbm{1}_{\{s > y/u\}}. 
\]
Then $\zeta(0)=0$, $\zeta(y/u)=y$, and $\zeta \in B_0$. 
Also, one can see that $I_y(\zeta)=(y/u)\Lambda^*(u)+\int_{0}^{\mathcal{T}(\xi)}\Lambda^*(\dot{\xi}(s))ds+\theta_+ \xi^{(u)}(\T(\xi))=(y/u)\Lambda^*(u)+ I_y(\xi)$. From the construction of $u$, 
\[
I_{0}(\zeta) \leq \beta y +\epsilon+I_{y}(\xi)
\]
as desired. This concludes the proof of part \emph{(iv)}. \eop

\section{Additional technical lemmas and proofs}
\label{additionalproofsection}

%In this section we provide the proofs of Lemma's and Propositions FILLIN, which were stated 
%in Section \ref{sec:firstproofsection}. These results are in turn based on several additional technical lemma's which are stated and proved here as well. 

\subsection{Proof of Lemma \ref{useful-properties-of-variational-By}}
\label{pf-useful-properties-of-variational-By}
%{Preliminary results on the variational problem $\mathcal{B}_y^*$}\label{preliminary-results-By}

We first state a lemma on the monotonicity of the reflection map which is useful in the proof of Lemma \ref{useful-properties-of-variational-By}.
%, which is needed to give a proof 
%of Lemma \ref{useful-properties-of-variational-By}.

\begin{lemma}
\label{monotonicity-of-reflection}
\linksinthm{monotonicity-of-reflection}
Suppose that $\alpha, \beta, \gamma \in \D[0,T]$, $\alpha(s) = \beta(s) + \gamma(s)$, and $\gamma(s)$ is non-negative and non-decreasing. 
Then, $\Psi(\alpha)(t) \geq \Psi(\beta)(t)$ for all $t\in[0,T]$.
\end{lemma}
\begin{proof}{Proof.}
\linksinpf{monotonicity-of-reflection}
The proof of this lemma is an immediate consequence of Theorem 14.2.2 of \cite{ramasubramanian2000subsidy}.
\eop
\end{proof}
% \begin{proof}{Proof.}
% Recall first that if $z\geq0$ then $x \wedge (y+z) \leq (x\wedge y) + z$ for any $x, y\in \R$.
% From the non-negativity and monotonicity assumptions on $\gamma$, we have that
% \begin{equation*}
% 0 \wedge \alpha(s) \leq 0 \wedge \beta(s) + \gamma(s) \leq 0 \wedge \beta(s) + \gamma(t),
% \quad0 \leq s \leq t
% \end{equation*}
% and hence,
% $$
% 0 \wedge \inf_{0\leq s\leq t} \alpha(s) \leq 0 \wedge \alpha(s) \leq 0 \wedge \beta(s) + \gamma(t),
% \quad0 \leq s \leq t.
% $$
% Taking infimum over $s\in[0,t]$, we get
% $
% 0\wedge \inf_{s\in[0,t]}\alpha (s) \leq 0 \wedge \inf_{s\in[0,t]} \beta (s) + \gamma(t).
% $
% Therefore, 
% $$
% \Psi(\alpha)(t) = \alpha(t) - 0\wedge \inf_{s\in[0,t]}\alpha (s) 
% \geq 
% \alpha(t) - 0\wedge\inf_{s\in[0,t]} \beta(s) - \gamma(t)
% = \beta(t) - 0\wedge \inf_{s\in[0,t]} \beta(s) = \Psi(\beta)(t).
% $$
% \eop
% \end{proof}

\medskip

\proof{Proof of Lemma \ref{useful-properties-of-variational-By}.}
\linksinpf{useful-properties-of-variational-By}
For \textit{(i)}, consider $\xi^*(t) \triangleq  y + \mu \cdot t$ and note that $I_y(\xi^*) = 0$. Also, $y\geq \bar y$ implies that $\xi^* \in B_{ y}$, and hence, $\B_{ y}^* = 0$. 
For part \textit{(ii)}, we first show that  $\mathcal B_y^*=\inf_{\xi \in B_y \cap J_+}I_y(\xi)$, where $\notationdef{nota-J-+}{J_+}\triangleq \{\xi \in \D[0,\infty): \xi^{(d)}\equiv0\}$.
 To prove this, let $\xi \in B_y$ be given arbitrarily. 
 We  construct a path $\xi_2$ so that $\xi_2 \in B_y \cap J_+$ and $I_y(\xi_2) \leq I_y(\xi)$.  Towards this, discard the downward jumps of $\xi$ and let $\xi_1\triangleq \xi^{(a)}+\xi^{(u)}$. 
 Then, set $\xi_2(t) \triangleq \xi_1(t \wedge \mathcal{T}(\xi)) + \mu([t-\mathcal{T}(\xi)]^+), \ t \geq 0 $. In view of Lemma~\ref{monotonicity-of-reflection}, $\Psi(\xi_2) \geq \Psi(\xi)$ over $[0,\mathcal{T}(\xi)]$, and hence, $\xi_2 \in B_y$  and $I_y(\xi) = I_y(\xi_1)+\theta_- \xi^{(d)}(\mathcal{T}(\xi))= I_y(\xi_2)+\theta_- \xi^{(d)}(\mathcal{T}(\xi)) \geq I_{y}(\xi_2)$. 
 This proves $\mathcal B_y^*=\inf_{\xi \in B_y \cap J_+}I_y(\xi)$.

Now, suppose that $\xi \in B_y\cap J_+$ is given. It is sufficient to prove that there exists a  constant $c_\xi \in (0,\infty)$ such that for any given $\epsilon>0$, one can construct a $\tilde{\xi}_\epsilon \in B_y^{\mathbbm{AC}}$ such that 
\begin{equation}
I_y(\tilde{\xi}_\epsilon) \leq I_y(\xi) + \epsilon \cdot c_\xi.
\label{target-inequality-of-lemma-3-2}
\end{equation}
In case $y \geq \bar y$, simply consider $\tilde{\xi}_\epsilon(t)\triangleq y+ \mu \cdot t$ regardless of $\xi$ and $\epsilon$. 
Then, $\tilde{\xi}_\epsilon \in B_y^{\AC}$ and $I_y(\tilde{\xi}_\epsilon)=0 \leq I_y(\xi) + \epsilon \cdot c_\xi$ with $c_\xi = 0$. 
In case $y < \bar y$, we first prove the following claim w.r.t. the local rate function $\Lambda^*$: for any $\epsilon>0$ there exists a $y_\epsilon \in \R_+$ such that
\begin{equation*}
\Lambda^*(y) \leq (\theta_+ + \epsilon)y, \quad \forall y \geq y_\epsilon.
%\label{claim-in-lem-3-2}
\end{equation*}
  %$$\liminf_{y \to \infty} \frac{-\theta y + \Lambda(\theta)}{y} \geq -\theta_+.$$
Since this is trivially true if $\theta_+ = \infty$, we prove the claim with the assumption that $\theta_+ < \infty$.
Note that, due to Jensen's inequality and the monotonicity of the map $x \mapsto \log x$, we see that $\Lambda(\theta)=\log(\E(e^{\theta U_1})) \geq \log e^{\theta \E(U_1)} = \mu\cdot\theta $. 
This implies, for any $\theta < \theta_+$, 
$$y \theta  - \Lambda(\theta) \leq y \theta  -  \mu \theta  \leq y\theta_+ -  \mu \theta_+ \leq (\theta_++\epsilon) y, \quad \forall y \geq -\mu\theta_+/\epsilon.$$ 
Therefore, 
$$\Lambda^* (y) = \sup_{\theta \leq \theta_+} \{y \theta - \Lambda(\theta)\} \leq (\theta_+ + \epsilon) y, \quad \forall y \geq -\mu \theta_+/\epsilon,$$
which is the claim with $y_\epsilon = -\mu \theta_+/\epsilon$.

To conclude the proof, pick  $\xi \in B_y \cap J_+$  and assume w.l.o.g.\ that $I_y(\xi)< \infty$ (since \eqref{target-inequality-of-lemma-3-2} is trivial otherwise). For any given $\epsilon>0$, we construct a path $\tilde{\xi}_\epsilon$ such that $\tilde{\xi}_\epsilon \in B_y^{\AC}$ and $I_y(\tilde{\xi}_\epsilon) \leq I_y(\xi) + \epsilon\cdot \xi^{(u)}(\cal T(\xi))$, which is \eqref{target-inequality-of-lemma-3-2} with $c_\xi = \xi^{(u)}(\cal T(\xi))$. 
Set $\tilde{\xi}_\epsilon(t)=y +y_\epsilon t$ for $t \in [0,\xi^{(u)}(\mathcal{T}(\xi))/y_\epsilon]$, and $\tilde{\xi}_\epsilon(t)= \xi^{(a)}(t-\xi^{(u)}(\mathcal{T}(\xi))/y_\epsilon)+ \xi^{(u)}(\mathcal{T}(\xi)) $ for $t > \xi^{(u)}(\mathcal{T}(\xi))/y_\epsilon$, with $y_\epsilon = -\mu \theta_+/\epsilon$. % $y_\epsilon$ is as in the above claim. BERT: "above" was not very clear  so I repeated the specific value
Then it is straightforward to check that $\tilde{\xi}_\epsilon \in B_y^{\mathbbm{AC}}$ from its construction and Lemma~\ref{monotonicity-of-reflection}. 
In addition,  $\tilde{\xi}_\epsilon\big( \mathcal{T}(\xi)+\xi^{(u)}(\mathcal{T}(\xi))/y_\epsilon \big)=\xi^{(a)}\big(\mathcal{T}(\xi)\big)+ \xi^{(u)}(\mathcal{T}(\xi))=0$ since $\xi \in B_y \cap J_+$, and hence, $\mathcal{T}(\tilde{\xi}_\epsilon) \leq \mathcal{T}(\xi)+\xi^{(u)}(\mathcal{T}(\xi))/y_\epsilon$. Consequently, 
  \begin{align*}
  I_y(\tilde{\xi}_\epsilon) & =\int_{0}^{\mathcal{T}(\tilde{\xi}_\epsilon)}\Lambda^*\Big(\dot{\tilde{\xi}}_\epsilon(s)\Big)ds
  =
  \frac{\xi^{(u)}(\mathcal{T}(\xi))}{y_\epsilon}\Lambda^*(y_\epsilon)+\int_{\frac{\xi^{(u)}(\mathcal{T}(\xi))}{y_\epsilon}}^{\mathcal{T}({\xi})+\frac{\xi^{(u)}(\mathcal{T}(\xi))}{y_\epsilon}} \Lambda^*\big(\dot\xi^{(a)}(s-\xi^{(u)}(\mathcal{T}(\xi))/y_\epsilon)\big)ds
  \\
  &
  \leq
  \xi^{(u)}(\mathcal{T}(\xi))(\theta_++\epsilon)+\int_{0}^{\mathcal{T}({\xi})}\Lambda^*\big(\dot\xi^{(a)}(s)\big)ds
  \leq I_y(\xi) + \epsilon \cdot \xi^{(u)}(\mathcal{T}(\xi)).
  \end{align*}
  We arrived at the desired inequality. 
  \eop
\endproof

\subsection{Proof of Proposition \ref{eq:ineqformonotonicity}}
\label{pf-eq:ineqformonotonicity}
\linksinpf{eq:ineqformonotonicity}
For part \emph{(i)}, let $x, y$ be such that $0 \leq x < y$. 
We will show that for any $\epsilon>0$, there exists $\zeta \in B_y^{\AC}$ such that $I_y(\zeta) < \B_x^* + \epsilon$.
To show this, note from Lemma~\ref{useful-properties-of-variational-By} and Lemma~\ref{we-can-consider-only-upto-M} that there exists $\xi \in B_x^{\AC}$ such that $I_x(\xi) < \B_x^* + \epsilon$ and $\T(\xi) < \infty$.
Set 
$$
\zeta(t) \triangleq 
(y-x) + \xi(t \wedge \T(\xi)) + \mu\cdot[t-\T(\xi)]^+.
$$
Since $\zeta(0) = y$, and $\zeta(t) - \xi(t)$ is non-negative and non-decreasing on $t\in[0,\T(\xi)]$, Lemma~\ref{monotonicity-of-reflection} implies that $\Psi(\zeta)(t) > \Psi(\xi)(t)$ on $t\in[0,\T(\xi)]$, and hence,  $\zeta \in B_y^{\AC}$. Moreover, note that $\T(\xi) \leq \T(\zeta)$, and $\Lambda^*(\dot \zeta(s)) = \Lambda^*(\mu) = 0$ on $s\in[\T(\xi), \T(\zeta)]$. Therefore, 
\begin{align*}
 I_y(\zeta) 
    = \int_0^{\T(\zeta)}\Lambda^*(\dot{\zeta}(s))ds 
    =\int_0^{\T(\xi)}\Lambda^*(\dot{\zeta}(s))ds 
    = \int_0^{\T(\xi)}\Lambda^*(\dot{\xi}(s))ds
     = I_x(\xi) < \B_x^* + \epsilon.
\end{align*}
For part \emph{(ii)}, note that we only need to prove one side of the inequality thanks to part \emph{(i)}. That is, it is enough to show that if $0\leq x < y$, then $\B_x^* \leq \B_y^* + (y-x)\Lambda^*(1)$.
Fix an $\epsilon>0$ and pick $\zeta\in B_y^{\AC}$ such that $I_y(\zeta) \leq  \B_y^* + \epsilon$, which is always possible due to part \emph{(ii)} of Lemma~\ref{useful-properties-of-variational-By}.
Set
$$
\xi(t)
\triangleq 
(x+t)\one_{[0,y-x]}(t) + \zeta(t-(y-x))\one_{[y-x,\infty)}(t).
$$
It follows that $\T(\xi) = \T(\zeta) + y-x$, and hence,
\begin{align*}
I_x(\xi) 
&= 
\int_0^{y-x} \Lambda^*(1)ds + \int_{y-x}^{\T(\xi)} \Lambda^*(\dot\xi(s))ds= 
(y-x) \Lambda^*(1) + \int_{y-x}^{\T(\zeta)+ y-x} \Lambda^*(\dot\xi(s))ds
\\
&=
(y-x) \Lambda^*(1) + \int_{0}^{\T(\zeta)} \Lambda^*(\dot\xi(s+(y-x)))ds
=
(y-x) \Lambda^*(1) + \int_{0}^{\T(\zeta)} \Lambda^*(\dot\zeta(s))ds
\\
&=
(y-x) \Lambda^*(1) + I_y(\zeta)
\leq 
(y-x) \Lambda^*(1) + \B_y^* + \epsilon.
\end{align*}
Since $\xi \in B_x$,  $\B_x^* \leq (y-x) \Lambda^*(1) + \B_y^* + \epsilon$.
Taking $\epsilon \to 0$ yields \emph{(ii)}.\\
For part \emph{(iii)}, we first note that $\mathcal B_y^* < \infty$: if we set $\xi(t) \triangleq (|\mu|(p+1))^{\alpha}(1+\mu(t-1))\I_{[1,\infty)}(t)$, then $\xi \in B_y^{\mathbbm{AC}}$ and $I_y(\xi) < \infty$. 
% we use the representation in \eqref{partitioned infimum}; that is,
%   $
% 	\mathcal{B}_y^{\ast} = \inf_{\xi \in B_y^{\CON}}I_{y}(\xi) = \inf_{z_0 \geq z_T} \inf_{T\geq 0} \inf_{\xi \in F^{y}_{z_0, z_T, T}} I_{y}(\xi).  
%   $
% Moreover, due to Lemma~\ref{F-z-T-ccompact}, $F^{y}_{z_0, z_T, T}=\{ \xi: \xi \in B_y^{\CON}, \dot \xi(0)=z, \dot \xi(T)=z_T, \ \xi(T)=0 \}$ is a compact set w.r.t. the $J_1$ topology.
To show that $\mathcal B_y^*>0$, 
%we recall some  useful properties of $\Lambda^*$:
%\begin{itemize}
%\item[A)] $\Lambda^*$ is nonnegative, 
%\item[B)] $\Lambda^*(\mu)=0$, and
%\item[C)] $\Lambda^*$ is convex and nonnegative, decreasing %on $(-\infty,\mu)$, and it is increasing on $(\mu,\infty)$.
%\end{itemize}
%\cmt{BZ: what follows here is a more direct approach, please %check it.}
we can take $y\in (0,\bar y)$ without loss of generality in view of (i). 
Define $\bar y (t) = \bar y + \mu t$ and set
 $T_\mu = \T(\bar y(\cdot))= \bar y/|\mu|$. Write 
\begin{equation}
\label{min-of-inf}
    \mathcal{B}_y^*= \min\{ \inf_{t\in [0, T_\mu]}\inf_{\xi \in B_y^{\mathbbm{AC}}, \T(\xi)=t} I_y(\xi), \inf_{t\in[T_\mu,\infty]}\inf_{\xi \in B_y^{\mathbbm{AC}}, \T(\xi) = t } I_y(\xi) \}.
\end{equation}
We show that each of the two infima in this minimum is strictly positive. Observe that
\begin{align*}
\inf_{t\in[T_\mu,\infty]}\inf_{\xi \in B_y^{\mathbbm{AC}}, \T(\xi) =t } I_y(\xi) 
%&\geq \inf_{\xi: \xi(0)=y, \xi \in \mathbbm{AC}, \T(\xi) \geq T_\mu } \int_0^{T_\mu} \Lambda^*(\dot \xi (s)) ds\\
&\geq\inf_{t\in[T_\mu,\infty]}  \inf_{\xi: \xi(0)=y, \xi \in \mathbbm{AC}, \T(\xi) =t } \int_0^{t} \Lambda^*(\dot \xi (s)) ds\\
&\geq \inf_{t\in[T_\mu,\infty]} t \Lambda^*(-y/t) \geq \frac{\bar y}{|\mu |} \Lambda^* \left(\mu \frac{y}{\bar y} \right)>0,
%    &= \inf_{s \leq |\mu | / \bar y} \frac 1s \sup_{\theta} \theta ys - \Lambda(\theta)
\end{align*}
using that $\Lambda^*$ is a convex non-negative function 
\BZA{should we add a reference or not?}%
\CRA{I think this is standard enough that we don't need to provide a reference.}%
with $\Lambda^*(\mu)=0$, and Jensen's inequality. 

To lower bound the double infimum in (\ref{min-of-inf}), observe that, 
since $y<\bar y$, the area constraint $\int_0^t \xi (s)^p ds\geq 1$ can only be valid 
if there exists an $s\in [0,t]$ such that $\xi(s) \geq \bar y + \mu s$. Consequently, 
\begin{align*}
    \inf_{\xi \in B_y^{\mathbbm{AC}}, \T(\xi)=t} I_y(\xi) &\geq 
    \inf_{s\leq t} \inf_{\xi \in \mathbbm{AC}, \xi(0)=y, \xi(s) \geq \bar y + \mu s} I_y(\xi)\\
    &\geq \inf_{s\leq t}  \inf_{\xi \in \mathbbm{AC}, \xi(0)=y, \xi(s) \geq \bar y + \mu s} \int_0^{s} \Lambda^*(\dot \xi (u)) du\\
    &= \inf_{s\leq t}  \inf_{\xi \in \mathbbm{AC}, \xi(s)-\xi(0) \geq \bar y - y + \mu s} \int_0^{s} \Lambda^*(\dot \xi (u)) du\\
    &\geq \inf_{s\leq t} s \Lambda^*\left(\frac{\bar y-y}{s}+\mu\right),
\end{align*}
where we applied Jensen's inequality in the last step. Consequently, 
\[
    \inf_{t\in [0, T_\mu]}\inf_{\xi \in B_y^{\mathbbm{AC}}, \T(\xi)=t} I_y(\xi) \geq 
    \inf_{s\leq T_{\mu}} s \Lambda^*\left(\frac{\bar y-y}{s}+\mu\right)= \frac{\bar y}{|\mu |} \Lambda^*\left( \mu + |\mu| \frac{\bar y-y}{\bar y}\right)>0.
\]
The equality holds since $\Lambda^*$ is a strictly convex non-negative function with $\Lambda^*(\mu)=0$ . \eop \\

\subsection{Proof of Lemma \ref{merging-actions-to-the-left-lemma}}
\label{pf-merging-actions-to-the-left-lemma}
\linksinpf{merging-actions-to-the-left-lemma}
For part \emph{(i)}, we first construct a new trajectory $\xi_1$ from $\xi$ by discarding the downward jumps, i.e., $\xi_1=\xi^{(a)}+\xi^{(u)}$. 
Obviously, $I_y^\BVT(\xi_1) \leq I_y^\BVT(\xi)$. 
Note that $\xi_1 = \xi + (-\xi^{(d)})$, where $-\xi^{(d)}$ is non-negative and non-decreasing. From Lemma~\ref{monotonicity-of-reflection} we have that $\Psi(\xi_1)(t) \geq \Psi(\xi)(t)$ for all $t\in[0,T]$, and hence,  $\Phi_T(\xi_1) \geq \Phi_T(\xi)$. 
For each $t \in [0,T]$, let $l(t) \triangleq \inf\{s\in[0,T]:\, \Psi(\xi)(u) > 0 \text{ for all } u \in [s,t]\}$,  $r(t) \triangleq \sup\{s\in[0,T]:\, \Psi(\xi)(u) > 0 \text{ for all } u \in [t, s]\}$, and $\sigma(t) \triangleq [l(t), r(t))$. 
Set $\mathcal C^+_1 \triangleq \{ \sigma(t)\subseteq [0,T]:  t\in [0, T]\}$.
%Note that from the right continuity of $\xi$, $[l(t), r(t))$ is non-empty for any $t$ such that $\xi(t) > 0$ and $t < T$.
Note that, by construction, the elements of $\mathcal C^+_1$ cannot overlap, and hence, there can be at most countable number of elements in $\mathcal C^+_1$.
In view of this, we write $\mathcal C^+_1 = \{ [l_i, r_i):\ i\in \N\}$ and let $\sigma_i \triangleq [l_i, r_i)$. 
The following observations are immediate from the construction of $\mathcal C^+_1$, the right continuity of $\xi$, and the fact that $\xi_1$ does not have any downward jumps.
\begin{itemize}
\item[O1.]
If $t\in[0,T)$ does not belong to any of the elements of $\mathcal C^+_1$, then $\Psi(\xi_1)(t) = 0$.
\item[O2.]
$\Psi(\xi_1)$ is continuous on the right end of the intervals $\sigma_i$ except for the case $r_i = T$.
\end{itemize}
%\begin{itemize}
%\item[O3.]
Note that O1 also implies that $\xi_1(t)=\xi_1(t-)$ for such $t$'s.
%\end{itemize}
%Also, we note that 
%Denote the left end-point of $\sigma_i$ with $l_i$ and the right end-point with $r_i$.
Let $s_n \triangleq \sum_{i=1}^{n-1} (r_i - l_i)$ for $n\in\N$.
Note that $s_n \to s_\infty \in [0, T]$ as $n\to\infty $.
Let $\dot\xi^{(a)}(t)$ denote the time derivative  $\frac d {dt}\xi^{(a)}(t)$ of $\xi^{(a)}$ at $t$, and set
$$\zeta_1(t) \triangleq y + \int_0^t \dot\zeta_1(s) ds +   \zeta_1^{(u)}(t),$$
where
$$\dot\zeta_1(t) \triangleq \sum_{i\in \N} \dot\xi^{(a)}(t-s_i+l_i)\ind_{[s_i, s_{i+1})}(t) + \mu \ind_{[s_\infty, T]}(t),$$
and
$$\zeta_1^{(u)}(t) \triangleq \sum_{i\in \N} \big(\xi^{(u)}(t\wedge s_{i+1}- s_i+l_i)- \xi^{(u)}(l_i-)\big)\ind_{[s_i, T]}(t).$$
That is, on the interval $[s_i, s_{i+1})$, $\zeta_1$ behaves the same way as $\xi_1$ does on the interval $[l_i, r_i)$, whereas $\zeta_1$
decreases linearly at the rate $|\mu|$ outside of those intervals.
Given this, it can be checked that
\begin{itemize}
\item[O3.]
$\int_{s_i}^{s_{i+1}}\big(\Psi(\zeta_1)(s)\big)^p \d s\geq \int_{l_i}^{r_i}\big(\Psi(\xi_1)(s)\big)^p\d s$,
\item[O4.]
$\int_{l_i}^{r_i} \Lambda^*(\dot\xi^{(a)}(s))\d s
=\int_{s_i}^{s_{i+1}} \Lambda^*(\dot\zeta_1^{(a)}(s))\d s$,
\item[O5.]
$\zeta_1^{(u)}(s_{i+1}-) - \zeta_1^{(u)} (s_i-) = 
\xi_1^{(u)}(r_{i}-) - \xi_1^{(u)} (l_i-)$.
\end{itemize}
Now, we verify the conditions \emph{i-1)}, \emph{i-2)}, \emph{i-3)}, \emph{i-4)}.
Note first that the conditions \emph{i-1)} and \emph{i-4)} are obvious from the construction of $\zeta_1$.
We can verify \emph{i-2)} as follows:
\begin{align*}
\Phi_T(\xi_2) 
&
= 
\int_0^T \big(\Psi(\zeta_1)(s)\big)^p \d s 
\geq 
\int_0^{s_\infty} \big(\Psi(\zeta_1)(s)\big)^p \d s 
= 
\sum_{i=1}^\infty \int_{s_i}^{s_{i+1}}\big(\Psi(\zeta_1)(s)\big)^p \d s 
\\
&\geq
\sum_{i=1}^\infty \int_{l_i}^{r_i}\big(\Psi(\xi_1)(s)\big)^p \d s 
=
\int_0^T \big(\Psi(\xi_1)(s)\big)^p \d s  
= 
\Phi_T(\xi_1),
\end{align*}
where the second inequality is from O3, and the second last equality is from O1.
Moving onto \emph{i-3)}, note that, due to the left continuity of $\zeta_1$, $s_n \to s_\infty$ implies that $\xi(s_n-) \to \xi(s_\infty-)$. 
Also, $\zeta_1^{(u)}(s_\infty) - \zeta_1^{(u)}(s_\infty-) =0$ and $\zeta_1^{(u)}$ is constant on $[s_\infty, T]$. 
Therefore, $\sum_{i=1}^\infty \Big(\zeta_1^{(u)}(s_{i+1}-)- \zeta_1^{(u)}(s_i-)\Big) = \lim_{n\to \infty} \zeta_1^{(u)} (s_{n+1}-) = \zeta_1^{(u)} (s_\infty-) = \zeta_1^{(u)}(T)$, 
where we adopted the convention that $\zeta_1^{(u)}(0-)=0$. 
From O4, O5, and this observation,
\begin{align*}
I_y^\BVT(\zeta_1) 
&
= 
\int_0^T \Lambda^*(\dot \zeta_1(t)) \d s + \theta_+ \cdot \zeta_1^{(u)}(T)
\\
&
=
\sum_{i=1}^\infty\int_{s_i}^{s_{i+1}} \Lambda^*(\dot \zeta_1(t)) \d s 
+ \theta_+ \cdot \sum_{i=1}^\infty \Big(\zeta_1^{(u)}(s_{i+1}-)- \zeta_1^{(u)}(s_i-)\Big) 
\\
&=
\sum_{i=1}^\infty\int_{l_i}^{r_i} \Lambda^*(\dot \xi^{(a)}(t)) \d s 
+ \theta_+ \cdot \sum_{i=1}^\infty \Big(\xi^{(u)}(r_i-)- \xi^{(u)}(l_i-)\Big) 
\\
&\leq
\int_0^T \Lambda^*(\dot \xi^{(a)}(t)) \d s 
+ \theta_+ \cdot \xi^{(u)}(T)
=
I_y^\BVT(\xi_1).
\end{align*}
%Note that 
%\begin{itemize}
%\item[O3.]
%$\xi_2(t) \geq 0$ for $t\in[0, s_\infty)$, and 
%\item[O4.]
%$\xi_2^{(u)}$ does not jump on $[s_\infty, T]$
%\end{itemize}

%Therefore, $\Phi_T(\xi_2) \geq \Phi_T(\xi_1)$ and $I_y^T(\xi_2) \leq I_y^T(\xi_1)$.
For part \emph{(ii)}, we construct $\zeta_2$ from $\zeta_1$
by moving all the jumps of $\xi^{(u)}$ to time 0. 
This neither increases $I_y^\BVT$ nor decreases $\Phi_T$.
That is, if we set 
$$
\zeta_2(t) \triangleq y +  \int_0^t \dot \zeta_1(s) ds + \zeta_1^{(u)}(T),  
$$
then 
$\Phi_T(\zeta_2) \geq \Phi_T(\zeta_1)$ obviously, and $\theta_+\cdot \zeta_1^{(u)}(T)+ I_{y+\zeta_1^{(u)}(T)}^T(\zeta_2) \leq I_y^T(\zeta_1)$.
Noting that $\zeta_1^{(u)}(T) \leq \xi^{(u)}(T)$, we see that $\zeta$ satisfies all the claims of the lemma.

For part \emph{(iii)},  
let $\zeta \in \AC[0,T]$ be a concave majorant of $\xi$.
Then there exists a non-increasing $\dot\zeta\in\D[0,T]$ such that $\zeta(t) = \xi(0) + \int_0^t \dot\zeta(s)ds$. (Due to the continuity of $\xi$, $\xi(0)$ and $\zeta(0)$ should coincide.)
Let $\zeta_3(t) \triangleq \xi(0)+ \int_0^t \mu \vee \dot \zeta(s)ds$.
Note that \emph{iii-1)}, \emph{iii-2)}, and \emph{iii-4)} are straightforward to check from the construction. 
To show that \emph{iii-3)} is also satisfied, 
we construct $\mathcal C^+_2\triangleq \{(l'_i, r'_i) \subseteq[0,T]: i\in \N\}$ in a similar way to $\mathcal C^+_1$ so that 
%\begin{itemize}
%    \item 
    the elements of $\mathcal C^+_2$ are non-overlapping, and
%    \item 
    $\xi(s) < \zeta_3(s)$ if and only if $s\in(l'_i, r'_i)$ for some $i\in\N$.
%\end{itemize}
Note that, due to the continuity of $\zeta$ and $\xi$, $\zeta(l'_i) = \xi(l'_i)$ and $\zeta(r'_i) = \xi(r'_i)$, and $\zeta$ has to be a straight line on $(l'_i, r'_i)$ for each $i\in \N$. 
Set $s_0 \triangleq 0\vee\sup\{t\in[0,T]: \dot\zeta(t)\geq \mu\}$.
Then, no interval in $\mathcal C^+_2$ contains $s_0$, because otherwise, $\zeta$ has to be a straight line in a neighborhood of $s_0$, and hence, $\dot \zeta$ has to be constant there, but this is contradictory to the definition of $s_0$.
Now, let $\dot\xi$ denote a derivative of $\xi$. 
Then $\int_{l'_i}^{r'_i} \Lambda^*(\mu\wedge \dot\zeta(s))ds = 
\int_{l'_i}^{r'_i} \Lambda^*(\mu)ds = 0$ for $i$'s such that $r'_i> s_0$, and hence,
\begin{align*}
I_y^\BVT(\xi) - I_y^\BVT(\zeta_3)
&
= 
\int_0^T \Lambda^*(\dot\xi(s)) ds - 
\int_0^T \Lambda^*(\mu \vee \dot\zeta(s)) ds
\\
&
\geq
\sum_{i\in\N:\ r'_i \leq s_0}\int_{l'_i}^{r'_i} \left(\Lambda^*(\dot\xi(s)) - \Lambda^*(\dot \zeta(s))\right) ds.
\end{align*}
Note that from the construction of $\mathcal C^+_2$, if $s\in [l'_i, r'_i]$ for some $i$ such that $r'_i\leq s_0$, we have that $\dot\zeta(s) = (\zeta_3(r'_i) - \zeta_3(l'_i))/(r'_i - l'_i)
=(\xi(r'_i) - \xi(l'_i))/(r'_i - l'_i)$,
%= \int_{l'_i}^{r'_i} \dot\xi(s)ds/(r'_i - l'_i)$ for $s\in(l'_i, r'_i)$
and hence, from Jensen's inequality,
\begin{align*}
\int_{l'_i}^{r'_i} \left(\Lambda^*(\dot\xi(s)) - \Lambda^*(\dot \zeta(s))\right) ds
&=
\int_{l'_i}^{r'_i} \Lambda^*(\dot\xi(s)) ds 
-\int_{l'_i}^{r'_i} \Lambda^*\big((\xi(r'_i) - \xi(l'_i))/(r'_i - l'_i)\big)ds
\\
&
=
\int_{l'_i}^{r'_i} \Lambda^*(\dot\xi(s)) ds
- (r'_i - l'_i)\cdot \Lambda^*\Big(\int_{l'_i}^{r'_i} \dot\xi(s)ds/(r'_i - l'_i)\Big)
\\
&
\geq 0.
\end{align*}
Therefore, $\zeta_3$ satisfies \emph{iii-3)} as well.
\eop

\subsection{Proof of Lemma~\ref{we-can-consider-only-upto-M}}
\label{pf-we-can-consider-only-upto-M}
\linksinpf{we-can-consider-only-upto-M}
Recall that $\bar{y}\triangleq (|\mu|(p+1))^{1/1+p}$. 
If $y \geq \bar y$, the equality in \eqref{eq:we-can-consider-only-upto-M} holds with the optimal values of the LHS and RHS both being zero:
to see this, we invoke Lemma~\ref{useful-properties-of-variational-By}---part \textit{i)}.
% set $M \triangleq -y/\mu$ and $\zeta(t) \triangleq y + \mu t$, and note that  
% $\int_0^{\T(\zeta)} R(\zeta) (s)^p ds \geq 1$ and $\T(\zeta)= M$, and hence, $\zeta \in B^{\AC;M}_{y}$ while $I_y(\zeta) = 0$. \todo{We are repeating the same argument as the one in the proof of Lemma~\ref{useful-properties-of-variational-By}. Separate the argument as a lemma?}
Moving on to the case $y< \bar y$, it is enough to show that there exists $M>0$ such that
\begin{equation}\label{wtp-lem-we-can-consider-only-upto-M}
\text{for any given $\xi \in B^{\AC}_{y} \setminus B^{\AC;M}_{y}$, one can find $\zeta \in B^{\AC;M}_{y}$ such that $I_y(\zeta) \leq I_y(\xi)$}.
\end{equation}
To construct such $M$, consider  $w$ and $z$  such that  $\mu < w < 0 < z$, $\Lambda^*(w) < \infty$ and $\Lambda^*(z) < \infty$. 
We consider a piece-wise linear path 
\[
\zeta(t)\triangleq (y + zt)\mathbbm{1}_{[0,(\bar{y}-y)/z]}(t) + (\bar{y}+\mu\big(t-(\bar{y}-y)/z)\big)\mathbbm{1}_{[(\bar{y}-y)/z,\,\infty)}(t)
\]
so that $\zeta \in B^{\AC;M}_{y}$ and $I_y(\zeta) = \Lambda^*(z) \frac{\bar y - y}{z}$. Let
\[
M\triangleq \max\left\{\frac{(\bar y - y)\Lambda^*(z)}{z\Lambda^*(w)}, \ (\bar y- y)/z-\bar y/\mu,\ -y/w\right\}.
\]
Suppose that $\xi \in B^{\AC}_{y}\setminus B^{\AC;M}_{y}$, and hence, $\T(\xi) > M$. 
Note that 
$\mu < w  < -y/\T(\xi) $ by the construction of $M$. We can now estimate
\begin{align*}
I_y(\xi) 
&=
\int_0^{\T(\xi)}\Lambda^*(\dot \xi(s)) ds 
\geq 
\T(\xi)\cdot \Lambda^*\Big(-y/\T(\xi)\Big)
\geq 
\T(\xi) \cdot \Lambda^*(w)
\geq 
M\cdot \Lambda^*(w)
\geq \frac {(\bar y- y)}  z\Lambda^*(z)
= I_{y}(\zeta).
\end{align*}
In this derivation, the first inequality is from Jensen's inequality. The second inequality follows since $\Lambda^*(x)$ is non-decreasing in $x\geq \mu$.
The third and the fourth inequalities are from the choice of $\xi$ and the construction of $M$, respectively.
This concludes the proof of \eqref{wtp-lem-we-can-consider-only-upto-M}. 

To see the existence of $c>0$ and $d>0$, note that for the case $y \geq \bar y$, our construction of $M(y)$ is linear in $y$, whereas $M(y)$ is bounded for the case $y < \bar y$.
\eop

\subsection{Proof of Lemma~\ref{B_y_star-equal-V_y_T_star}}
\label{pf-B_y_star-equal-V_y_T_star}
\linksinpf{B_y_star-equal-V_y_T_star}
The proof that  $\mathcal{B}_y^*={\mathcal{V}_y^*}^T$ for sufficiently large  $t$'s follows immediately from the following claims along with Lemma~\ref{useful-properties-of-variational-By}.

{\bf Claim 1:} ${\V_y^T}^*$ is nonincreasing in $T$. \\
\emph{Proof of Claim 1.} Let $t_1 < t_2$. For each $\xi_1 \in V_y^{t_1}$, consider $\xi_2(s) \triangleq \xi_1(s\wedge t_1)+  \mu(s - t_1)\ind_{(t_1,t_2]}(t)$. Then, $\xi_2 \in V_y^{t_2}$ and $I_y^{\BV[0,t_1]} (\xi_1) = I_y^{\BV[0,t_2]} (\xi_2)$. 
Therefore, ${\V_y^{t_2}}^*$ is at least as small as ${\V_y^{t_1}}^*$.

{\bf Claim 2:} If $M>0$ is such that $\inf_{\xi\in B_y^{\AC; M}} I_y(\xi) = \inf_{\xi\in B_y^\AC} I_y(\xi)$ as in Lemma~\ref{we-can-consider-only-upto-M}, then
$$\inf_{\xi\in B_y^{\AC;M}} I_y(\xi) \geq {\V_y^M}^*.$$
\emph{Proof of Claim 2.} Given an $\epsilon>0$, consider $\xi_\epsilon\in B_y^{\AC;M}$ such that $I_y(\xi_\epsilon) \leq \inf_{\xi\in B_y^{\AC;M}} I_y(\xi)+\epsilon$. 
Set 
$
\zeta_\epsilon(t) \triangleq \xi_\epsilon(t\wedge \T(\xi_\epsilon)) + \mu (t-\T(\xi_\epsilon)) \ind_{(\T(\xi_\epsilon),M]}(t).
$
%$
%\zeta_\epsilon(t) \triangleq 0\vee \xi_\epsilon(t\wedge \T(\xi_\epsilon)) + \mu (t-\T(\xi_\epsilon)) \ind_{(\T(\xi_\epsilon),M]}(t).
%$
%(The only difference between $\zeta_\epsilon$ and $\xi_\epsilon$ on $[0, \T(\xi_\epsilon)]$ is that $\zeta_\epsilon$ the downward jump size of $\zeta_\epsilon$ at $\T(\xi_\epsilon)$ may be smaller than $\xi_\epsilon$.)
Then, $\zeta_\epsilon \in V_y^M$ and hence,
$${\V_y^M}^* = \inf_{\xi\in V_y^M} I_y^\BVM(\xi) \leq I_y^\BVM(\zeta_\epsilon)= I_y(\xi_\epsilon) \leq \inf_{\xi\in B_y^M} I_y(\xi)+\epsilon.$$
Taking $\epsilon \to 0$, we arrive at Claim 2.

{\bf Claim 3:} 
For any $T>0$, 
$$\inf_{\xi\in B_y} I_y(\xi) \leq {\V_y^T}^*.$$
\emph{Proof of Claim 3.}
(Throughout this proof, we interpret $\theta_+\cdot z$ as 0 if $\theta_+ = \infty$ and $z = 0$. Likewise, we interpret $\epsilon/z = \infty$ if $z=0$.)
First note that the claim is trivial if ${\V_y^T}^* = \infty$, and hence, we only consider the case that ${\V_y^T}^* < \infty$.
Fix an $\epsilon>0$ and consider $\xi_1\in V_y^T$ such that $I_y^\BVT(\xi_1) < {\V_y^T}^*+\epsilon$. 
(Note that this implies that $\xi_1^{(u)}(T) = 0$ if $\theta_+ = \infty$.)
Due to  Lemma~\ref{merging-actions-to-the-left-lemma} part--\textit{(ii)}, there exists a path  $\xi_2 \in \AC[0,T]$ such that $\xi_2(0) = y + z$, $0 \leq z \leq \xi_1^{(u)}(T)$,  
$\xi_2$ is nonnegative over $[0,t_0]$ for some $t_0 \in [0, T]$, $\xi_2$ is affine with slope $\mu$ over $[t_0, T]$,  $\Phi_T(\xi_2)  \geq \Phi_T(\xi_1) \geq 1$, and 
\begin{equation}\label{1inequ}\theta_+ \cdot z + I_{y+z}^\BVT(\xi_2) \leq I_y^\BVT(\xi_1)\leq {\V_y^T}^* + \epsilon.\end{equation}
Recall the well known property of $\Lambda^*$ that  $\lim_{x \to \infty}\frac{\Lambda^*(x)}{x} = \theta_+$ (see, for example, equation (5.5) in \cite{mogulskii1993large}). Consequently, we can choose a $u> 0$ large enough so that
\begin{equation}\label{usefulclaim}
\Lambda^*(u)/u \leq \theta_+ + \epsilon/z.
\end{equation}
Set $\hat T \triangleq z/u+T$ and consider $\xi_3 \in \AC[0,\hat T]$ such that
\[
\xi_3(s) = (y + us)\mathbbm{1}_{[0, z/u]}(s)+\xi_2\left(s-z/u \right)\mathbbm{1}_{(z/u, \hat T]}(s), \ s \in [0,\hat T].
%+ \big(\xi_2(T) + \mu(s-z/u-T)\big)\ind_{(z/u+T, \infty)}(s). 
\]
Then, $\xi_3(0)=y$ and $\Phi_{\hat T}(\xi_3) \geq \Phi_T(\xi_2) \geq 1$. 
Moreover, in view of \eqref{1inequ} and \eqref{usefulclaim},
\begin{align*}
I_y^{\BV[0,\hat T]}(\xi_3) 
&
= (z/u)\Lambda^*(u)+\int_{z/u}^{\hat T}\Lambda^*(\dot{\xi_3}(s))ds
\\
&= (z/u)\Lambda^*(u)+\int_{0}^{T}\Lambda^*(\dot{\xi_2}(s))ds
\\
&\leq \theta_+\cdot z + \epsilon + I_{y+z}^\BVT(\xi_2) 
\\
&\leq  {\V_y^T}^* + 2\epsilon.
\end{align*}
Next, from part \emph{(iii)} of Lemma~\ref{merging-actions-to-the-left-lemma}, we know that there exists a path $\xi_4 \in \AC[0,\hat{T}]$ such that $ \xi_4(0) = y$, $\Phi_{\hat T}(\xi_4) \geq \Phi_{\hat T}(\xi_3) \geq 1$, $I_y^{\BV[0,\hat T]}(\xi_4) \leq I_y^{\BV[0,\hat T]}(\xi_3) \leq {\V_y^T}^* + 2\epsilon$, $\xi_4$ is concave on $[0,\hat T]$, and $\dot{\xi}_{4} $ is bounded by $\mu$ from below.
Finally, define $\xi \in \D[0,\infty)$ as
$$\xi(t) \triangleq \xi_4(t\wedge \hat T) + \mu([t - \hat{T}]^+),\quad t\geq 0.$$
Note that if $\cal T(\xi) \leq \hat T$, due to the concavity of $\xi$, $\Psi(\xi)$ and $\Psi(\xi_4)$ are zero after $\cal T(\xi)$. Therefore, $\Phi(\xi) = \Phi_{\hat T}(\xi_4) \geq 1$. 
If $\cal T(\xi) > \hat T$, 
then $\Phi(\xi) > \Phi_{\hat T}(\xi_4) \geq 1$.
That is, $\xi \in B_y$ in all cases. 
Moreover, since $\dot\xi(t) = \mu$ for $t\geq \hat T$,  $I_y(\xi) = I_y^{\BV[0,\hat{T}]}( \xi_4) \leq {\V_y^T}^* + 2\epsilon$. 
Therefore,
\[
\inf_{\xi \in B_y}I_y(\xi) \leq {\V_y^T}^*+2\epsilon.
\]
Since $\epsilon$ is arbitrary, this proves Claim 3.
%%%%%%%%%%%%%%%%%%%%%%%%%%%%%%%%%%%%%%%%%
%%%%%%%%%%%%%%%%%%%%%%%%%%%%%%%%%%%%%%%%%%%%%%%%%%%
%%%%%%%%%%%%%%%%%%%%%%%%%%%%%%%%%%%%%%%%%%%%%%%%%%%%%%%%%%
%%%%%%%%%%%%%%%%%%%%%%%%%%%%%%%%%%%%%%%%%%%%%%%%%%%%%%%%%%%%%%%%%%
\eop
\endproof

\subsection{Uniform continuity:  proof of Lemma~\ref{uniformcontinuity2}.}
\label{pf-uniformcontinuity2}

We first state two preparatory lemmas. 
Let \notationdef{nota-TV}{$\mathrm{TV}(\xi)$} be the total variation of $\xi$.

\begin{lemma}
    \label{uniformcontinuity}
    \linksinthm{uniformcontinuity}
    For any $M<\infty$, the function 
    $H: \D[0,T] \rightarrow [0,\infty)$ given by 
    $\notationdef{nota-H}{H(\xi)} \triangleq \int_0^T \xi (s)ds$
    is Lipschitz continuous on the set $\{\xi: \mathrm{TV}(\xi) + \xi(0) \leq M\}$ w.r.t.\ the $M_1'$ metric. 
\end{lemma}

\begin{proof}{Proof.}
\linksinpf{uniformcontinuity}
Let $\xi$ be such that $\mathrm{TV}(\xi) + \xi(0) \leq M$ and let $\zeta$ be such that  $d_{M_1'}(\xi, \zeta) \leq \epsilon$. 
Set $\eta_-(t) \triangleq \inf\{ x: d((x,t), \Gamma(\xi)) \leq \epsilon\}$, where $\Gamma(\xi)$ is the completed graph of $\xi$ and $d$ is the $L_1$ distance in $\R^2$, i.e., $d((x,t), (u,s)) = |x-y|+ |t-s|$. Then $d_{M_1'}(\xi, \zeta) \leq \epsilon$ implies that $\zeta(t) \geq \eta_-(t)$ for all $t\in[0, T]$.
Note that, if we denote the arc length of $\Gamma(\xi)$ with $\mathrm{len}(\Gamma(\xi))$, then $\mathrm{len}(\Gamma(\xi))$ is bounded by $T+\mathrm{TV}(\xi) + \xi(0)$.
Due to the construction of $\eta$ and the fact that $L_1$ balls are contained in $L_2$ balls of the same radius, the difference between the area below $\xi$ and the area below $\eta$ is bounded by $\mathrm{len}(\Gamma(\xi)) \cdot \epsilon$. Putting everything together, we conclude that 
\begin{equation*}
\label{integral-lowerbound}
    \int_0^T \xi(s) ds - \int_0^T \zeta (s) \leq \int_0^T \xi(s) ds - \int_0^T \eta_-(s) \leq \mathrm{len}(\Gamma(\xi)) \cdot \epsilon \leq \big(T+ \mathrm{TV}(\xi) + \xi(0)\big) \cdot \epsilon \leq (T+ M)\cdot \epsilon.
\end{equation*}
Similarly, by majorizing $\zeta$ with $\eta_+(t) \triangleq \inf\{ x: d((x,t), \Gamma(\xi)) \leq \epsilon\}$, we also get 
\begin{equation*}
    \int_0^T \xi(s) ds - \int_0^T \zeta (s) \geq  (T+ M)\cdot \epsilon,
\end{equation*}
proving the Lipschitz continuity of $H$ with Lipschitz constant $(T+M)$.
%\cmt{BZ: Doesn't this proof extend to $M_1'$?}\cmt{CR: Revised so that it is now w.r.t.\ to $M1'$}
\eop 
 \end{proof}

\begin{lemma}
\label{reflectionlipschitz}
\linksinthm{reflectionlipschitz}
The reflection map $\Psi$ is a Lipschitz continuous map from $\D[0,T]$  to $\D[0,T]$ w.r.t.\ the ${M}_1'$ topology with Lipschitz constant 2. 
\end{lemma}

\begin{proof}{Proof.}
\linksinpf{reflectionlipschitz}
The proof is a straightforward adaptation of the proof of Lipschitz continuity given in Theorem 13.5.1 of \cite{whitt2002stochastic} for the $M_1$ topology. That theorem is based on elementary estimates, and the key Lemma 13.5.3 in \cite{whitt2002stochastic},
which establishes that parametric representation of a path $\xi \in \D[0,T]$ is preserved under taking projections. The proof of this property for $M_1$, given in \cite{whitt2002stochastic}, p.\ 440 extends to $M_1'$ by using the definition of extended completed graph $\Gamma'$ for $M_1'$, rather than the completed graph $\Gamma$ for $M_1$. 
Along this lines, it follows that if $(u,t)$ is a parametric representation of $\Gamma'(\xi)$, then $(\Psi(u), t)$ is a parametric representation of $\Gamma'(\Psi(\xi))$.
Using this result, the steps in the proof of Theorem 13.5.1 of \cite{whitt2002stochastic} follow verbatim
for $M_1'$.  We omit the details. 
\eop
\end{proof}

Now we are ready to prove Lemma~\ref{uniformcontinuity2}.

\begin{proof}{Proof of Lemma~\ref{uniformcontinuity2}.}
\linksinpf{uniformcontinuity2}
Suppose that $I_0^{\BV[0,T]}(\xi), I_0^{\BV[0,T]}(\zeta) \in [0,\gamma]$. Then, from the inequality (13) in \cite{vysotsky}, we know that there exists $\gamma'$ such that $TV(\xi), TV(\zeta) \in [0, \gamma']$. To prove the uniform continuity, suppose that $d_{M_1'} (\xi, \zeta) < \epsilon$. Then, $d_{M_1'} (\Psi(\xi), \Psi(\zeta) ) < 2\epsilon$ by Lemma~\ref{reflectionlipschitz}, and $TV(\Psi(\xi)), TV(\Psi(\zeta)) \in[0, 2\gamma']$.
In turn, we have that $d_{M_1'} (\Psi(\xi)\vee \epsilon, \Psi(\zeta)\vee\epsilon ) < 2\epsilon$ and $TV(\Psi(\xi)\vee\epsilon), TV(\Psi(\zeta)\vee\epsilon) \in [0, 2\gamma']$.
Using  the mean value theorem, we obtain the following inequality for $x,y,a,b \in [0,\infty)$ such that $x, y \in [a,b]$:
\begin{equation}\label{gen-inequal-p-th-power}
|x^p-y^p| \leq p (a^{p-1} \vee b^{p-1}) |x-y|.
\end{equation} 
Now, suppose that $(u,t)$ and $(v,s)$ are the parametric representations of $\Psi(\xi)\vee\epsilon$ and $\Psi(\zeta)\vee\epsilon$, respectively. Then, 
there exists $r_\xi\in[0,1]$ such that $u(r) \leq \epsilon$ for $r\leq r_\xi$ and $u(r) \geq \epsilon$ for $r\geq r_\xi$. Likewise, there exists $r_\zeta\in[0,1]$ such that $v(r) \leq \epsilon$ for $r\leq r_\zeta$ and $v(r) \geq \epsilon$ for $r\geq r_\zeta$. We assume w.l.o.g.\ that $r_\xi \leq r_\zeta$. 
Note that, since $u(r), v(r) \in [0,\epsilon]$ on $r\in[0, r_\xi]$, we get
\begin{align*}
\sup_{[0,r_\xi]}|u^p(r) - v^p(r)|  \leq \epsilon^p.
\end{align*}
Also, since $u(r) \in [\epsilon, 2\gamma']$ and $v(r) \in [0, \epsilon]$ on $r\in[r_\xi,r_\zeta]$, we get from \eqref{gen-inequal-p-th-power}
\begin{align*}
\sup_{[r_\xi,r_\zeta]}|u^p(r) - v^p(r)|  
&\leq \sup_{[r_\xi,r_\zeta]}\{|u^p(r) - \epsilon^p|+|\epsilon^p - v^p(r)|\}  \leq \sup_{[r_\xi,r_\zeta]}\{|u^p(r) - \epsilon^p|+\epsilon^p\} 
\\
&\leq \sup_{[r_\xi,r_\zeta]}\{p(\epsilon^{p-1}\vee (2\gamma')^{p-1}) |u(r) - \epsilon|+\epsilon^p\} 
\\
&\leq \sup_{[r_\xi,r_\zeta]}\big\{p(\epsilon^{p-1}\vee (2\gamma')^{p-1})\big( |u(r)-v(r)| + |v(r) - \epsilon| \big)+\epsilon^p\big\} 
\\
&\leq \sup_{[r_\xi,r_\zeta]}\big\{p(\epsilon^{p-1}\vee (2\gamma')^{p-1})\big( |u(r)-v(r)| + \epsilon \big)+\epsilon^p\big\}
\\
&\leq p(\epsilon^{p-1}\vee (2\gamma')^{p-1})\big( \|u-v\|_\infty + \epsilon \big)+\epsilon^p.
\end{align*}
Finally, since $u(r),v(r) \in [\epsilon, 2\gamma']$ on $r\in[r_\zeta, 1]$, we get again from \eqref{gen-inequal-p-th-power}
\begin{align*}
\sup_{[r_\zeta,1]}|u^p(r) - v^p(r)| \leq \sup_{[r_\zeta,1]}p(\epsilon^{p-1} \vee (2\gamma')^{p-1}) |u(r) - v(r)| \leq p(\epsilon^{p-1} \vee (2\gamma')^{p-1})\|u-v\|_\infty. 
\end{align*}
From these inequalities, we see that if $(u,t) \in \Gamma'(\Psi(\xi)\vee\epsilon)$ and $(v,s) \in \Gamma'(\Psi(\zeta)\vee\epsilon)$,
\begin{align*}
\|u^p - v^p\|_\infty 
&= \sup_{[0,r_\xi]}|u^p(r) - v^p(r)| \ \vee\  \sup_{[r_\xi,r_\zeta]}|u^p(r) - v^p(r)| \ \vee\ \sup_{[r_\zeta,1]}|u^p(r) - v^p(r)|
\\
&\leq p(\epsilon^{p-1}\vee (2\gamma')^{p-1})\big( \|u-v\|_\infty + \epsilon \big)+\epsilon^p.
\end{align*}
Now, we can bound the $M_1'$ distance between $\big(\Psi(\xi)\vee \epsilon)^p$ and $\big(\Psi(\zeta)\vee \epsilon\big)^p$ as follows:
\begin{align*}
&d_{M_1'} \Big(\big(\Psi(\xi)\vee \epsilon)^p, \big(\Psi(\zeta)\vee \epsilon\big)^p\Big) 
\\
&= \inf_{\substack{(u,t) \in \Gamma'((\Psi(\xi)\vee\epsilon)^p) \\ (v,s) \in \Gamma'((\Psi(\zeta)\vee\epsilon)^p)}} \{\|u-v\|_\infty + \|t-s\|_\infty\}
\\
&= \inf_{\substack{(u,t) \in \Gamma'(\Psi(\xi)\vee\epsilon) \\ (v,s) \in \Gamma'(\Psi(\zeta)\vee\epsilon)}} \{\|u^p-v^p\|_\infty + \|t-s\|_\infty\}
\\
&\leq \inf_{\substack{(u,t) \in \Gamma'(\Psi(\xi)\vee\epsilon) \\ (v,s) \in \Gamma'(\Psi(\zeta)\vee\epsilon)}} \{p(\epsilon^{p-1}\vee (2\gamma')^{p-1})\big( \|u-v\|_\infty + \epsilon \big)+\epsilon^p + \|t-s\|_\infty\}
\\
&\leq \Big(1\vee \big(p(\epsilon^{p-1}\vee (2\gamma')^{p-1})\big)\Big)\inf_{\substack{(u,t) \in \Gamma'(\Psi(\xi)\vee\epsilon) \\ (v,s) \in \Gamma'(\Psi(\zeta)\vee\epsilon)}} \{ \|u-v\|_\infty + \|t-s\|_\infty\} + p(\epsilon^{p-1}\vee (2\gamma')^{p-1})\epsilon + \epsilon^p
\\
&= \Big(1\vee \big(p(\epsilon^{p-1}\vee (2\gamma')^{p-1})\big)\Big)d_{M_1'} (\Psi(\xi)\vee\epsilon, \Psi(\zeta)\vee\epsilon) + p(\epsilon^{p-1}\vee (2\gamma')^{p-1})\epsilon + \epsilon^p
\\
&\leq \Big(1\vee \big(p(\epsilon^{p-1}\vee (2\gamma')^{p-1})\big)\Big)2\epsilon + p(\epsilon^{p-1}\vee (2\gamma')^{p-1})\epsilon + \epsilon^p
\\
&\triangleq \delta(\epsilon).
\end{align*}
Note that $\delta(\epsilon) \to 0$ as $\epsilon \to 0$.

To apply Lemma~\ref{uniformcontinuity}, we examine the total variations of $(\Psi(\xi) \vee \epsilon)^p$ and $(\Psi(\zeta) \vee \epsilon)^p$. 
Recall the notation $\mathbf T=\cup_{d=1}^\infty\{(t_1,\ldots,t_d): 0\leq t_1 < t_2 <\cdots < t_d \leq 1\}$.
From \eqref{gen-inequal-p-th-power}, 
\begin{align*}
TV((\Psi(\xi) \vee \epsilon)^p) 
&= \sup_{\bf t \in \cal P} \sum_{i=1}^{n_{\bf t}} |(\Psi(\xi)\vee \epsilon)^p(t_i) - (\Psi(\xi)\vee \epsilon)^p(t_{i-1})|
\\
&\leq p (\epsilon^{p-1}\vee (2\gamma')^{p-1}) \sup_{(t_0,\ldots,t_k) \in \bf T} \sum_{i=1}^{k} |(\Psi(\xi)\vee \epsilon)(t_{i}) - (\Psi(\xi)\vee \epsilon)(t_{i-1})|
\\
&= p (\epsilon^{p-1}\vee (2\gamma')^{p-1}) TV( \Psi(\xi) \vee \epsilon) 
\\
&\leq p (\epsilon^{p-1}\vee (2\gamma')^{p-1}) 2\gamma'.
\end{align*}
\BZA{didnt we use other notation before for this or am i mistaken? Check}\CRA{We used $\mathbf T$ \hyperlink{temp_pin_1}{[here]} for a similar object. I think we defined $\cal P$ as above to facilitate the notation in the summation. I changed it to $\mathbf T$.}
Similarly, $$TV((\Psi(\zeta) \vee \epsilon)^p)  \leq p (\epsilon^{p-1}\vee (2\gamma')^{p-1}) 2\gamma'.$$
These two bounds allow us to apply Lemma~\ref{uniformcontinuity} to $H$ to obtain 
\begin{align*}
d_{M_1'} \Big(H\big((\Psi(\xi)\vee\epsilon)^p\big), H\big((\Psi(\zeta)\vee\epsilon)^p\big)\Big) 
&\leq (T+2\gamma' p (\epsilon^{p-1}\vee (2\gamma')^{p-1})) \cdot d_{M_1'} \big(G(\xi)\vee\epsilon^p, G(\zeta)\vee\epsilon^p\big)
\\
&\leq (T+2\gamma' p (\epsilon^{p-1}\vee (2\gamma')^{p-1})) \delta (\epsilon).
\end{align*}
Therefore,
\begin{align*}
&d_{M_1'}(\Phi_T(\xi), \Phi_T(\zeta)) 
\\
&= d_{M_1'} \big(H(G(\xi)), H(G(\zeta))\big) 
\\
&\leq d_{M_1'} \big(H(G(\xi)), H(G(\xi)\vee\epsilon^p)\big)  + d_{M_1'} \big(H(G(\xi)\vee\epsilon^p), H(G(\zeta)\vee\epsilon^p)\big)  + d_{M_1'} \big(H(G(\zeta)\vee\epsilon^p), H(G(\zeta))\big)  
\\
&\leq \epsilon^p T + \delta(\epsilon) + \epsilon^p T.
\end{align*}
This concludes the proof of the desired uniform continuity.
% \cmt{CR: Instead of citing (14) of \cite{vysotsky}, we can keep the argument below so that we make $R'$ explicit.}
% Let $\xi$ be such that $I_0^{\BV[0,T]}(\xi) \leq R$.
% Let $\delta \in (0, \min \{ \theta_+, |\theta_-|\})$. 
% Observe that $\Lambda^*(\dot{\xi}(s)) \geq \delta |\dot{\xi}(s)| - \Lambda (\delta)\vee\Lambda(-\delta)$.
% Consequently, 
% \[
% \int_0^1 |\dot{\xi}(s)| ds + \xi^u(1) + |\xi^d(1)| \leq (R + \Lambda (\delta)\vee\Lambda(-\delta))/\delta \triangleq  R'.
% \]
\eop 
%
%\noindent
%\cmt{%
%CHR: Sketch of the proof for the uniform continuity of $\xi %\mapsto \int_0^T \xi(s) ds$:\\
%1) The total variation is bounded in a sub-level set of %$I_0^{\BV[0,T]}$ (by the same argument as the one for boundedness of %the paths themselves). That is, there exists a constant %$M_\alpha$ such that $\mathrm{TV}(\xi) \leq M_\alpha$ for %all $\xi$ such that $I_0^{\BV[0,T]}(\xi) \leq \alpha$, where %$\mathrm{TV}(\xi)$ denotes the total variation of $\xi$.\\
%2) Let $\xi$ be any given path in the sub-level set.\\
%3) Set $\eta(t) \triangleq \inf\{ x: d((t,x), \Gamma(\xi)) %\leq \epsilon\}$ where $\Gamma(\xi)$ is the parametrization %of $\xi$, and $d$ is the $L_1$ distance in $\R^2$, i.e., %$d((t,x), (s,y)) = |t-s|+|x-y|$. \\
%4) Then $d_{M_1}(\xi, \zeta) \leq \epsilon$ implies that %$\zeta(t) \geq \eta(t)$ for all $t\in[0, T]$. \\
%5) The difference between the area below $\xi$ and the area %below $\eta$ is bounded by $\mathrm{len}(\Gamma(\xi)) \times \epsilon$, where $\mathrm{len}(\Gamma(\xi))$ denotes the length of $\Gamma(\xi)$. (This is because of the construction of $\eta$ and the fact that $L_1$ balls are contained in the $L_2$ ball of the same radius.)\\
%6) $\mathrm{len}(\Gamma(\xi))$ is bounded by %$T+\mathrm{TV}(\xi)$.\\
%7) From 1), 5), and 6), we know that $\int_0^T \xi(s) ds - %\int_0^T \zeta(s) \geq \int_0^T \xi(s) ds - \int_0^T %\eta(s) \geq (T+ M_\alpha) \epsilon$.\\
%8) The upper bound can be established in the same way. 
%}
\end{proof}

%\label{subsection-proofsec2}

%Now we are ready to prove Proposition~\ref{bp-ison-bzero}.

\newpage

\begin{APPENDICES}
	\section{Results on the theory of Markov chains.}\label{genresults}
	Let $\{X_m, -\infty < m < \infty\}$ be a geometrically ergodic stationary Markov chain on the state space $S$ (which includes an element $0$) and  invariant distribution $\pi$, such that $\pi(\{0\}) = \pi(0) > 0$. 
	Let $\{X_m^*, -\infty < m < \infty\}$ be the time-reversed stationary version of $\{X_m, -\infty < m < \infty\}$. 
	It is well known that $(X^*_{0},...,X^*_{k})\stackrel{d}{=} (X_{k},..., X_{0})$ for any $k \geq 0$; see, \cite{Sforzo05} and references therein for a discussion on reversibility for general state-space Markov chains. The following lemma follows directly by applying this identity. %In fact, the identity can be seen to hold path-wise since we can define $X_{n}=X^*_{-n}$, assuming that $X_{0}$ follows $\pi$.	
	%Let $Q$ be the transition kernel of $X_n$.
	%It is well known that the transition kernel $(Q^{\gets})$ of $X_n^*$ is given by
	%\begin{equation}\label{ridi}
	%$	Q^{\gets}(y,dx)=\frac{\pi(dx)}{\pi(dy)} Q(x,dy).$
	%	\end{equation}
	\CRA{Should we provide a reference for Lemma~\ref{timereversibility}?}\JBA{Let me try to find a reference, but it appears to me that this is immediate from the identity: $(X^*_{0},...,X^*_{k})\stackrel{d}{=} (X_{k},..., X_{0})$ for any $k \geq 0$ as stated.}\CRA{Thanks!}%
 
	\begin{lemma}
	\label{timereversibility}
	\linksinthmwopf{timereversibility}
		It holds that
		\begin{align}
		\P_0\left(X_i \in A_i: 1 \leq i \leq k\right) &=\frac{1}{\pi(0)}\P_{\pi}\left(X_i^* \in A_{k-i}: 0 \leq i \leq k-1, X_k^*=0 \right),\\
		\E_0\left[g(X_0,X_1,...,X_k) \right] &= \frac{1}{\pi(0)}\E_\pi\left[g(X^*_k,X^*_{k-1},...,X^*_0)I(X_k^*=0) \right],
		\end{align}
		for any non-negative integer $k$ and measurable $g:\R_+^{k+1} \to \R$.
	\end{lemma}
	
%	\proof{Proof.}
%	The second claim follows from the first claim. For the first claim observe that
%	\begin{align*}
%	&
%	\P_0\left(X_i \in A_i: 0 \leq i \leq n\right) 
%	=  \int_{A_0}\ldots \int_{A_n} Q(0,dx_1)\ldots Q(x_{n-1},dy_n) \\
%	&
%	=
%	\int_{A_0}\ldots \int_{A_n}\frac{1}{\pi(0)}\pi(dx_1)Q^{\gets}(x_1,0)\frac{\pi(dx_2)}{\pi(dx_1)}Q^{\gets}(x_2,dx_1) \ldots  \frac{\pi(dx_n)}{\pi(dx_{n-1})}Q^{\gets}(x_n,dx_{n-1}) dx_0\ldots dx_n \nonumber
%	\\
%	&
%	=
%	\frac{1}{\pi(0)}\int_{A_0}\ldots \int_{A_n}\pi(dx_n)Q^{\gets}(x_1,0)Q^{\gets}(x_2,dx_1) \ldots  Q^{\gets}(x_n,dx_{n-1}) dx_0\ldots dx_n.\nonumber 

%	& 	=
%	\frac{1}{\pi(0)}\P_{\pi}\left(X_i^* \in A_{n-i}: 1 \leq i \leq n; X_n^* = 0\right).
%	\end{align*}	\Halmos
%	\endproof

	Using the previous result, we can now establish the following lemma.
	
	\begin{lemma}\label{lastcycleequivalence}
	\linksinthm{lastcycleequivalence}
%		Suppose that $\{X_k: k \in (-\infty,\infty)\}$ is a non-negative time stationary Markov process with an associated time reversed Markov process $\left(X^*_k\right)$.
Define $T=\inf\{n \geq 1: X_n=0\}$, $T^*=\inf\{n \geq 1: X^*_n=0\}$,
and suppose that $\P_{0}\left(T >n\right)=\bigO(e^{-cn})$ for  some  $c>0$. In addition, let $n_0$ be such that $\inf_{k\geq n_0} \P_0(X_k=0) \geq \pi (0)/2$. Then, 
\begin{equation}
 \P_{\pi}\left(\sum_{k=0}^{T^*-1}(X_k^*)^p \geq x, X_n^*=0\right)
		\leq  (n+1) \P_{\pi}\left(\sum_{k=0}^{T-1}X_k^p \geq x\right), 
		\end{equation}
and
		\begin{equation}
		\P_{\pi}\left(\sum_{k=0}^{T^*-1}(X_k^*)^p \geq x, X_n^*=0\right)\geq 
		   (\pi(0)^2/2)  \P_{0}\left(\sum_{k=1}^{T}X_k^p \geq x\right)-\bigO(e^{-cn}).
		\end{equation}
	\end{lemma}

	\proof{Proof.}
	\linksinpf{lastcycleequivalence}
	We first derive the upper bound by noting that
	\begin{align*}
	\P_{\pi}\left(\sum_{k=0}^{T^*-1}(X_k^*)^p \geq x,\ X^*_n=0\right)
	&=
	\sum_{m=0}^{n}\P_{\pi}\left(\sum_{k=0}^{m-1}(X_k^*)^p \geq x,\ T^*=m,\ X^*_n=0\right)
	\\
	&=
	\sum_{m=0}^{n}\P_{\pi}\left(\sum_{k=0}^{m-1}(X_k^*)^p \geq x,\ X^*_1>0, \ldots, X^*_{m-1}>0,\ X^*_m = 0,\ X^*_n=0\right)
%	\\
%	&=
%	\sum_{m=0}^{n}\P_{\pi}\left(\sum_{k=0}^{m}(X_k^*)^p \geq x,\ X^*_1>0, \ldots, X^*_{m-1}>0,\ X^*_m = 0,\ X^*_n=0\right)
%	\\
%	&\leq
%	\sum_{m=0}^{n}\P_{\pi}\left(\sum_{k=0}^{m}(X_k^*)^p \geq x,\ X^*_1>0, \ldots, X^*_{m-1}>0, X^*_m = 0\right)
	\\
	&\leq
	\sum_{m=0}^{n}\P_{\pi}\left(\sum_{k=0}^{m-1}(X_k^*)^p \geq x,\ X^*_1>0, \ldots, X^*_{m-1}>0\right)
	\\
%	&
%	\leq
%	\sum_{m=0}^{n}\P_{\pi}\left(\sum_{k=0}^{m-1}(X^*_{k-(m-1)})^p\geq x,\ X^*_{-(m-1)}>0,\ldots, X^*_{-1}>0 \right) %, X^*_1=0\right)
%	\\
	&
	=
	\sum_{m=0}^{n}\P_{\pi}\left(\sum_{k=0}^{m-1}X_{m-1-k}^p\geq x, X_{m-1}>0,\ldots,X_1>0 \right) %, X_{-1}=0\right)
	\\
	&
	=
	\sum_{m=0}^{n}\P_{\pi}\left(\sum_{k=0}^{m-1} X_{k}^p\geq x,\ T\geq m\right)
	\leq
	\sum_{m=0}^{n}\P_{\pi}\left(\sum_{k=0}^{T-1}X_{k}^p\geq x,\ T\geq m\right)\\
	&\leq (n+1) \P_{\pi}\left(\sum_{k=0}^{T-1}X_{k}^p\geq x\right).
	\end{align*}
%	Thus,
%	\begin{align}\label{eq:1decisiveineq}
%	&
%	\P_{\pi}\left(\sum_{k=0}^{T^*} X^*_k \geq x, X^*_n=0\right)
%	\leq \E\left(\mathbf{1}\left(\sum_{k=0}^{T}X_k \geq x\right)\min\{(T+1),n\}\right) 
%	\leq
%	n\P_{\pi}\left(\sum_{k=0}^{T}X_k \geq x\right).
%	\end{align}
For the lower bound, first write 
	\begin{align*}
	\P_{\pi}\left(\sum_{k=0}^{T^*-1}(X_k^*)^p \geq x,\ X_n^*=0\right)
	&=
	\sum_{m=1}^{n} \P_{\pi}\left(\sum_{k=0}^{m-1}(X_k^*)^p \geq x,\ T^*=m,\ X_n^*=0\right)
	\\
	&=
	\sum_{m=1}^{n} \P_{\pi}\left(\sum_{k=0}^{m-1}(X_k^*)^p \geq x,\ T^*=m \right)\P_{0}\left(X_{n-m}=0\right).
	\end{align*}
	Applying Lemma~\ref{timereversibility} with $k=m$ and $g(y_0,...,y_m) = I(\sum_{i=0}^{m}y_i^p > x,\ y_0> 0,\ y_1>0,\ldots,\ y_{m-1}>0)$, we obtain
	\begin{align*}
	    \P_{\pi}\left(\sum_{k=0}^{m-1} (X_k^*)^p \geq x, T^*=m \right) 
	    &=  \P_\pi\left(\sum_{k=0}^{m-1} (X_k^*)^p \geq x,\ X^*_0>0,\ldots,\ X^*_{m-1}>0,\ X^*_m=0\right) \\
	    &=  \P_\pi\left(\sum_{k=0}^{m} (X_k^*)^p \geq x,\ X^*_0>0,\ldots,\ X^*_{m-1}>0,\ X^*_m=0\right) \\
	    &= \pi(0) \P_0\left(\sum_{k=1}^{m} X_k^p \geq x, X^*_i>0, i=1,\ldots, m-1 \right) \\
	    &= \pi(0) \P_0\left(\sum_{k=1}^{m} X_k^p \geq x, T \geq m \right) \\
	    &\geq \pi(0) \P_0\left(\sum_{k=1}^{m} X_k^p \geq x, T =m \right).
%	    \\
%	    &= \pi(0) \P_0\left(\sum_{k=0}^{T} X_k^p \geq x, T =m \right).
	\end{align*}
	Consequently, for every fixed $n_0$ such that $\inf_{k\geq n_0} \P_0(X_k=0) \geq \pi (0)/2$,
	\begin{align*}
	\P_{\pi}\left(\sum_{k=0}^{T^*-1}(X_k^*)^p \geq x, X_n^*=0\right)
	&\geq  \pi(0) \sum_{m=1}^{n-n_0} \P_0\left(\sum_{k=1}^{T} X_k^p \geq x, T =m \right)\P_{0}\left(X_{n-m}=0\right)
	\\
	&\geq
     \pi(0) \P_0\left(\sum_{k=1}^{T} X_k^p \geq x, T\leq n-n_0 \right)\inf_{k\geq n_0} \P_0(X_k=0).\\
	&\geq (\pi(0)^2/2) \P_0\left(\sum_{k=1}^{T} X_k^p \geq x \right)-\bigO(e^{-cn}).
	\end{align*}
	
	\Halmos
	
	\endproof

\section{ LDP results.}
%A useful result on large deviations on renewal processes can be found in \cite{puhalskii1997functional}. %
%\begin{result}\label{LDP-for-regeneration-times}
%	Assume that $\E(e^{\alpha \tau_1})< \infty $ for some $\alpha > 0$. Let $\alpha^*=\sup\{\alpha: %\E(e^{\alpha \tau_1})< \infty\}$.  \\
%	$\textbf{(i)}$  $\{N_n, \ n \geq 1\}$ obeys the LDP in $E^{\uparrow}$ for the $M_1'$ topology with the rate function
%	\[
%	I_{N}(\xi) \triangleq \int_{0}^{\infty}\sup_{\alpha < \alpha^*}\{\alpha-\dot{\xi}_1^{l}(t)\log\E\left(e^{\alpha \tau_1}\right)\}dt-\log\P\left(T_1=0\right)\xi_2^l(\infty),
%	\]
%	where $\xi=\xi^l_1+\xi_2^l$ is the Lebesque decomposition of $\xi$ with respect to Lebesque measure; $%\xi_1^l$ is the absolutely continuous component with $\xi_1^l(0)=0$, $\xi_2^l$ is the singular component and $\dot{\xi}_1^l$ is the derivative. \\
%	$\textbf{(ii)}$ If, in addition, $\P(\tau_1=0)=0$, then the LDP holds for the $J_1$ topology with $I_N(\xi)=\infty$ if $\xi$ is not absolutely continuous.
%\end{result}
We collect some LDP results that have appeared in the literature.
A straightforward adaptation of Corollary 3.2 in \cite{bazhba2017sample} to our  context is the following:
\begin{result}
\label{sp-ldp-m1'-top}
\linksinthmwopf{sp-ldp-m1'-top}
Let $\bar{\mathcal K}_n \triangleq \frac{1}{n}\sum_{i=1}^{\lfloor nt \rfloor}\mathcal U_i, \ t \in [0,1]$ be the scaled random walk
% $\{K_m, m\geq 0\}$ %THIS NOTATION HAS BEEN USED FOR OTHER PURPOSES IN THIS PAPER 
% is a random walk 
driven by an i.i.d.\ sequence $\{\mathcal U_i, i\geq 1\}$.
Assume that $\E[e^{s\mathcal U_1}]<\infty$ for some $s<0$, and $\P(\mathcal U_1 \geq x) =e^{-L(x)x^{\alpha}}$ for $\alpha \in (0,1)$. 
Suppose  that $L$ is a slowly-varying function, and $L(x)x^{\alpha-1}$ is eventually decreasing. 
Then, \rvtxt{5}{5)}{\rvout{$\bar K_n^0$}\rvin{$\bar{\mathcal{K}}_n$}} satisfies the LDP in $(\D[0,T],\mathcal{T}_{M_1'})$ with the speed $ L(n) n^{\alpha}$ and the rate function $I_{M_1'}:\mathbb D[0,T]\to [0,\infty]$,
	\begin{equation*}\label{eq:LDP-ratefunction-M1'}
	I_{M_{1}'}(\xi)
	\triangleq
	\begin{cases}
	\sum_{t\in[0,1]}\left(\xi(t)-\xi(t-)\right)^{\alpha}
	&
	\text{if $\xi \in \D^{(\E \mathcal U_1)}[0,T]$ with $\xi(0) \geq 0$,}
	\\
	\infty
	&
	\text{otherwise}.
	\end{cases}
	\end{equation*}
\end{result}

%	\begin{result}[\cite{bazhba2017sample}]\label{lemma:multi-d-equivalence-J}
%		Let $ \bar X_n^{(i)}, \ i=1,\ldots,d, $ satisfy an LDP on $(\D,\mathcal{T}_{M_1'})$ with speed $c_i L(n) n^{\frac{1}{1+p}}$, where $c_i \in (0,\infty)$. Then, $(\bar X_n^{(1)},\ldots,\bar X_n^{(d)})$ satisfies the LDP in $\big(\prod_{i=1}^d\mathbb D, \prod_{i=1}^d\mathcal T_{M_1'}\big)$ with speed $L(n)n^\frac{1}{1+p}$ and the rate function $I^d:\prod_{i=1}^d\mathbb D\to [0,\infty]$
%		\begin{equation}
%		I^{(d)}(\xi_1,\ldots,\xi_d)
%		\triangleq
%		\begin{cases}
%		\sum_{i=1}^d c_i \cdot \sum_{t\in[0,1]}\left(\xi_i(t)-\xi_i(t-)\right)^\frac{1}{1+p}
%		&
%		\text{if \ $\xi_i\in \mathbb D^{\uparrow}_p$ for $i=1,\ldots,d$}
%		\\
%		\infty
%		&
%		o.w.
%		\end{cases}
%		.
%		\end{equation}
%		
%	\end{result}
The following result, which is folklore, but explicitly stated in \cite{nuyenszwart}, provides the logarithmic asymptotics for the invariant distribution $\pi$ of $\{X_n\}_{n \geq 0}$ with $ X_{n+1} = [X_n+U_{n+1}]^{+}, \, n\geq 0$, 
	with $\{U_i, i\geq 1\}$  are i.i.d.\ such that $\E(U_1)=\mu <0$.

%\todo{in case $\beta = \theta_+$ and $\E e^{\beta U} < 1$, the cited result (Theorem 5 in \cite{borovkov1997large}) seems to require extra conditions.}
\begin{result}[\cite{nuyenszwart}]\label{steadystateasymptotics}
\linksinthmwopf{steadystateasymptotics}
Recall that $\beta = \sup \{s: \E[e^{sU}]\leq 1\}$. It holds that,
	\[
	\lim_{n \to \infty}\frac{\log\pi([n,\infty))}{n}=-\beta.
	\]
\end{result}
Finally, we mention a recent sample-path LDP for random walks with light-tailed increments, developed in \cite{vysotsky}, that we use in this paper. 
%Recall that The Skorokhod $J_1$ topology is induced by the Skorokhod metric, $d_{J_1}(\xi,\zeta) \triangleq \inf_{\lambda \in \Lambda}\left\{\max\{\|\xi(t)-\zeta(\lambda(t))\|_{\infty}, \|\lambda(t)-e(t)\|_{\infty}\}\right\},$ where $\xi, \zeta \in \D$, $\Lambda$ is the set of non-decreasing  homeomorphism  of $[0,1]$ onto itself, $e(t)=t$ is the identity map, and $\|\xi\|_{\infty}=\sup_{t \in [0,1]}|\xi(t)|$ is the sup-norm.  
%Let $\mathcal{T}_{J_1}[0,1]$ denote the topology induced by $d_{J_1}$.
%Denote with $\xi^{(u)}(1) \ \text{and} \ \xi^{(l)}(1)$, the total size of the upward and downward jumps respectively.  
\begin{result}
\label{extensionmogulskildp}
\linksinthmwopf{extensionmogulskildp}
%\todo{Wrong}
Recall that $\bar{K}_x^y = y + \frac{1}{x}\sum_{i=1}^{\lfloor xt \rfloor}U_i, \ t \in [0,1]$ is the scaled random walk
% $\{K_m, m\geq 0\}$ %THIS NOTATION HAS BEEN USED FOR OTHER PURPOSES IN THIS PAPER 
% is a random walk 
driven by the i.i.d.\ sequence $\{U_i, i\geq 1\}$, which satisfies {\bf Assumptions}~\ref{finitemgf} and \ref{regularitycondition}. 
	Recall also that
	\begin{equation}\label{eq:ratefunctioninmrw}
	I_y^{\BV[0,T]}(\xi) =
	\begin{cases}
	\int_{0}^{T}\Lambda^*(\dot{\xi}^{(a)}(s))ds+ \theta_+(\xi^{(u)}(T))+\theta_-|\xi^{(l)}(T)| & \ \text{if} \ \xi \in \BV[0,1] \text{ and } \ \xi(0)=y,\\
	\infty & \ \text{otherwise}.
	\end{cases}
	\end{equation}
	\begin{itemize}
	\item[(i)] (\cite{borovkov2013large, borovkov2014large}) $\bar K_x^y$ (as $x\to\infty$) satisfies a large deviations lower bound in the $ M_1$ topology with the rate function $I_y^{\BV[0,T]}$. 
	\item[(ii)]  (\cite{vysotsky}) Let $\phi$ be a real-valued function on $\D[0,T]$ which is uniformly continuous in the $M_1'$ topology on the level sets 
	$\{\xi: I_y^{\BV[0,T]}(\xi) \leq \alpha\}, \alpha <\infty$. Then $\phi(\bar K_x^0)$ satisfies an LDP with the rate function $cl(J_\phi)$, where $cl(J_\phi)$ is the lower-semicontinuous regularization of $J_\phi(u) = \inf_{\xi: \phi(\xi)=u}I_0^{\BV[0,T]}(\xi)$.
	\end{itemize}
\end{result}

\section{Computing \texorpdfstring{$\mathcal B_0^*$}{Lg}: finding a smooth minimizer}
%Proof of Proposition~\ref{bp-ison-bzero}, and Lemma~\ref{F-z-T-ccompact}}\label{reducetoconcave-N-compactenessFYZT}
\label{appendix-C}
In this appendix, we provide some details which could facilitate the computation of $\mathcal B_0^*$. Note that it is not straightforward that the infimum in the representation of $\mathcal{B}_0^*$ in \eqref{By} is attained, since the associated objective function does not have compact level sets, unless the moment generating function of $U_1$ is finite everywhere, cf.\ \cite{Lynch}. 
The following proposition, however, facilitates the characterization of  $\mathcal B_0^*$.
\begin{proposition}
\label{bp-ison-bzero} 
\linksinthm{bp-ison-bzero} 
Let $\notationdef{nota-B-y-AC}{B_y^\AC} \triangleq B_y \cap \AC[0,\infty)$, $\notationdef{nota-B-y-CNCV}{B_y^{\CON}} \triangleq B_y^\AC \cap \{\xi \in \AC[0,\infty): \xi \text{ is concave}\}$, and recall that $\mathcal B_y^* = \inf_{\xi \in B_y} I_{y}(\xi)$.
	Then,
	$$
	\mathcal{B}_y^{\ast} = \inf_{\xi \in B_y^{\AC}}I_{y}(\xi) = \inf_{\xi \in B_y^{\CON}}I_{y}(\xi).
	$$
\end{proposition}
We defer the proof of this proposition to the end of this appendix. %is proposition to Appendix~\ref{reducetoconcave-N-compactenessFYZT}.
We apply Proposition \ref{bp-ison-bzero} to write
\begin{equation}
\label{partitioned infimum}
    	\mathcal{B}_y^{\ast} = \inf_{\xi \in B_y^{\CON}}I_{y}(\xi) = \inf_{z_0, z_T: z_0 \geq z_T} \inf_{T\geq 0} \inf_{\xi \in F^{y}_{z_0, z_T, T}} I_{y}(\xi),
\end{equation}
where  $\notationdef{nota-F-z-T}{F^{y}_{z_0, z_T, T}}\triangleq\{ \xi: \xi \in B_y^{\CON}, \dot \xi(0)=z_0, \dot \xi(T)=z_T, \xi(T)=0 \}$. 
We next ensure compactness of $F^{y}_{z_0, z_T, T}$ in the following lemma, which is also proven at the end of this appendix. 
%of which we present the proof in Appendix~\ref{reducetoconcave-N-compactenessFYZT}.
%{\bf The variational problem $(\mathcal{B}_y)$ is arbitrarily defined and we are free to choose the topology. Later on in the document, in Lemma~\ref{B_y_star-equal-V_y_T_star}, we show that $B_y^*={V_y^*}^T$ . Only the variational problem $V_y^T$ depends on the topology used by Vysotsky through the contraction principle  and the SPLDP for trajectories of random walks w.r.t. the $M_1'$ topology. Using the $J_1$ topology the rate function $I_y$ is lower-semicontinuous on $F^{y}_{z_0, z_T, T}$: The function $I_y$  is lower-semicontinuous on $F^y{z_0,z_T,T}$ when $F^y{z_0,z_T,T}$ is endowed with the weak$^*$-topology (see \cite{Lynch}). Since the $J_1$ topology is stronger than the weak topology the function $I_y$ is lower-semicontinuous on $F^y_{z_0,z_T,T}$ w.r.t. the $J_1$ topology. The function $I_y$  is lower-semicontinuous on $F^y{z_0,z_T,T}$ when $F^y{z_0,z_T,T}$ is endowed with the weak$^*$-topology (see \cite{Lynch}). Since the $J_1$ topology is stronger than the weak topology the function $I_y$ is lower-semicontinuous on $F^y_{z_0,z_T,T}$ w.r.t. the $J_1$ topology.  } % I agree so I edited away this comment - B
\begin{lemma}
\label{F-z-T-ccompact}
\linksinthm{F-z-T-ccompact}
$F^{y}_{z_0, z_T, T}$ is a compact set with respect to the $J_1$ topology.
\end{lemma}
%{\bf Is this topology  appropriate? }
% Every element in the set $F_{z,T}$ can be written as $\xi(t) = \int_0^t \dot \xi(s) ds$ with $\dot{\xi}(s) \in [\mu, z]$. Using this, it is shown in Lemma~\ref{F-z-T-ccompact} that $F_{z,T}$ is compact.
 %, which is sometimes neglected in queueing applications of large deviations theory. %The phenomenon that the minimization over absolutely continuous functions may be strictly smaller than the minimum taken over Lipschitz continuous functions is
%known as the {\em Lavrentiev phenomenon}. We refer to \cite{Loewen} for an accessible overview of this topic.  
%Note that this proposition cannot be directly obtained by writing down Euler-Lagrange equations, because, again, it is not clear that the infimum is attained. The reason is that the associated rate functions are good only if the moment generating function of $U_1$ is finite everywhere, cf.\ \cite{Lynch}.
On the set $F^{y}_{z_0, z_T, T}$, the conditions $\xi(0)=y, \xi(T)=0$,  and concavity imply that $\xi(s)>0$ for $s\in (0,T)$. Therefore, the identity $I_{y}(\xi)= \int_0^T \Lambda^*(\dot \xi(s))ds$ holds. The RHS of this identity is lower-semi continuous in $\xi$ on the compact set $F^{y}_{z_0, z_T, T}$, and therefore attains a minimum $\xi^*$, as long as the set  $F^{y}_{z_0, z_T, T}$ is non-empty. The latter property holds if the solution $\xi$ with slope $z_0>0$ on $[0,T |z_T|/ (z_0+ |z_T|)]$ and slope $z_T<0$ on $[T |z_T|/ (z_0+ |z_T|, T]$ yields an area of at least 1. 

We now characterize the minimizer of the inner infimum of the RHS in Lemma \ref{partitioned infimum} through Euler-Lagrange equations. Such a characterization is usually only possible if the minimum is sufficiently smooth (e.g.\ not just AC, but $C^1$). This requires additional assumptions. We use Lagrange duality, and note that the feasible region of admissible paths is only convex (and hence, the absence of a duality gap is only guaranteed) if $p\leq 1$. In addition, 
we utilize sufficient conditions for smoothness of optimal solutions of variational problems developed in \cite{Cesari}, Chapter 2.6, which seems to exclude the case $p<1$, so in what follows we assume $p=1$.
We make the additional assumption that $\Lambda(\theta)= \log E[\exp \{\theta U \}]$ is steep at $\theta_+$ and $\theta_-$, i.e.\ $\lim_{\theta\uparrow \theta_+}
\nabla \Lambda (\theta) = \infty$ and  $\lim_{\theta\downarrow \theta_-}
\nabla \Lambda (\theta) = -\infty$. Under these assumptions, $\Lambda^*(z) = z \theta (z) - \Lambda (\theta (z))$, with $\theta(z)$ 
the unique solution of $z= \nabla \Lambda (\theta)$. The steepness assumptions make $\Lambda^*$ a smooth ($C^\infty$) function on the entire real line. Its derivative satisfies $\nabla \Lambda^*(z) = (\nabla \Lambda)^{-1}(z)$.
To reduce our setting to the framework in \cite{Cesari}, Chapter 2.6, we first incorporate our area constraint into a Lagrangian. Fix $\ell \geq 0$ and define 
\begin{equation}
    f_\ell(\xi(t), \dot \xi(t)) = \Lambda^*(\dot \xi (t)) - \ell [\xi(s) - 1/T].
\end{equation}
The Lagrangian $L_\ell (\xi)$ of our problem w.r.t.\ the constraint $\int_0^T \xi(s) ds \geq 1$ is $\int_0^T f_\ell (\xi (s), \dot \xi(s)) ds$.
We show that the problem of minimizing $L_\ell (\xi)$ over the set of concave absolutely continuous paths 
$\xi$ such that $\xi(0)=\xi(T)=0$, $\xi'(0)=z_0, \xi'(T)=z_T, z_0>0>z_T$ produces a solution which is $C^1$ for every $\ell \geq 0$.
Because $\nabla \Lambda^*$ is strictly increasing, $f_\ell$
satisfies (2.6.1) of \cite{Cesari} which demands that $f_\ell^{(y)}(x,y) = \frac d{dy} f_\ell^{(y)}(x,y)= \nabla\Lambda^*(y)$ is strictly increasing. Property (2.6.4) in \cite{Cesari}, which demands that $|f_\ell^{(y)}(x,y)|=|\nabla\Lambda^*(y)|\rightarrow\infty$ as $|y|\rightarrow\infty$, uniform in $x\geq 0$, follows from the fact that $\Lambda$ is steep.
We next propose an AC candidate solution $\xi^*$ which satisfies the Euler-Lagrange equation almost everywhere. 
In our setting, this equation is given by, for some constant $c$, 
% Applying the standard variational method of Lagrange multipliers (see, for example, \cite{luenberger1997optimization}), there exist constants $c$ and $\ell\geq 0$ such that $\xi^*$ satisfies the differential equation
\begin{equation}
	\nabla \Lambda^*(\dot{\xi^*}(s))=c- \ell s, s\in [0,T].
\end{equation}
	Since $\dot \xi(0) = z_0$, $c=\nabla \Lambda^*(z_0)$. Apply
	 %Because $\xi^*$ is concave and $\nabla \Lambda^*$ is non-decreasing, $\ell \geq 0$. Furthermore, s
	  $\nabla \Lambda^*(z) = (\nabla \Lambda)^{-1}(z)$ to write, for almost every $s$,
	  	\begin{equation}\label{eq:optimal path2}
	\dot{\xi^*}(s)=\nabla \Lambda \left( \nabla \Lambda^*(z) - \ell s\right), s\in [0,T].
	\end{equation}
Since $f_\ell$ is $C^1$, and $\dot \xi^*(t) \in [z_T, z_0]$ on $[0,T]$, 
we can now conclude from Theorem 2.6.ii in \cite{Cesari} that $\xi^*$ is $C^1$ on $[0,T]$, so that (\ref{eq:optimal path2}) is valid for for all $t \in [0,T]$. This expression can now be substituted into the Lagrangian $L_\ell(\xi)$. Maximizing this over $\ell$ gives an expression for the inner infimum in (\ref{partitioned infimum}), which can then be optimized further over
 $z_0, z_T, T$.\\

\begin{proof}{Proof of Proposition~\ref{bp-ison-bzero}.}
\linksinpf{bp-ison-bzero}
Due to Lemma~\ref{useful-properties-of-variational-By} \emph{(ii)}, $\mathcal{B}_y^*= \inf_{\xi \in B_y^{\mathbbm{AC}}} I_y(\xi)$. 
Since $B_y^\CON \subseteq B_y^\AC$, we only have to prove that $\inf_{\xi \in B_y^{\mathbbm{AC}}} I_y(\xi) \geq \inf_{\xi \in B_y^\CON} I_y(\xi)$. 
For this, we show that for any given $\xi \in B_y^{\AC}$, there is $\zeta \in B_y^\CON$ such that $I_y(\zeta) \leq I_y(\xi)$.
To construct such $\zeta$, we first note that we can find $\xi_1 \in B_y^{\AC}$ such that $\T(\xi_1)<\infty$ and $I_y(\xi_1) \leq I_y(\xi)$ thanks to Lemma~\ref{we-can-consider-only-upto-M}.
Now set $T = \T(\xi_1)$ and denote the restriction of $\xi_1$ on $[0,T]$ with $\check\xi_1$---i.e., $\check\xi_1\in \AC[0,T]$ and $\check\xi_1(t) = \xi_1(t)$ for $t\in[0,T]$.
We appeal to \emph{(iii)} of Lemma~\ref{merging-actions-to-the-left-lemma} to pick a path  $\xi_2 \in \AC[0,T]$ such that $\xi_2(0) = y$, $\Phi_T(\xi_2)  \geq \Phi_T(\check\xi_1) \geq 1$, $I_{y}^\BVT(\xi_2) \leq I_y^\BVT(\check\xi_1)  = I_y(\xi_1) \leq I_y(\xi)$, and $\xi_2$ concave on $[0,T]$ with the derivative bounded by $\mu$ from below. Now, set $\zeta=\xi_2(t \wedge T)+\mu([t-T]^+), \ t \geq 0$.
Then, $\zeta \in B_y^\CON$ and $I_y(\zeta)=I_y^{\BV[0,T]}(\xi_2) \leq I_y(\xi)$.
\eop
\end{proof}

We end this section with the proof of compactness of $F^y_{z_0,z_T,T}$.

\begin{proof}{Proof of Lemma~\ref{F-z-T-ccompact}.}
\linksinpf{F-z-T-ccompact}

Recall that $F^{y}_{z_0, z_T, T}=\{ \xi: \xi \in B_y^{\CON}, \dot \xi(0)=z_0, \dot \xi(T)=z_T, \xi(T)=0 \}$.
Let $\cal P([0,T])$ denote the space of all Borel probability measures.
Consider $\varphi:\cal P([0,T])\to\D[0,T]$, a mapping defined for each $\nu \in \cal P([0,T])$ as follows:
$$\varphi(\nu)(t) \triangleq y + \int_0^t\big(z_0+(z_T-z_0)\nu([0,s])\big)ds.$$ 
We consider the weak topology on $\cal P([0,T])$ and the $J_1$ topology on $\D[0,T]$. 
Note that $\cal P([0,T])$ is a compact space due to Prokhorov's theorem.
Now, observe that 
$$F^{y}_{z_0, z_T, T} = \varphi(\cal P([0,T]))\cap \{\xi \in \D[0,T]: \xi(T) = 0\}.$$
Since $\{\xi \in \D[0,T]: \xi(T) = 0\}$ is closed, the proof of this lemma is complete if $\varphi$ is continuous.
To confirm the continuity of $\varphi$, consider a sequence $\nu_n$ that converges to $\nu$ in $\cal P([0,T])$. 
Note that 
\begin{align}
d_{J_1}(\varphi(\nu_n), \varphi(\nu)) 
&\leq \|\varphi(\nu_n)- \varphi(\nu)\|_\infty 
= \sup_{t\in[0,T]}\left| (z_T-z_0)\int_0^t\big( \nu_n([0,s]) - \nu([0,s])\big) ds \right|
\nonumber\\
&\leq |z_0 - z_T| \int_0^T \Big| \nu_n([0,s]) - \nu([0,s])\Big|ds,\label{bound-of-J1-distance-between-varphi-nu-s}
\end{align}
and $\big|\nu_n([0,s])- \nu([0,s])\big|$ is bounded by 2 for each $n$. 
Since the (weak) convergence of $\nu_n$ to $\nu$ implies that $\big|\nu_n([0,s]) - \nu([0,s])\big|\to0$ for almost every $s\in[0,T]$, the bounded convergence theorem guarantees that \eqref{bound-of-J1-distance-between-varphi-nu-s} converges to zero as $n\to\infty$, proving the continuity of $\varphi$.
This concludes the proof of the desired compactness of $F^{y}_{z_0, z_T, T}$. 
\eop
\end{proof}

\section*{Acknowledgements}
The research of
BZ and MB is supported by NWO grant 639.033.413. The research of JB is supported by NSF
grants 1915967, 1820942, 1838576 as well as DARPA award N660011824028. The research of C-HR is supported by NSF grant CMMI-2146530. The authors are grateful to Oliver Tse and Vladislav Vysotsky for useful discussions and to the referees for comments that helped to improve the exposition of the paper.

%{\bf The references need to be cleaned up: using et al is not appropriate, and some of the journal names are abbreviated, while others not. This needs to be made consistent. When adding a new bibtex item, please do not make it by hand, but go to the website where the paper officially appeared and export the bibtex item. }

\end{APPENDICES}

\bibliographystyle{informs2014}
\bibliography{BBRZDVE}

\ifnotationindex
\newpage
\linkdest{location of notation index}
\noindent{\large \bf Notations}

\begin{itemize}[leftmargin=*]

\linkdest{\notationindexlocationphrase A}
    \item \notationidx{nota-A-t-epsilon}{$A_{t,\epsilon}$}: $A_{t,\epsilon}=\{\xi \in \D[0,t]: \xi(0)=\epsilon, \  \int_{0}^{t}\Psi(\xi)(s)^p ds > 1, \, \xi(s) >0,\ \forall \, s \, \in [0,t]\}$
	\item \notationidx{nota-tilde-A-t-epsilon}{$\tilde A_{t,\epsilon}$}: $\tilde A_{t,\epsilon}=\{\xi \in \D[0,t]: \xi(0)=\epsilon, \ \int_{0}^{t}\Psi(\xi)(s)^p ds > 1, \, \xi(s) >\epsilon/2, \, \forall s \, \in [0,t]\}$
    \item \notationidx{nota-alpha}{$\alpha$}: $\alpha \triangleq 1/(1+p)$

\linkdest{\notationindexlocationphrase B}
    \item \notationidx{nota-B-y-AC}{$B_y^\AC$}: $B_y^\AC \triangleq B_y \cap \AC[0,\infty)$
    \item \notationidx{nota-B-y-CNCV}{$B_y^\CON$}: $B_y^{\CON} \triangleq B_y^\AC \cap \{\xi \in \AC[0,\infty): \xi \text{ is concave}\}$
    \item \notationidx{nota-B-y-AC-M}{$B^{\AC;M}_{y}$}: $B^{\AC;M}_{y} \triangleq B^{\AC}_y \cap \{\xi \in \D[0,\infty): \T(\xi) \leq M\} $
    \item \notationidx{nota-B-y-M}{$B_y^M$}: $B_y^M \triangleq B_y \cap \{\xi \in \D[0,\infty): \T(\xi) \leq M\}$

    \item \notationidx{nota-By}{$B_y$}: $B_y \triangleq \left\{ \xi \in \BV[0,\infty): \xi(0)=y, \ \int_{0}^{\T(\xi)}\Psi(\xi)(s)^pds \geq 1 \right\}$
    \item \notationidx{nota-cal-B-y-star}{$\mathcal{B}_{y}^{\ast}$}: $\B_y^* \triangleq \inf_{\xi \in B_{y}} I_y(\xi)$
    \item \notationidx{nota-cal-B-y}{$(\mathcal B_y)$}: optimization problem associated with $\mathcal{B}_{y}^{\ast}$
    \item \notationidx{nota-cal-B-pi-star}{$\mathcal{B}_{\pi}^{\ast}$}: $\B_\pi^* \triangleq \inf_{y \in [0,\infty),\, \xi \in B_{y}} \left\{ \beta y+ I_y(\xi) \right\}$
    \item \notationidx{nota-cal-B-pi}{$(\mathcal B_\pi)$}: optimization problem associated with $\mathcal{B}_{\pi}^{\ast}$
    
    \item \notationidx{nota-BV-T}{$\BV[0,T]$}:  subspace of $\D[0,T]$ consisting of paths that are of bounded variation
    \item \notationidx{nota-BV-infty}{$\BV[0,\infty)$}:  subspace of $\D[0,\infty)$ consisting of paths that are of bounded variation on any compact interval
    
    \item \notationidx{nota-beta}{$\beta$}: $\beta \triangleq \sup\left\{\theta \geq 0: \E(e^{\theta  U}) \leq 1 \right\}$

\linkdest{\notationindexlocationphrase C}
\linkdest{\notationindexlocationphrase D}
    \item \notationidx{nota-D}{$\D[0,T]$}: space of c\'adl\'ag paths from $[0,T]$ to $\R$
    \item \notationidx{nota-D-leqslant-infty-uparrow}{$\D_{\leqslant \infty}^{\uparrow}[0,T]$}: subspace of $\D[0,1]$ consisting of non-decreasing pure jump paths that assume non-negative values at the origin
    \item \notationidx{nota-D-lambda}{$\D^{(\lambda)}[0,T]$}: $\D^{(\lambda)}[0,T]= \{\xi \in \D[0,T]: \xi(t) =\lambda t + \zeta(t),\,\forall t\in[0,T],\, \zeta \in \D_{\leqslant \infty}^{\uparrow}[0,T]\}.$
    \item \notationidx{nota-D-Lambda}{$D_\Lambda$}: $D_{\Lambda}=\{x: \Lambda(x)< \infty\}$
    \item \notationidx{nota-dM1'}{$d_{M_1'}$}: $d_{M_1'}(\xi,\zeta)
		\triangleq
		\inf_{\substack{(u,t) \in \Gamma'(\xi)\\(v,r) \in \Gamma'(\zeta)}}
		\{
		\|u-v\|_\infty + \|t-r\|_\infty
		\}$

\linkdest{\notationindexlocationphrase E}
\linkdest{\notationindexlocationphrase F}
    \item \notationidx{nota-F-z-T}{$F_{z,T}$}: $F_{z,T}=\{ \xi: \xi \in B_0^{\CON}, \dot \xi(0)=z, \xi(T)=0 \}$

    % \item \notationidx{nota-fx}{$f(x)$}: $f(x)=x^p$

\linkdest{\notationindexlocationphrase G}
    \item \notationidx{nota-Gamma'}{$\Gamma'(\xi)$}: $\Gamma'(\xi) \triangleq \{(u,t)\in \R\times [0,1]: u\in [\xi(t-)\wedge \xi(t),\ \xi(t-)\vee \xi(t)]\}$ \hfill (extended completed graph)

\linkdest{\notationindexlocationphrase H}
    \item \notationidx{nota-H}{$H: \D[0,T] \rightarrow [0,\infty)$}: $H(\xi) = \int_0^T \xi (s)ds$
\\[-10pt]

\linkdest{\notationindexlocationphrase I}
    \item \notationidx{nota-I-alpha}{$I_\alpha(\xi)$}: 
        $
        I_\alpha(\xi) =
        \begin{cases}
        \sum_{t: \xi(t) \neq \xi(t-)}(\xi(t) - \xi(t-))^{\alpha} & \mathrm{for} \ \xi \in \D^{(\lambda)}[0,T], \\
        \infty & \ \text{otherwise}. 
        \end{cases}
	$ \CR{local}
	\\[0pt]
   \item \notationidx{nota-I-S}{$I_S(\zeta)$}: 
        $
    		I_S(\zeta) \triangleq
    		\begin{cases}
    		{\mathcal{B}^*_{0}}(\zeta(1)-\zeta(1-))^{\alpha}  &  \mathrm{if} \ \zeta \in \D^{\leqslant 1}[0,T], \\
    		\infty & \ \text{otherwise}.
    		\end{cases}
	$
	\\[0pt]
    \item \notationidx{nota-I-t}{$I_{\mathbf t}(x_1,\ldots,x_k)$}: 
        $
	I_{\mathbf t}(x_1,\ldots,x_k)=
	\begin{cases}
	{\mathcal{B}_0^*}\sum_{i=1}^k (x_i-\lambda\Delta t_i)^{\alpha} & \text{if} \ x_i\geq \lambda \Delta t_i,\ \forall i=1,\ldots, k,\\
	\infty & \ \text{otherwise}.  
	\end{cases}
	$ \CR{local}
	\\[0pt]
    \item \notationidx{nota-I-tilde-t}{$\tilde I_{\mathbf t}(x_1,\ldots,x_k)$}: 
        $
	\tilde{I}_{\mathbf t}(x_1,\ldots,x_k)=
	\begin{cases}
	{\mathcal{B}_0^*}\sum_{i=1}^k (x_i-x_{i-1}-\lambda \Delta t_i)^{\alpha} & \mathrm{if} \ x_i-x_{i-1} \geq \lambda \Delta t_i, \mathrm{\ for\ } i=1,\ldots, k \\
	\infty & \ \text{otherwise}.  
	\end{cases}
	$ \CR{local}
	\\[-10pt]
    \item \notationidx{nota-I-V}{$I_V(x)$}: $I_V(x) = \mathcal B_0^*\cdot x^\alpha$ \CR{local}\\[-10pt]

    \item \notationidx{nota-I-W}{$I_{\mathbf{W}}(\zeta)$}: $I_{\mathbf{W}}(\zeta) \triangleq \inf\left\{ I_{Z,S}(\xi_{1},\xi_{2}): \zeta =\xi_1+ \xi_2,  \ \xi_1 \in \D^{(\lambda)}[0,T], \, \xi_2 \in \D^{\leqslant 1}[0,T] \right\}$\\[-10pt]
    \item \notationidx{nota-I-Y}{$I_{\textbf{Y}}(\zeta)$}: 
            	$
            	I_{\textbf{Y}}(\zeta)
            	\triangleq
            	\begin{cases}
            	{\mathcal{B}^*_{0}} \sum_{t: \zeta(t) \neq \zeta(t-)}\left(\zeta(t)-\zeta(t-)\right)^{\alpha}
            	&
            	\text{if } \zeta \in  \D^{(\lambda)}[0,T], 
            	\\
            	\infty
            	&
            	\text{otherwise}.
            	\end{cases}
            	$\\[0pt]

    \item \notationidx{nota-Iy}{$I_y(\xi)$}: 
            $
            I_{y} (\xi)
            \triangleq
            \begin{cases}
            \int_{0}^{ \T(\xi)}  \Lambda^*(\dot{\xi}^{(a)}(s))ds + \theta_+\xi^{(u)}(\T(\xi)) + \theta_-\xi^{(d)}(\T(\xi))
            & \text{if } \xi(0) = y \text{ and }\xi \in \BV[0,\infty)
            \\
            \infty
            & \text{otherwise}
            \end{cases}
            $
            \\[0pt]
    \item \notationidx{nota-I-y-BV}{$I_y^\BVT(\xi)$}: $I_y^\BVT(\xi)\triangleq \begin{cases}
            \int_{0}^{T} \Lambda^{*}(\dot{\xi}^{(a)}(s)) d s+\theta_{+} \xi^{(u)}(T)+\theta_{-} \xi^{(d)}(T)
            &
            \text{if } \xi(0) = y \text{ and } \xi \in \BV[0,T] ,
            \\
            \infty
            &
            \text{otherwise.}
            \end{cases}$\\[0pt]
    \item \notationidx{nota-I-Z}{$I_Z(\xi)$}: 
        $
        I_Z(\xi)\triangleq
	\begin{cases}
	{\mathcal{B}_0^*}\sum_{t: \xi(t) \neq \xi(t-)}(\xi(t) - \xi(t-))^{\alpha} & \mathrm{if} \ \xi \in \D^{(\lambda)}[0,T], \\
	\infty & \ \text{otherwise}.
	\end{cases}
	$ 
	\\[-8pt]
    \item \notationidx{nota-I-ZS}{$I_{Z,S}(\zeta,\xi)$}: $I_{Z,S}(\zeta,\xi)=I_{Z}(\zeta)+I_{S}(\xi)$

\linkdest{\notationindexlocationphrase J}
    \item \notationidx{nota-J-+}{$J_+$}: $J_+ \triangleq \{\xi \in \D[0,\infty): \xi^{(d)}\equiv0\}$

    \item \notationidx{nota-J-y}{$J_{y}(a)$}: $J_{y}(a)= \inf\{I_0^{\BV[0,T]}(\xi): \xi \in \D[0,T], \xi(0)=y, \Phi_T(\xi)=a\}$.

\linkdest{\notationindexlocationphrase K}
    \item \notationidx{nota-K-n}{$K_n$}: $K_n =\sum_{i=1}^n U_i$
    \item \notationidx{nota-bar-K-n-y}{$\bar K_n^y$}: $\bar K_n^y (t)= y + \frac{1}{n}K_{\lfloor nt \rfloor}$, \quad where \quad $K_n =\sum_{i=1}^n U_i$  
    \item \notationidx{nota-K-t}{$K_t$}: $K_t\triangleq \left\{\xi \in \D[0, t]: \xi(0)=0, \ \int_{0}^{t}(\Psi(\xi)(s))^p ds \geq 1, \ \xi(s) \geq 0 \ \text{for} \ s \in [0, t] \right\}$

\linkdest{\notationindexlocationphrase L}
    \item \notationidx{nota-Lambda-star}{$\Lambda^*(y)$}: $\Lambda^*(y)\triangleq \sup_{\theta \in \mathbbm{R}}\{\theta y -\log \E(e^{\theta U})\}.$

    \item \notationidx{nota-lambda}{$\lambda$}: $\lambda= \E(\sum_{i=1}^{T_1}X_i^p)/\E(T_1)$

\linkdest{\notationindexlocationphrase M}
    \item \notationidx{nota-mu}{$\mu$}:  $\mu= \E(U_1)<0$

\linkdest{\notationindexlocationphrase N}
    \item \notationidx{nota-N-t}{$N(t)$}: $N(t) =\max\{k \geq 0: T_k \leq t\}$

\linkdest{\notationindexlocationphrase O}
\linkdest{\notationindexlocationphrase P}
    \item \notationidx{nota-Phi}{$\Phi(\xi)$}: $\Phi(\xi)\triangleq \int_{0}^{\mathcal{T}(\xi)}(\Psi(\xi)(s))^p ds$

    \item \notationidx{nota-Phi-T}{$\Phi_T$}: $\Phi_T(\xi)=\int_{0}^{T}(\Psi(\xi)(s))^pds$

    \item \notationidx{nota-Pi'}{$\Pi'(\xi)$}: Set of all parametrizations of $\xi$

    \item \notationidx{nota-Psi}{$\Psi(\xi)$}: $\Psi(\xi)(t)=\xi(t)-\inf_{0 \leq s \leq t}\{\xi(s) \wedge 0\}, \ \forall  t \geq 0$

\linkdest{\notationindexlocationphrase Q}
\linkdest{\notationindexlocationphrase R}
    \item \notationidx{nota-bar-cal-R-n}{$\bar R_n$}: $\bar{R}_n(t) = \frac{1}{n}\sum_{i=T_{N(n)}+1}^{\lfloor nt \rfloor}X_i^p$

\linkdest{\notationindexlocationphrase S}
    \item \notationidx{nota-bar-S-n}{$\bar S_n$}: $\bar{S}_n(t) = \frac1n V_n\mathbbm{1}_{\{1\}}(t)$

\linkdest{\notationindexlocationphrase T}
    \item \notationidx{nota-theta-pm}{$\theta_+, \theta_-$}: $\theta_+=\sup \{\theta: \E(e^{\theta U_1})< \infty \}$, $\theta_-=\inf \{\theta: \E(e^{\theta U_1})< \infty \}$

    \item \notationidx{nota-T-i}{$T_i$}: $T_i \triangleq \inf\{ k> T_{i-1}: X_k = 0\}$ for $i\geq 1$, \quad $T_0 \triangleq 0$

    \item \notationidx{nota-T}{$T(x)$}: $T(s)\triangleq x \cdot \I_{\{1\}}$ \CR{local} 

    \item \notationidx{nota-cal-T}{$\T(\xi)$}: $\T(\xi) = \inf\{t>0: \Psi(\xi)(t) \leq 0\}$

    \item \notationidx{nota-TM1'}{$\mathcal{T}_{M_1'}$}: $M_1'$ Skorokhod topology

    \item \notationidx{nota-TV}{$\mathrm{TV}(\xi)$}: Total variation of $\xi$

    \item \notationidx{nota-tau-j}{$\tau_j$}: $\tau_j = T_{j}-T_{j-1}$, the inter-arrival times of the renewal process $N$

\linkdest{\notationindexlocationphrase U}
    \item \notationidx{nota-Ui}{$\{U_i, i \geq 1\}$}: i.i.d.\ random variables s.t.\ $\E(U_1)=\mu <0$

\linkdest{\notationindexlocationphrase V}
    \item \notationidx{nota-V-n}{$V_n$}: $V_n = \sum_{i=T_{N(n)}+1}^{n}X_i^p$
    \item \notationidx{nota-cal-V-y-T-*}{${\mathcal V_y^T}^*$}: ${\mathcal V_y^T}^*\triangleq  \inf_{\xi \in V_{y}^T}I_{y}^\BVT(\xi)$
    \item \notationidx{nota-cal-V-y-T}{$(\mathcal{V}_y^{T})$}: optimization problem associated with ${\mathcal V_y^T}^*$
    \item \notationidx{nota-V-y-T}{$V_y^T$}: $V_y^T \triangleq \left\{ \xi \in \D[0,T]: \xi(0)=y,   \ \Phi_T(\xi) \geq 1 \right\}$

\linkdest{\notationindexlocationphrase W}
    \item \notationidx{nota-cal-W}{$\mathcal{W}$}: pointwise convergence

    \item \notationidx{nota-W-j}{$W_j$}: $   W_{j} = \sum_{i=T_{j-1}+1}^{T_j}X_i^p$

\linkdest{\notationindexlocationphrase X}
    \item \notationidx{nota-Xk}{$X_n$}: $X_{n+1} \triangleq [X_n+U_{n+1}]^{+}, \, n\geq 0$, \quad $X_0=0$
    
    \item \notationidx{nota-xia}{$\xi^{(a)}$}: absolutely continuous part of $\xi$
    \item \notationidx{nota-xis}{$\xi^{(s)}$}: sigular part of $\xi$; $\xi^{(s)}(0) = 0$
    \item \notationidx{nota-xiu}{$\xi^{(u)}$}: non-decreasing singular part of $\xi$; $\xi^{(u)}(0) = 0$
    \item \notationidx{nota-xid}{$\xi^{(d)}$}: non-increasing singular part of $\xi$; $\xi^{(d)}(0) = 0$

\linkdest{\notationindexlocationphrase Y}
    \item \notationidx{nota-bar-y}{$\bar{y}$}: $\bar{y} = (|\mu|(p+1))^{1/1+p}$
    
    \item \notationidx{nota-Y-bar}{$\bar Y_n$}: $\bar Y_n(\cdot) \triangleq \frac{1}{n}\sum_{i=1}^{\lfloor n \cdot \rfloor}X_i^p$

\linkdest{\notationindexlocationphrase Z}
    \item \notationidx{nota-bar-Z-n}{$\bar Z_n$}: $\bar{Z}_n(t) = \frac 1 n\sum_{j=1}^{N(nt)}W_j$

\end{itemize}
\fi

%%%%%%%%%%%%%%%%%%%%%%%%%%%%%%%%%%%%%%%%%%%%%%%%%%%%%%%%%%%%%
%%%%%%%%%%%%%%%%%%%% \ifnavigationlinks %%%%%%%%%%%%%%%%%%%%%
%%%%%%%%%%%%% end of \ifnavigationlinks %%%%%%%%%%%%%%%%%%%%%

\end{document}